\documentclass[11pt]{article}

\usepackage{amsmath, amssymb, bbm} 
\usepackage{verbatim}
\usepackage{hyperref}
\usepackage{color}
\usepackage{graphicx}
\usepackage{amsthm}
\usepackage{enumerate}
\usepackage{dsfont,bm}
\usepackage{esint} 

\newtheorem{theorem}{Theorem}[section]
\newtheorem{thm}[theorem]{Theorem}
\newtheorem{corollary}[theorem]{Corollary}
\newtheorem{cor}[theorem]{Corollary}
\newtheorem{lemma}[theorem]{Lemma}
\newtheorem{lem}[theorem]{Lemma}
\newtheorem{proposition}[theorem]{Proposition}
\newtheorem{prop}[theorem]{Proposition}
\newtheorem{definition}[theorem]{Definition}
\newtheorem{defn}[theorem]{Definition}
\newtheorem{assumption}[theorem]{Assumption}
\newtheorem{remark}[theorem]{Remark}
\newtheorem{example}[theorem]{Example}
\numberwithin{equation}{section}

\def\be{\begin{equation}}
\def\ee{\end{equation}}
\def\bes{\begin{equation*}}
\def\ees{\end{equation*}}

\def\<{\langle}
\def\>{\rangle}

\newcommand\restr[2]{{
		\left.\kern-\nulldelimiterspace 
		#1 
		\vphantom{\big|} 
		\right|_{#2} 
}} 
\newcommand{\abs}[1]{{\left\vert\kern-0.25ex #1
		\kern-0.25ex\right\vert}}
\newcommand{\on}[1]{\operatorname{ #1}}
\newcommand{\set}[1]{\left\{ #1 \right\}}

\newcommand{\Sett}[2]{\left\{ #1  : \, #2 \right\}}
\newcommand{\diag}[0]{\operatorname{diag}}

\newcommand{\BC}[0]{\hyperlink{bc}{$\operatorname{(MD)}$}}
\newcommand{\one}{\mathds{1}} 
\DeclareMathOperator*{\esssup}{ess\,sup}
\DeclareMathOperator*{\essinf}{ess\,inf}

\def\pf{{\medskip\noindent {\bf Proof. }}}

\def\sA {{\mathcal A}} \def\sB {{\mathcal B}} 
 \def\sE {{\mathcal E}} \def\sF {{\mathcal F}}

 \def\sN {{\mathcal N}}

  \def\sX {{\mathcal X}}

 \def\bE {{\mathbb E}} 
\def\bG {{\mathbb G}}  
  
 \def\bN {{\mathbb N}} 
\def\bP {{\mathbb P}} \def\bQ {{\mathbb Q}} \def\bR {{\mathbb R}}
  
\def\bV {{\mathbb V}}  
 \def\bZ {{\mathbb Z}}

\def\med{\medbreak\noindent}

\def\sms{\smallskip}
\def\ms{\medskip}

\def\sm{\smallskip\noindent}

\def\ignore#1{}

\def\ol{\overline}           

\def\eps{\varepsilon}
 
\def\vp{\varphi}
\def\Gam{\Gamma} \def\gam{\gamma}

\def\to {\rightarrow}

\def\pd {\partial}
\def\q{\quad} \def\qq{\qquad}
\def\dint{\int\kern-.6em\int}

\def\diam{\mathop{{\rm diam }}}

\def\div{{\mathop {{\rm div\, }}}}

\def \half {{\tfrac12}}
\def\fract{\tfrac}

\def\wt{\widetilde}

\def\be{\begin{equation}}
\def\ee{\end{equation}}
\def\bes{\begin{equation*}}
\def\ees{\end{equation*}}
\def\ba{\begin{align}}
\def\ea{\end{align}}
\def\xxea{\end{align}}
\def\bas{\begin{align*}}
\def\eas{\end{align*}}

\def\proof{{\smallskip\noindent {\em Proof. }}}
\def\qed{{\hfill $\square$ \bigskip}}

\definecolor{dgreen}{rgb}{0, 0.6, 0.1}
\definecolor{dblue}{rgb}{0, 0.0, 0.6}
\definecolor{vdblue}{rgb}{0,.08, 0.45}
\definecolor{dred}{rgb}{0.7, 0.0, 0.0}
\definecolor{vdblue}{rgb}{0,.08, 0.45}

\definecolor{purple}{rgb}{0.6, 0.0, 0.6}
\definecolor{mytext}{rgb}{0.1, 0.1, 0.1}

\def\red{\color{red}}

\begin{document}
	
	\font\titlefont=cmbx14 scaled\magstep1
	\title{\titlefont  \vspace{-5ex}  Stability of EHI and regularity of MMD spaces  }
	
	\author{
 	Martin T. Barlow, \ Zhen-Qing Chen \ and \ Mathav Murugan
	}
	
	\maketitle

 \begin{abstract}
 	We show that elliptic Harnack inequality is stable under  	form-bounded 
	perturbations   for strongly local symmetric Dirichlet forms on complete locally compact separable 
	metric spaces that satisfy metric doubling property (or equivalently, relative   ball connectedness property).

\vskip.2cm

	
\vskip.2cm


\end{abstract}
	
	\section{Introduction} \label{sec:intro}

Let $X=\{ X_t, t \in [0,\infty); \bP^x, x \in \sX\}$ be a diffusion process on 
a locally compact separable
metric space $(\sX,d)$. 
A function $h$ on a ball $B=B(x,r)$ is {\em harmonic} if $h(X_{t \wedge \tau_B})$ is a local martingale
 under $\bP^x$ for every $x\in B$; here $\tau_B$ is the exit time from $B$ by the process $X$ 
and the filtration is the minimal augmented filtration generated by $X$. 
The (scale-invariant) {\em elliptic Harnack inequality} (EHI) holds for $X$ if there exist constants 
 $C\geq 1$ and $\delta \in (0, 1)$ 
such that whenever $h$ is non-negative and harmonic on a ball
 $B=B(x,r)$, then  
\be
 \sup_{B(x, \delta r)} h \le C \inf_{B(x, \delta r)} h. 
 \ee
If it holds, the EHI is a valuable tool for the study of the process $X$ and its   associated heat kernel. 
A well-known theorem of Moser \cite{Mo1} is that the EHI holds if $X$ is the symmetric diffusion associated 
with a uniformly elliptic divergence form operator  
 $\sA = {\rm div} (A(x) \nabla)$.
 Associated with such a process is the symmetric Dirichlet form $(\sE, \sF)$ on $L^2(\bR^d; dx)$,
where 
$$
\sF=W^{1,2}(\bR^d)= \left\{ f\in L^2(\bR^d; dx): \nabla f \in L^2(\bR^d; dx) \right\}
$$ 
 is the Sobolev space on $\bR^d$ of order $(1, 2)$ and 
$$ 
\sE(f,f) = \int_{\bR^d}  \nabla f (x)\cdot A(x) \nabla f (x) \, dx, \quad f\in \sF.
$$ 
 We say two symmetric 
Dirichlet forms  $\sE^{(1)}$ and $\sE^{(2)}$ 
on $L^2(\bR^d; dx)$ with common domain $\sF$ are {\em comparable} if
there exists $C \ge 1$ such that
$$ 
C^{-1}  \sE^{(1)} (f,f) \le  \sE^{(2)} (f,f) \le C  \sE^{(1)} (f,f)  \quad \mbox{ for all  } f \in \sF.
$$
 Moser's result gives the stability of the EHI, in the sense that if $\sE^{(1)}$ and $\sE^{(2)}$ 
are comparable symmetric Dirichlet forms on $L^2(\bR^d; dx)$, associated with uniformly elliptic divergence form operators
$\sA_i$, then the EHI holds for $\sE^{(2)}$ 
if and only if it holds for $\sE^{(1)}$.

A few years later, Moser \cite{Mo2, Mo3} proved a parabolic Harnack inequality PHI, which 
holds for non-negative solutions to the heat equation associated with a  uniformly elliptic divergence form operator $\sA$. 
In particular, if $u$ is any non-negative solution to the heat equation $\frac{\partial u}{\partial t}= \sA u$ in a time-space cylinder $Q=(0,r^2) \times B(x,r)$, then writing 
 $T=r^2, 
Q_-= (T/4,T/2)\times B(x,\delta r), Q_+=(3T/4,T)\times B(x,\delta r)$, we have
\[
\esssup_{Q_-}u \le C_P \essinf_{Q_+} u,
\]
where the constants $C_P >1$ and $\delta\in(0,1)$ do not depend on $x,r$ or $u$.
Subsequently Grigor'yan and Saloff-Coste in \cite{Gr0, Sal92} gave a characterization of the PHI,
and the stability of the PHI follows immediately from this characterization. 
The methods of these papers 
 are very robust, and this characterization of the PHI 
was extended to diffusions on locally compact separable metric  spaces \cite{St}, and to random walks on graphs \cite{De1}.  

For a number of years the stability of the apparently simpler EHI remained an open
problem.  Stability on a large class of 
  unbounded spaces 
(including  Riemannian manifolds and graphs) was proved by two of us recently in \cite{BM1}. However,
the result there relied on the metric space being geodesic and satisfying 
 some strong local regularity conditions;  
one key use of this regularity was to ensure the existence of Green's functions. 

The natural context for the study of the EHI is that of 
 locally compact separable metric measure  spaces with strongly local regular (symmetric) Dirichlet forms,
which we call   MMD spaces.  Examples include 
Riemannian manifolds, the cable systems of graphs \cite{V}, as well as various
classes of fractals. Not only do MMD spaces provide a common framework for all these examples, 
but also certain transformations (change of measure, quasisymmetric change
of metric) 
 which are natural for MMD spaces but
 are not so natural for manifolds and graphs.
   These transformations are key to the argument in \cite{BM1}. 

\sms
This paper has three main goals: 
 
  \begin{description}
\item{(i)}  We give a  weak sufficient condition (a local Harnack inequality) 
for a MMD space to have Green's functions.
This improves significantly the results of earlier papers, such as \cite{BM1, BM2}, which
needed some parabolic regularity. 
 In particular, it allows us to drop the Green function assumption (\cite[Assumption 2.3]{BM1}) made in \cite{BM1}. 
  
\item{(ii)} We carry through the program of \cite{BM1} in the context of a MMD space satisfying
these weak regularity conditions. 
 In particular, we drop  the {\it bounded geometry} assumption (see \cite[Assumption 2.5]{BM1} for its definition)
on the MMD space $(\sX, d, m, \sE, \sF)$,
   and relax the condition
that $(\sX,d)$ is a length (or geodesic) space;  both are needed  in \cite{BM1}. 
 We assume that $(\sX, d)$ is a complete metric space and satisfies metric 
doubling---see Definition \ref{D:MD}, which is equivalent
under EHI to the  `relatively ball connectedness' condition (see Definition \ref{d:metric}).  
The latter
  has the advantage that  it is preserved by quasisymmetric changes of metric.
 Example \ref{E:instab1} shows that some regularity  of the metric is needed 
if we are to have stability of the EHI.  

\item{(iii)}  
The  metric space  $(\sX, d)$ may either be  of bounded diameter or of  infinite diameter. 
\end{description}

For a metric space $(\sX,d)$, we use $B(x, r)$ to denote the open ball centered at $x\in \sX$ with radius $r$.  The closure and the boundary of the ball $B(x, r)$ will be denoted as 
$\overline{B(x, r)}$ and $\partial {B(x, r)} $, respectively.   
 
 \begin{definition}   \label{D:MD}  	
 A metric space $(\sX,d)$ is said to be  {\em metric doubling} (MD) if there exists $N\geq 2 $ such that for any $x \in \sX$,
$R>0$ there exist $z_1 , \dots, z_N\in \sX$ such that $B(x,R) \subset \cup_{i=1}^N B(z_i, R/2)$. 
We call a metric space that satisfies the metric doubling property a \emph{doubling metric space}.
	 \end{definition}

 \begin{remark} \label{R:1.2}   \rm 
 	Assouad \cite{As} showed that if $(\sX, d)$ is (MD), then for every $\alpha \in (0, 1)$, 
the metric space  $(\sX, d^\alpha)$  has a bi-Lipschitz embedding into ${\mathbb R}^n$ for some $n\geq 1$. 
This in particular implies that a  doubling metric space $(\sX, d)$ is always separable.
One can also deduce the separability of a  doubling metric space from its definition.
Indeed, every ball in a  doubling metric space is totally bounded by
 \cite[Exercise 10.17]{Hei}, which implies the separability.
 Furthermore, if $(\sX, d)$ is complete and (MD),  then  every ball in $(\sX, d)$ is relatively compact by the aforementioned totally boundedness property. Consequently, 
 any complete doubling metric  space $(\sX, d)$  is    separable and locally compact. 
 \end{remark}

 \medskip
 
 When $(\sX, d) $ is a locally compact metric space, a Dirichlet form $(\sE, \sF)$ on $L^2(\sX, m)$ 	
is said to be {\it strongly local}  if $\sE(u, v)=0$ whenever $u, v\in \sF$ have compact supports with $v$ being constant in
an open neighborhood of ${\rm supp}[u]$.  See Proposition \ref{P:2.2} below for its equivalent characterizations.

\ms
 The main result of this paper is the following stability result on the (scale invariant) EHI. 
See Definition \ref{D:ehi} for a precise definition of the EHI.

	\begin{thm} \label{T:main-stable}
		Let $(\sX,d)$ be a complete    doubling  metric space, 
		and let $m$ be a  
		Radon measure on $\sX$ with full support.
		Let  
		$(\sE,\sF)$ be a strongly local  symmetric regular Dirichlet form on $L^2(\sX; m)$. 
		Suppose that  $(\sX,d,m,\sE,\sF)$ satisfies the EHI. 
		Let $(\sE', \sF)$ be another strongly local symmetric regular Dirichlet form on $L^2(\sX; m)$ such that
		$$ 
			C^{-1}  \sE(f,f) \le  \sE'(f,f) \le C  \sE(f,f)   \quad \hbox{ for all  } f \in \sF 
		$$
		for some $C \ge 1$.
		Then $(\sX,d,m,\sE',\sF)$ satisfies the EHI. 
	\end{thm}

 \ms
  
 Theorem \ref{T:main-stable} is established based on  
   characterizations of the EHI given in Theorem \ref{T:main-new}.
 This stability result  is further extended in Theorem \ref{T:7.11} to  strongly local MMD spaces that may have different
 symmetrizing measures. 
 
 \medskip
 
 Suppose that $(\sX, d)$ is a locally compact complete length space. By \cite[Theorem 3.11]{BM1}, $(\sX, d)$ is metric doubling
if   $(\sX,d,m,\sE,\sF)$ satisfies the EHI and has regular Green functions.
 Theorem \ref{T:main-stable} together with Theorems \ref{T:3.6} and \ref{T:Rgreen} readily implies the following corollary, which 
substantially improves the main result, Theorem 1.3, of \cite{BM1}. Alternately, Corollary \ref{C:main-stable} follows from Theorems \ref{T:main-stable}
	and \ref{T:ehitomd} since the condition (c) in Theorem \ref{T:ehitomd} follows from \cite[Theorem 2.5.28]{BBI}.

\begin{corollary} \label{C:main-stable}
Let  $(\sX,d)$ be a complete locally compact length space, 
		and  let $m$ is a  
		Radon measure on $\sX$ with full support.
		Let  
		$(\sE,\sF)$ be a strongly local  symmetric regular Dirichlet form on $L^2(\sX; m)$. 
		Suppose that  $(\sX,d,m,\sE,\sF)$ satisfies the EHI. 
		Let $(\sE', \sF)$ be another strongly local symmetric regular  Dirichlet form on $L^2(\sX; m)$ such that
		for some $C \ge 1$,
		$$ 
			C^{-1}  \sE(f,f) \le  \sE'(f,f) \le C  \sE(f,f)   \quad \hbox{ for all  } f \in \sF. 
		$$
		Then $(\sX,d,m,\sE',\sF)$ satisfies the EHI. 
\end{corollary} 

  \sms
The above corollary implies a generalized version of Moser's EHI \cite{Mo1} from $\bR^n$ to an arbitrary Riemannian manifold. Let $(M,g)$ be a Riemannian manifold and let $\operatorname{Sym}(TM)$ denote the bundle of symmetric endomorphisms of the tangent bundle $TM$. We say that $\sA$ is a \emph{uniformly elliptic operator in divergence form} if there exists $A: M \to \operatorname{Sym}(TM)$ a measurable section of $\operatorname{Sym}(TM)$ and a constant $K \ge 1$ such that 
\[
K^{-1} g(\xi,\xi) \le g(A\xi,\xi) \le K g(\xi,\xi)\quad \mbox{for all $\xi \in TM$,}
\]
such that $\mathcal{A}(\cdot)= \div(A \nabla(\cdot))$, where $\div$ and $\nabla$ denote the Riemannian divergence and gradient respectively. This extends \cite[Theorem 1.4(a)]{BM1} without any assumption on curvature. The description of the Dirichlet form in this context is given in \cite[\textsection 2.2.5]{CF}. The following generalized Moser's EHI follows  from Corollary \ref{C:main-stable}.
 \begin{cor} \label{c:moser}
 	Let $(M,g)$ be a Riemannian manifold and let $\Delta$ denote the corresponding Laplace-Beltrami operator. If $(M,g)$ satisfies EHI for non-negative  solutions on $\Delta u=0$, then it satisfies EHI for non-negative solutions of $\sA u=0$, where $\sA$ is any uniformly elliptic operator in divergence form.
 \end{cor}
	
The remainder of this paper is organized as follows. 	In Section \ref{S:2}, we present definitions and terminology associated
with regular symmetric Dirichlet forms 
as well as  some basic facts that will be used in this paper.  Existence and regularity of Green functions
are given in Section \ref{S:lreg} for transient regular Dirichlet forms. 
The transience condition is removed in Section \ref{S:4}. It is shown there that any strongly local regular Dirichlet form $(\sE, \sF)$ 
on a connected locally compact separable metric space $\sX$ 
that satisfies the local EHI is irreducible and has regular Green functions.  
Various consequences of the EHI are presented in Section \ref{S:5}.
In particular, it is shown that  for a complete locally compact separable metric space $(\sX, d)$,  under the EHI, 
relatively ball connected, metric doubling and quasi-arc connected properties are all mutually equivalent.
 In     Section \ref{S:6},  a good doubling  measure  $\mu$  is constructed on a MMD space $(\sX, d, m, \sE, \sF)$   
 that satisfies the EHI and is relatively ball connected.  
This measure relates well with capacities and is a smooth measure with full  quasi support on $\sX$. 
It is shown in Section \ref{S:7} that the Dirichlet form time-changed by the positive continuous additive functional generated by this doubling measure $\mu$ 
is a MMD space $(\sX, d, \mu, \sE, \sF^\mu)$  that satisfies Poincar\'e inequality ${\rm PI}(\Psi)$, the cutoff energy inequality 
${\rm CS}(\Psi)$ and a capacity estimate ${\rm cap}(\Psi)$, where $\Psi$ is a suitable regular scale function.
From this we  can obtain  equivalent characterizations of the EHI   in Theorem \ref{T:main-new}, and 
deduce the stability result of the EHI stated in Theorem \ref{T:main-stable}.   
 The   scale function $\Psi$ varies both in space and in time;  
\cite{Telcs} first studied such location dependent scaling functions in detail.
 An extension of Theorem \ref{T:main-stable} is given at the end of Section \ref{S:7}  
that the second Dirichlet form $\sE'$ may have symmetrizing measure $\mu$ different from $m$;
see Theorem \ref{T:7.11}. 
Three examples are given in Section \ref{S:8}. The first example shows that without certain regularity of the metric,
the stability of the EHI may fail. 
The second example is of a strongly local regular Dirichlet form that fails to satisfy the
non-scale-invariant Harnack inequality. 
The third one is of a space which satisfies the EHI and is covered by the results of this paper,
but fails to satisfy the local regularity required in \cite{BM1}.

\section{Preliminaries}\label{S:2} 

In this section, we give definitions of some terminology from Dirichlet form theory that are used in this paper
and some basic facts.  We refer the reader to \cite{CF, FOT} for more details on the theory of symmetric Dirichlet forms. 
We use $:=$ as a way of definition.

 Let $(\sX, \sB (\sX))$  be a  measurable space and $m$ a $\sigma$-finite measure on $\sX$.
  A   bilinear form $(\sE, \sF)$ on
 $L^2 (\sX; m)$ is said to be a  {\it symmetric Dirichlet form}  if 
 \begin{description}
 \item{(i)} $\sF$ is a dense linear subspace of $L^2(\sX; m)$; 
 
 \item{(ii)} $\sE$ is symmetric and bilinear on $\sF\times \sF$ such that    $\sE(f, f)\geq 0$ for every $f\in \sF$;

 \item{(iii)} $\sF$ is a Hilbert space with inner product $\sE_1(f, g):=\sE(f, g) + \int_{\sX} f(x) g(x) m(dx)$; 
 
 \item{(iv)}  For every $f\in \sF$, $g:=(0\vee f) \wedge 1$ is in $\sF$ and $\sE(g, g)\leq \sE(f, f)$.
 \end{description}
 A bilinear form $(\sE, \sF)$ on $L^2(\sX; m)$ satisfying properties (i)-(iii) above is called a {\it symmetric closed form}.
 Any symmetric closed form is in one-to-one correspondence with a strongly continuous symmetric 
 contraction semigroup $\{T_t; t\geq 0\}$
 on $L^2(\sX; m)$.
 Property (iv) above is called the Markovian property which is equivalent to the corresponding semigroup $\{T_t; t\geq 0\}$   being Markovian; that is, $0\leq T_tf\leq 1$ for any $f\in L^2(\sX; m)$ with $0\leq f\leq 1$. 
A real-valued function $f$ is said to be in the {\it extended Dirichlet space} $\sF_e$ if there is an $\sE$-Cauchy sequence
$\{f_k; k\geq 1\} \subset \sF$ so  that $\lim_{k\to \infty} f_k =f $ $m$-a.e. on $\sX$, and we define 
$\sE(f, f)= \lim_{k\to \infty}  \sE(f_k, f_k)$.  Clearly, $\sF \subset \sF_e$. 
It is known that $\sF=\sF_e \cap L^2(\sX; m)$; see \cite[Theorem 1.1.5(iii)]{CF}. 

\medskip

 The Dirichlet form $(\sE, \sF)$ on $L^2(\sX; m)$ is said to be {\it transient} if 
  there exists a bounded $g \in L^1(\sX; m)$ that is strictly positive on $\sX$ 
  so that 
  $$
   \int_{\sX} |u(x)| g(x) m(dx) \leq \sE(u, u)^{1/2} \quad \hbox{for every } u\in \sF.
  $$
   Clearly, if $(\sE, \sF)$ is transient, then $(\sF_e, \sE)$ is a Hilbert space.
   The Dirichlet form $(\sE, \sF)$ on $L^2(\sX; m)$ is said to be {\it recurrent} if $1\in \sF_e$ and $\sE(1, 1)=0$. 
   Denote by $\{T_t; t\geq 0\}$ the semigroup on $L^2(\sX; m)$ corresponding to the Dirichlet form $(\sE, \sF)$.
   By Theorem 2.1.5 and Theorem 2.1.8 of \cite{CF}, 
   $(\sE, \sF)$ is transient if and only if there is some $g \in L^1(\sX; m)$ that is strictly positive on $\sX$
   and satisfies $G g:=\int_0^\infty T_t g \, dt <\infty$ $m$-a.e. on $\sX$;
   and $(\sE, \sF)$ is recurrent if and only if for any non-negative $g$ on $\sX$ with $\int_{\sX} g(x) m(dx)<\infty$, 
      $Gg \in \{0,\infty\}$  $m$-a.e. on $\sX$. 
   
 \medskip
   
  Denote by ${\cal B}^m(\sX)$  the completion of  the field ${\cal B}(\sX)$   under the measure 
	$m$.  A set $A\in {\cal B}^m(\sX)$ is said to be {\it $\{T_t\}_{t\geq 0}$-invariant}  if $T_t (1_{A^c} f )=0$
	$m$-a.e. on $A$ for all $t>0$ and $f\in L^2(\sX; m)$. 
	By \cite[Proposition 2.1.6]{CF}, $A\in {\cal B}^m(\sX)$ is    $\{T_t\}_{t\geq 0}$-invariant if and only if 
	  $1_A u \in \sF$ for every $u\in \sF$ and 
	 \begin{equation}\label{e:2.1}
	 \sE (u, v)= \sE (1_Au, 1_{A} v) + \sE(1_{A^c} u, 1_{A^c} v)  \quad \hbox{for every } u, v\in \sF.
	 \end{equation} 
The Dirichlet form $(\sE, \sF)$ on $L^2(\sX; m)$ is said to be {\it irreducible} if for any  $\{T_t\}_{t\geq 0}$-invariant
set $A$, either $m(A)=0$ or $m(A^c)=0$.   An irreducible Dirichlet form is either transient or recurrent; see
\cite[Propositions 2.1.3(iii) and 2.1.6]{CF}.   

\medskip 

A Dirichlet form $(\sE, \sF)$ on $L^2(\sX; m)$ is said to be {\it regular}  if 
 \begin{description}
 \item{(i)} $(\sX, d) $ is a locally compact separable  metric space and $m$ is a Radon measure on $\sX$ with full support;
   
  \item{(ii)}  $\sF \cap C_c(\sX) $ is $\sqrt{\sE_1}$-dense in $\sF$, where 
   $C_c(\sX )$  is the space of continuous functions on $\sX$ having compact support;
  
   \item{(iii)}  
  $\sF\cap C_c( \sX ) $ is dense   in $ C_c(\sX)$ with respect to the uniform norm
   $\| f\|_\infty= \sup_{x\in \sX} |f(x)|$ . 
  
 \end{description} 

\medskip

For a regular Dirichlet form $(\mathcal{E},\mathcal{F})$ on $L^2(\mathcal{X},m)$, an increasing sequence $\{F_k; k\geq 1\}$ of closed subsets   of $\sX$ is said to be an  {\it $\sE$-nest} if
$\cup_{k\geq 1} \sF_{F_k} $ is $\sqrt{\sE_1}$-dense in $\sF$, where $\sF_{F_k}:=\{f\in \sF: f=0 \  m \hbox{-a.e. on } 
\sX\setminus F_k \}$. 
A  set $N\subset \sX $ is said to be \emph{$\sE$-polar} if there is an $\sE$-nest  $\{F_k; k\geq 1\}$ so that $N\subset \sX \setminus 
\cup_{k\geq 1}  F_k$.  An $\sE$-polar set $A$ always has $m(A)=0$. 
$\sE$-polar sets can also be characterized by using capacity. 
Given a regular Dirichlet form  $(\sE, \sF)$ on $L^2(\sX; m)$, we can define 1-capacity ${\rm Cap}_1$  as follows.
For any open  subset $U\subset \sX$, 
\begin{equation}\label{e:2.2}
{\rm Cap}_1 (U):= \inf\{ \sE_1 (f, f): f\in \sF,  \ f\geq 1 \hbox{ $m$-a.e. on } U\}
\end{equation}
with the convention that  $\inf \emptyset :=\infty$, and for any subset $A\subset \sX$,
\begin{equation}\label{e:2.3} 
{\rm Cap}_1 (A):= \inf\{ {\rm Cap}_1 (U): U \subset \mathcal{ X} \mbox{ open}, U\supset A\}. 
\end{equation}
It is known (see \cite[Theorem 1.3.14]{CF}) that for a regular Dirichlet form  $(\sE, \sF)$ on $L^2(\sX; m)$, $A\subset \sX$ is $\sE$-polar if and only if
it has zero 1-capacity. 
A statement depending on $x\in A$ is said to hold {\it $\sE$-quasi-everywhere}  ($\sE$-q.e. in abbreviation) if 
there is an $\sE$-polar set $N\subset A$ so that the statement is true for every $x\in A\setminus N$.
A function $f$ is said to be {\it $\sE$-quasi-continuous}  on $  \sX$ if there is an $\sE$-nest $\{F_k; k\geq 1\}$  so that 
$\restr{f}{F_k}\in C(F_k )$ for every $k\geq 1$, where $C(F_k):= \{u:F_k \to \bR | \mbox{$u$ is continuous}\}$.  When there is no possible ambiguity, we often drop  ``$\sE$-"  from 
$\sE$-quasi-everywhere and $\sE$-quasi-continuous. 
For a regular  Dirichlet form  $(\sE, \sF)$ on $L^2(\sX; m)$, every $f\in \sF_e$ has an $m$-version that 
is quasi-continuous on $\sX$, which is unique up to an $\sE$-polar set; see \cite[Theorem 2.3.4]{CF} or 
\cite[Theorem 2.1.7]{FOT}. We always take a function 
$f$ in  $\sF_e$ to be represented by its quasi-continuous version.

Recall that a Hunt process $X=\{X_t, t\geq 0; \bP^x, x\in \sX\}$ on a locally compact separable metric space $\sX$
is a strong Markov process that is right continuous and quasi-left continuous on the one-point compactification
$\sX_\partial:=\sX \cup\{\partial \}$ of $\sX$.  
A set $C\subset \sX_\partial$ is said to be {\it nearly Borel measurable} if for any probability measure $\mu$ on $\sX$
there are Borel sets $A_1, A_2$ such that $A_1\subset C\subset A_2$ and 
$$ 
\bP^\mu ( \hbox{there is some   $t\geq 0$ such that } X_t\in A_2\setminus A_1 )=0.
$$ 
   Let $m$ be a Radon measure with full support on $\sX$.
A Hunt process $X$ is said to be {\it $m$-symmetric}  if the transition semigroup is symmetric on $L^2(\sX; m)$.
For an $m$-symmetric Hunt process $X$ on $\sX$, 
a set $\sN\subset \sX$ is said to be {\it properly exceptional} for $X$ 
if $\sN$ is nearly Borel measurable, $m(\sN) =0$ and 
$$
\bP^x ( X_t\in \sX_\partial \setminus \sN \hbox{ and } X_{t-} \in \sX_\partial \setminus \sN \hbox{ for all } t>0)=1
\quad \hbox{for every } x\in \sX \setminus \sN.
$$

In 1971, Fukushima showed that any  symmetric regular  Dirichlet form  $(\sE, \sF)$ on $L^2(\sX; m)$
has an $m$-symmetric  Hunt process $X=\{X_t, t\geq 0; \bP^x, x\in \sX\}$ on $\sX$ associated with it
in the sense that the transition semigroup of $X$ is  a version of  the  strongly continuous semigroup
$\{T_t; t\geq 0\}$ on $L^2(\sX; m)$ corresponding to $(\sE, \sF)$,  see \cite[Theorem 7.2.1]{FOT}. 
Furthermore, for any non-negative Borel measurable $f\in L^2(\sX; m)$ and $t>0$, 
$$ P_t f(x):= \bE^x[ f(X_t)]
$$
 is a quasi-continuous version of $T_tf$ on $\sX$.
The Hunt process $X$ associated with a regular  Dirichlet form  $(\sE, \sF)$ on $L^2(\sX; m)$ 
is unique in the following sense (see \cite[Theorem 4.2.8]{FOT}): 
 if $X'$ is another Hunt process associated 
with the regular  Dirichlet form  $(\sE, \sF)$ on $L^2(\sX; m)$, then there is a common properly exceptional set
outside which these two Hunt processes have the same transition functions.
 We say the $m$-symmetric Hunt process $X$ on $\sX$ is {\it transient}, {\it recurrent}, and {\it irreducible} 
 if so is its associated Dirichlet form $(\sE, \sF)$ on $L^2(\sX; m)$.

 \medskip
 
 In the remainder of this section,   $(\sE, \sF)$ is a regular Dirichlet form on $L^2(\sX; m)$
 and   $X=\{X_t, t\geq 0; \bP^x, x\in \sX\}$ is the  Hunt process associated with it. 
 Let $\zeta$ denote the lifetime of $X$, and $\{\sF_t; t\geq 0\}$ be the minimum augmented filtration generated by $X$.

  A subset $\sN\subset \sX$ is said to be
{\it $m$-polar}    if there is a nearly Borel set $\sN_1\supset \sN$ so that 
$\bP^x (\sigma_{\sN_1}<\infty)=0$ for $m$-a.e. $x\in \sX$, where $\sigma_{\sN_1}=\inf\{t>0: X_t \in \sN_1\}$.  
It is known that a subset $\sN\subset \sX$ is $\sE$-polar if and only if
it is $m$-polar, and any $\sE$-polar set is contained in a Borel properly exceptional set for $X$;
see \cite[Theorems 3.1.3 and 3.1.5]{CF}.

If $(\sE, \sF)$ is irreducible, then (see \cite[Theorem 3.5.6]{CF}) 
for any non-$\sE$-polar nearly Borel measurable set $A$,
\begin{equation}\label{e:2.4} 
\bP^x(\sigma_A<\infty) >0 \quad \hbox{for  $\sE$-q.e. } x\in \sX.
\end{equation} 
 Let $D$ be an open subset of $\sX$. The part process $X^D$ of $X$ killed upon exiting $D$ is a Hunt process
 on $D$ whose associated Dirichlet form 
   $(\sE^D, \sF^D)$ 
on $L^2(D; m|_D)$ is   regular.
Here $m|_D$ is the measure $m$ restricted to the open set $D$,  
\begin{equation}\label{e:FD} 
\sF^D=\{f\in \sF: f=0 \hbox{ $\sE$-q.e. on } D^c\},
\end{equation} 
 and $\sE^D=\sE$ on $\sF^D$; 
see, e.g.,  Exercise 3.3.7 and  Theorem 3.3.9 of \cite{CF}. 
 It is well known (see, e.g.,  \cite[Theorem 3.3.8]{CF}) that $A\subset D$ is $\sE^D$-polar if and only if it is $\sE$-polar.  
 Property \eqref{e:2.4} combined with \cite[Proposition 2.1.10]{CF} yields the following.

\begin{proposition}\label{P:2.1}  
If $(  \sE, \sF)$ 
is irreducible and $D^c$ is not $\sE$-polar, then
   the regular Dirichlet form   $(\sE^D, \sF^D)$  on $L^2(D; m|_D)$ is transient.  
 \end{proposition} 
 
   \medskip
   
For $u\in \sF_e$, the following Fukushima decomposition holds (see \cite[Theorem 4.2.6]{CF} or \cite[Theorem 5.2.2]{FOT}): for  $\sE$-q.e.~$x \in \mathcal{X}, \mathbb{P}^x$-a.s.,
\begin{equation}\label{e:Fu} 
u(X_t)-u(X_0)=M^u_t+N^u_t,  \quad t\geq 0,
\end{equation}
where $M^{u}$ is a martingale additive functional of $X$ having finite energy and $N^{u}$ is a continuous additive
	functional of $X$ having zero energy. The predictable quadratic variation $\<M^u\> $ of  the square-integrable martingale 
	$M^u$ is a positive continuous
	additive functional of $X$, whose corresponding Revuz measure is denoted by $\mu_{\<u\>}$. We call $\mu_{\<u\>}$
	the energy measure of $u\in \sF_e$. It is known that 
 $$
 \frac12 \mu_{\<u\>}(\sX) \leq \sE( u, u) \leq   \mu_{\<u\>}(\sX)  \quad \hbox{for } u\in \sF_e.
 $$
When $(\sE, \sF)$ admits no 
 killing inside $\sX$, which is equivalent to the Hunt process  admitting no 
 killing  inside $\sX$
(that is, $\bP^x (X_{\zeta-} \in \sX, \zeta <\infty)=0$ for $\sE$-q.e. $x\in \sX$),
we have 
	\begin{equation}\label{e:2.6}
	 \sE( u, u)  =\frac12   \mu_{\<u\>}(\sX)  \quad \hbox{for } u\in \sF_e.
	 \end{equation} 
	 When $u\in \sF_e$ is bounded, its energy measure $\mu_{\<u\>}$ 
	 can be computed by the formula 
\begin{equation}\label{e:2.7}
\int_{\sX} v (x) \mu_{\<u\>} (dx) = 2 \mathcal{E}(u,uv)- \mathcal{E}(u^{2}, v) \qquad \hbox{for all bounded } 
v \in \mathcal{F}. 
\end{equation} 
For general $u\in \sF_e$,   $\mu_{\<u\>}$ is the increasing limit of $\mu_{\<u_n\>}$ 
as $n\to \infty$, where $u_n:= (-n)\vee(u\wedge n) \in \sF_e$. 
See (4.3.12)-(4.3.13),  and Theorems 4.3.10 and 4.3.11 of \cite{CF} for the above stated properties  of  $\mu_{\<u\>}$. 

\medskip

The following is taken from Theorem 2.4.3 and Theorem 4.3.4 of \cite{CF}.

\begin{proposition}\label{P:2.2}   The following are equivalent.
\begin{description}
\item{\rm (i)} $(\sE, \sF)$ is strongly local  (i.e., $\mathcal{E}(u,v)=0$ whenever $u,v \in \mathcal{F}$, the support of $u$ is compact and $v$ is constant on a neighborhood of the support of $u$); 

\item{\rm (ii)} $\sE(u, v)=0$ whenever $u, v\in \sF$ with $u(v-c)=0$ $m$-a.e. on $\sX$ for some constant $c$;

\item{\rm (iii)} The associated Hunt process $X$ is a diffusion with no 
 killing
 inside $\sX$; that is, there is a Borel 
properly exceptional set $\sN_0\subset \sX$ so that for every $x\in \sX \setminus \sN_0$, 
\begin{equation} \label{e:2.8}
\bP^x(X_t  \hbox{ is continuous in } t\in [0, \zeta)) =1 \hbox{ and }  \quad 
\bP^x (X_{\zeta-} \in \sX, \zeta <\infty) =0 . 
\end{equation} 
 \end{description}
\end{proposition} 
   
  \medskip
  
In Theorem \ref{T:3.5} below,  a new  criterion for irreducibility will be given for   strongly local regular
Dirichlet forms. 

\medskip
  
	We use notation $V \Subset D$ for $V$ being a relatively compact open subset of $D$.
	 	For any open set $U$, we define
\begin{equation}\label{e:Floc}
 \mathcal{F}_{\rm loc}^U
:= \Biggl\{ f \Biggm|
\begin{minipage}{290pt}
$f$ is an $m$-equivalence class of $\mathbb{R}$-valued Borel measurable functions
on $\sX$ such that 
 for each  open $V\Subset U$, there is some $g\in\mathcal{F}^U$
so that $f=g$ $m$-a.e. on $V$. 
\end{minipage}
\Biggr\}.
\end{equation}
 Note that each $f\in \mathcal{F}_{\rm loc}^U$ admits an $m$-version that is  $\sE$-quasi-continuous on $U$,
which is unique modulo an $\sE$-polar set.
We always let a function in $\mathcal{F}_{\rm loc}^U$ be represented by its quasi-continuous version.
When $U=\sX$, we simply write $\mathcal{F}_{\rm loc}$ for $\mathcal{F}_{\rm loc}^{\sX}$. 

\medskip

When the Dirichlet form $(\sE, \sF)$ is strongly local, 
the energy  measure $\mu_{\<u\>}$ has the following  strong local property;   see
 \cite[Proposition 4.3.1 and Theorem 4.3.10(i)]{CF}. 

\begin{proposition}\label{P:2.3}
Suppose the Dirichlet form  $(\sE, \sF)$ is strongly local and  $D$ is an open subset of $\sX$.
Then 
\begin{description}
\item{\rm (i)} $\mu_{\<u\>} (D)=0$ if $u\in \sF_e$ and $u$ is constant $\sE$-q.e. on $D$;
 
\item{\rm (ii)} $\mu_{\<u\>}= \mu_{\<v\>}$ on $D$  for every $u, v\in \sF$ so that 
 $u-v$ is a constant $\sE$-q.e. on $D$. 
 \end{description}
 \end{proposition}

   Let $\{U_k; k\geq 1\}$ be an increasing sequence of relatively compact open subsets whose union is $\sX$.
	For $u\in \mathcal{F}_{\rm loc}$, there is some $u_k\in \sF$ so that $u_k=u$ $m$-a.e. on $U_k$.
      	Define $\mu_{\<u\>}=  \mu_{\<u_k\>}$ on $U_k$.
	  Since $(\sE, \sF)$ is strongly local, $\mu_{\< u\>}$ is uniquely defined by Proposition \ref{P:2.3}(ii).
	  In view of \eqref{e:2.6}, this allows us to extend the definition of $\sE$ to $\mathcal{F}_{\rm loc}$ by setting 
\begin{equation}\label{e:2.10}
	\sE (u, u)  := \frac12 \mu_{\<u\>} (\sX), \quad u\in \mathcal{F}_{\rm loc}.
\end{equation} 
    
In this paper, we will use time change of Dirichlet form and its associated Hunt process
so we need the notion of smooth measure. The following definition is  from \cite[Definition 2.3.13]{CF}. 
 
	\begin{definition} \label{d:smooth}  \rm 
	Let $(\sE,\sF)$ be a  regular Dirichlet form on $L^2(\sX; m)$. 
		A  (positive) Borel measure $\mu$ on $\sX$ is \emph{smooth} if it satisfies the following conditions:
		\begin{enumerate}[(a)]
              \item $\mu$ charges no $\sE$-polar set; 
	\item there exists an $\sE$-nest $\set{F_k}$   such that 
			$	\mu(F_k) < \infty$ for all $k \geq 1$.
		 \end{enumerate}
	\end{definition}

	By \cite[Theorem 1.2.14]{CF}, the above definition of smooth measure is equivalent to that defined in
	\cite[p.83]{FOT}.	  
	Clearly every positive Radon measure charging no $\sE$-polar set  is smooth, as in this case we can take the $\sE$-nest 
	$\{F_k\}$ to be the closures of an increasing sequence of relatively compact open sets whose union is $\sX$.   
We say $D \subset \sX$ is \emph{quasi open} if there exists 
	 	an $\sE$-nest 
	 $\set{F_n}$ such that $D \cap F_n$ is an open subset of $F_n$ in the relative topology for each $n \in \bN$. The complement of a quasi open set is called \emph{quasi closed}.

	\begin{definition} \label{d:qsupp} \rm   (See \cite[Definition 3.3.4]{CF} or \cite[p.190]{FOT}.)
			Let $\mu$ be a  smooth Borel measure on $\sX$. A set $F \subset \sX$ is called a \emph{quasi support} of $\mu$ if it satisfies the following: \\
(a) $F$ is quasi closed and $\mu (\sX \setminus F)=0$. \\
(b) 			If $\tilde{F}$ is another set that satisfies (a), then 
 $F \setminus \tilde{F}$ is $\sE$-polar. \\
 Such a set $F$ exists by \cite[Theorem 4.6.3]{FOT}.
We say that $\mu$ has \emph{full quasi support} if $\sX$ is a quasi support of $\mu$.
	 \end{definition}

\medskip

 We assume in the remaining of this paper  that
 $(\sE, \sF)$ is a symmetric  strongly local  regular Dirichlet form on $L^2(\sX; m)$. We call
$(\sX, d, m, \sE, \sF)$ a \emph{metric measure Dirichlet  (MMD) space}. 
Sometimes, to emphasize its dependence on the symmetrizing measure, we write $\sF^m$ for $\sF$. 
 Let $X=\{X_t, t\geq 0; \bP^x, x\in \sX\}$ be the diffusion process associated with $(\sX, d, m, \sE, \sF)$,
  whose lifetime is denoted
by $\zeta$ and whose shift operators are denoted by $\{\theta_t\}_{t \ge 0}$. 
The one-point compactification of the locally compact metric space $(\sX, d)$ is denoted as 
 $\sX_\partial :=\sX \cup \{\partial\}$. 
 	 For a nearly Borel measurable set $A \subset \sX_\partial$,  define the stopping times
	$$ \sigma_A= \inf\{t > 0: X_t \in A \}, \quad \tau_{A} = \sigma_{\sX_\partial \setminus A}= \inf\{t>0: X_t \notin A\}, $$ 
	and write $\tau_{x}$ for $\tau_{ \{x\} }$.	Note that by definition, $\tau_A \leq \zeta$ if $A \subset \sX$. 
 
\medskip    
 
 Here and in the following, we use the convention 
that $X_\infty :=\partial$,
and that any function $u$ defined on a subset of $\sX$
is extended to $\{\partial\}$ by taking $u(\partial)=0$.

	The following definition is taken from \cite[Definition 2.1]{C}.

	\begin{defn}\label{D:2.6}
\rm  Let $D$ be an open subset of $\sX$. We say a universally measurable function
$u$ defined  $\sE$-q.e.\ on $D$ is {\it harmonic} in $D$ (with respect to the process $X$) if
  for every relatively compact open subset $U$ of $D$,
  $t\mapsto u(X_{t\wedge \tau_U})$
 is a uniformly integrable $\bP^x$-martingale for  $\sE$-q.e. $x\in U$.
We say that a universally measurable function $u$ on $\overline{D}$  is {\it regular harmonic} in $D$ (with respect to the process $X$) if
$\bE^x [ |u(X_{\tau_D})|] <\infty$ and 
$u(x)=\bE^x [ u(X_{\tau_D})]$  for  $\sE$-q.e. $x\in D$.
\end{defn}
 
\begin{remark}\label{R:2.7} \rm
\begin{enumerate} [(i)]
\item Observe that if $u$ is a universally measurable function on $\overline{U}$ and if $\bP^x (\tau_U <\infty)=1$ for  $\sE$-q.e. $x\in U$, then $t\mapsto u(X_{t\wedge \tau_U})$
 is a uniformly integrable $\bP^x$-martingale for  $\sE$-q.e. $x\in U$ if and only if $\bE^x [ |u(X_{\tau_U})|] <\infty$ and
$u(x)= \bE^x \left[ u (X_{\tau_U}) \right]$ for  $\sE$-q.e. $x\in U$. 
  Sufficient conditions on
 $\bP^x (\tau_U <\infty)=1$ for  $\sE$-q.e. $x\in U$ are given in  Propositions \ref{P:3.1} and \ref{P:3.2} below. 

\item  By the Markov property of $X$, it is clear that any regular harmonic function in $D$ is harmonic in $D$. 
\end{enumerate} 
\end{remark}

The relation  of the above probabilistic notion of harmonicity to the analytic notion of harmonicity
has been investigated in \cite{C} for general symmetric regular Dirichlet forms.
In the setting of strongly local symmetric regular Dirichlet forms, one direction becomes much easier to analyze; cf. \cite[Theorem 2.7 and Remark 2.8]{C}. 

	\begin{defn}\label{D:2.7}
\rm  Let $D$ be an open subset of $\sX$. We say a function
$u$ is {\it $\sE$-harmonic} in $D$  if  $u\in \sF^D_{\rm loc}$ 
  and 
	\begin{equation}\label{e:2.12} 
	\sE(u, v)=0  \quad \hbox{for every } v \in C_c(D)\cap \sF .
  \end{equation}
\end{defn}

Note that as explained in \eqref{e:2.10}, 
the bilinear form $\sE(u, v)$ in \eqref{e:2.12} is well defined due to the strong locality of
$(\sE, \sF)$.

\begin{prop}\label{P:2.9}  Let $D$ be an open subset of $\sX$.
\begin{enumerate} [{  \rm (i)}]
\item Any $\sE$-harmonic  and $\sE$-quasi-continuous function in $D$ is harmonic in $D$. 
 
 \item
Moreover, any $\sE$-harmonic and $\sE$-quasi-continuous function  in $D$
 is regular harmonic in any relatively compact open subset $U$ of $D$ such that 
\begin{equation} \label{e:exit}
	\bP^x(\tau_U < \infty) >0 \quad \mbox{for  $\sE$-q.e.~$x \in U$.}
\end{equation}  

\item If $u$ is locally bounded and harmonic in $D$, 
then $u$ is $\sE$-quasi-continuous on $D$ and $\sE$-harmonic in $D$.
\end{enumerate} 
\end{prop}

\pf     The proof for (i) and (ii)   is essentially a particular  case of \cite[Theorem 2.7]{C}.
For  strongly local regular Dirichlet forms, the proof can be much simplified.
For the reader's convenience, we present a proof here. 

\sm   (i) Suppose that $h$ is an $\sE$-harmonic function in $D$. 
Let $U$ be any relatively compact open subset of $D$. 
  There is a relatively compact open subset $V$ of $D$ so that 
$\overline U \subset V$. 
Since $h\in \sF^D_{\rm loc}$,
there is some $f\in \sF^D$ so that $f=h$ $m$-a.e. on $V$ and hence $\sE$-q.e. on $V$ as they are all represented by
their $\sE$-quasi-continuous versions. Since $h$ is  $\sE$-harmonic   in $D$ and $(\sE, \sF)$ is strongly local,  we have
\begin{equation}\label{e:2.13} 
\sE (f, v)=\sE(h, v)=0 \quad \hbox{for every } v \in C_c(U)\cap \sF .
\end{equation} 
 Set $g(x):= \bE^{x}[f(X_{\tau_U})]$. Then by \cite[Theorem 3.4.8]{CF} or \cite[Theorem 4.6.5]{FOT}, $g \in \sF_e$, $g=f$  $\sE$-q.e.~on $U^c$ and 
 $g$ is $\sE$-harmonic in $U$. Since $f-g \in \sF_e$ with $f-g=0$ $\sE$-q.e. on $U^c$, we have $\sE(f-g,f-g)=0$. By \cite[Lemma 2.2]{C}, 
\be \label{e:const}
\bP^x(f(X_t)-g(X_t)= f(X_0)-g(X_0) \mbox{ for every $t \ge 0$})=1
\quad \sE \hbox{-q.e. on } \sX. 
\ee
Therefore 
$t \mapsto h(X_{t \wedge \tau_U}) = f(X_{t \wedge \tau_U}) = f(X_0)-g(X_0) +g( X_{t \wedge \tau_{U}})$ is a uniformly integrable $\bP^x$-martingale for $\sE$-q.e.~$x \in U$. So  $h$ is harmonic in $D$. This proves (i).

 \sm (ii) We continue with the notation of part (i), and assume that \eqref{e:exit} holds.
  Note that by \cite[Theorems 3.3.9 and 3.4.9]{CF}, $f-g$ is in the extended Dirichlet space
for the part Dirichlet form $(\sE^U, \sF^U)$ with $\sE^U (f-g,f-g) = \sE(f-g,f-g)=0$.
 Under condition \eqref{e:exit}, we have by 
 \cite[Lemma 2.2]{C}  applied to  $(\sE^U, \sF^U)$  that  
 $h=f=g$ $\sE$-q.e. on $U$. Thus  for $\sE$-q.e. $x\in U$, 
 $h(x)= g(x)= \bE^{x}[f(X_{\tau_U})]=\bE^{x}[h(X_{\tau_U})]$, where the last equality is due to the fact that $h=f$ $\sE$-q.e. on $V\supset \overline U$.
 This proves that 
 $h$ is regular harmonic  in $U$.

\sm (iii)    The proof of this part is the same as that of  \cite[Theorem 2.9]{C} except that the first sentence in the proof there
  needs some details as follows. 
Without loss of generality, we assume that the Hunt process $X$ associated with the strongly local regular DIrichlet form $(\sE, \sF)$ is
defined on the canonical sample space $\Omega$, which is the space of continuous functions on $[0, \infty)$ taking values in 
$\sX_\partial:=\sX\cup \{\partial\}$.  Denote by $\{\theta_t; t\geq 0\}$  time shifting operators on $\Omega$. 
 Let $\{\sF_t; t\geq 0\}$ be the minimal augmented filtration generated by $X$, which is known to be right continuous.
Suppose that  $u$ is locally bounded and harmonic in $D$. 
Then for every relatively compact open subset $U$ of $D$, there is a properly exceptional set $\sN\subset U$
so that 
  $t\mapsto u(X_{t\wedge \tau_U})$
 is a uniformly integrable (in fact bounded) $\bP^x$-martingale for every $x\in U \setminus \sN$.
  Let $\Omega_0$ be the collections of those $\omega \in \Omega$ so that 
 \begin{equation}\label{e:2.16}
\xi  (\omega) :=\lim_{n \in \bN, n\to \infty}   u(X_{n \wedge \tau_U (\omega)} ) (\omega)  \hbox{ exists as a finite number}. 
\end{equation}
 From the discrete martingale theory,  we know that for $x\in U\setminus \sN$,  
  $\bP^x (\Omega_0)=1$ and  $u(X_{n \wedge \tau_U})$ converges to $\xi$ in $L^1(\bP^x)$   as $n\to \infty$. 
  This together with the strong Markov property of $X$ implies  that for each  $x\in U\setminus \sN$ and 
  $t\geq 0$, $\bP^x$-a.s.
  \begin{equation}\label{e:2.17a}
  u(X_{t\wedge \tau_U}) = \lim_{n\to \infty} \bE^x \left[ u(X_{n\wedge \tau_U}) | \sF_{t\wedge \tau_U} \right]= 
  \lim_{n\to \infty} \bE^x \left[ u(X_{n\wedge \tau_U}) | \sF_{t } \right]
  = \bE^x \left[ \xi | \sF_{t } \right] . 
  \end{equation} 
  Since $\xi$ is $\sF_\infty$-measurable, and the filtration $\{\sF; t\geq 0\}$ is right continuous, 
  there exists a right continuous martingale $(M_t)$ such that $M_t$ is a version of
 $\bE^x(\xi| \sF_t)$ for each $t$. Using \eqref{e:2.16} and the martingale convergence theorem, 
 $M_t$  converges to 
 $\bE^x \left[ \xi | \sF_\infty \right] = \xi$ both $\bP$-a.s. and in $L^2(\bP^x)$  as $t\to \infty$. 
 Hence for every   $x\in U \setminus \sN$ and $t\geq 0$, we have by  \eqref{e:2.17a} and \eqref{e:2.16} that $\bP^x$-a.s. on $\{t<\tau_U\}$, 
\begin{equation}\label{e:2.18}
\xi \circ \theta_t =  \left( \lim_{n \to \infty }   \bE^x \left[ \xi | \sF_{n-t } \right] \right) \circ \theta_t
= \lim_{n\to \infty  }   u(X_{(n-t)\wedge \tau_U}) \circ \theta_t 
= \lim_{n\to \infty }   u(X_{n\wedge \tau_U}) =\xi.
\end{equation} 
 By \eqref{e:2.16}, \eqref{e:2.17a} and \eqref{e:2.18}, the random variable $\xi$ has the property that  
  \begin{align}
\label{e:hrep}	 u(x)&= \lim_{t\to \infty}  \bE^x[   u(X_{t \wedge \tau_U } ) ]
 =\bE^x [ \xi] \quad \mbox{for $x\in U\setminus \sN$,} \\
\label{e:ter}	  
 \bP^x( &\xi \neq  \xi\circ \theta_{t}, t \le \tau_U ) =0   \quad  \hbox{for every  } x \in U \setminus \sN. 
\end{align}
 Moreover,  it follows from \eqref{e:2.16} that 
\begin{equation} 
\bP^x  ( \xi   = u(X_{\tau_U})  \hbox{ on } \{\tau_U <\infty\} )=1 \quad \mbox{for } x\in U\setminus \sN.
\end{equation}
 We can then use the same proof as that for \cite[Theorem 2.9]{C}  to conclude that $u$ is $\sE$-harmonic in $U$ and hence in $D$. 
From \eqref{e:hrep}-\eqref{e:ter}  and \cite[Theorem 6.1.6]{CF}, we conclude that $u$ is $\sE$-quasi-continuous on $D$.   
\qed

\smallskip 

 As noted in \eqref{e:2.4}, condition 
\eqref{e:exit} is satisfied if $(\sE, \sF)$ is irreducible and $U^c$ is non-$\sE$-polar.  
It follows from Proposition \ref{P:2.9} that a locally bounded universally measurable function $u$ in $D$
is harmonic in $D$ if and only if it is  $\sE$-quasi-continuous on $D$ and $\sE$-harmonic in $D$.  

\begin{prop}\label{P:2.10} Suppose $D$ is an open subset of $\sX$.
\begin{enumerate} [\rm (i)]
\item    Any non-negative regular harmonic function $u$ in $D$ is the   $\sE$-q.e.   limit  on $\overline{D}$ of an increasing sequence of bounded non-negative regular harmonic functions $\{u_n; n\geq 1\}$ in $D$.

\item Suppose that $u$ is a  non-negative harmonic function in $D$ and 
that $U$ is a  relatively compact open subset of $ D$.
Then $\restr{u}{U}$  is the increasing limit $\sE$-q.e.~on $U$
of a sequence of bounded harmonic functions on $U$. 
\end{enumerate} 
\end{prop}

\pf  (i) For $n\geq 1$, define $u_n(x)= \bE^x \left[ (n\wedge u)(X_{\tau_D})\right]$ for $x\in \ol D\setminus
\sN_0$. Clearly, by \eqref{e:2.8} and the fact that $\mathbb{P}^x(\tau_D=0)=1$ for  $\sE$-q.e.~$x \in \partial D$ by \cite[Theorem A.2.6(i), 4.1.3 and 4.2.1(ii)]{FOT},  $u_n$ is a bounded non-negative regular harmonic function in $D$ that increases to
$u(x)=\bE^x \left[ u(X_{\tau_D})\right]$ for    $x\in \ol D\setminus \sN_0$, as $n\to \infty$.

(ii) Let $U$ be any relatively compact open subset of $D$. 
 By definition, there is a property exceptional set $\sN\subset U$ so that 
$t\mapsto u(X_{t\wedge \tau_U})$ is a $\bP^x$-uniformly integrable martingale for every $x\in U\setminus \sN$. 
As noted in the proof for Proposition \ref{P:2.9}(iii), 
from the martingale theory, we know that for $x\in U\setminus \sN$, $\bP^x$-a.s,
\begin{equation}\label{e:2.17}
M_t:=\lim_{s\in \bQ: s\to t+} u(X_{s\wedge \tau_U}) \quad \hbox{exists for every } t\geq 0
\end{equation}
 and $M_t$ is a right-continuous modification 
of $u(X_{t\wedge \tau_D})$. 
 By the martingale convergence theorem.
 $M_t$ converges  to some random variable $\xi\geq 0$   in $L^1(\bP^x)$ and $\bP^x$-a.s. for every  $x\in U \setminus \sN$
as $t\to \infty$.  Note that in view of \eqref{e:2.16}, $\xi$ has the property that for every $x\in U\setminus \sN$, 
$u(x)= \bE^x [ \xi]$  and $\bP^x$-a.s.
 $\xi= \xi\circ \theta_{t\wedge \tau_U}$ for every $t\in [0, \tau_U)$. 
For each $n\geq 1$, define $u_n(x)= \bE^x [\xi \wedge n]$. 
Then $u_n$ is bounded and harmonic in $U$ and it increases to $u$ on $U\setminus \sN$ as $n\to \infty$.
 \qed

	\section{ Local regularity  for transient spaces} \label{S:lreg}

		Since $(\sE, \sF)$ is strongly local,  by Proposition \ref{P:2.2}, its corresponding Hunt process 
		$X$ is a diffusion that admits no   killing 
		inside $\sX$. 
	Thus there exists  a  Borel properly exceptional set $\sN_0$ so that   the Hunt process $X$, whose lifetime is denoted by $\zeta$, 
	can start from every point in $\sX\setminus \sN_0$ and 
 \eqref{e:2.8} holds. 
 	It follows that
	\begin{equation}\label{e:2}
	\bP^x (X_t=x \hbox{ for all } t\in [0, \zeta) \hbox{ and } \zeta<\infty)=0
	\qquad \hbox{for every } x\in \sX \setminus \sN_0.
	\end{equation}

	  In the remainder of this  section  unless otherwise specified, 
	we  assume in addition that $(\sE, \sF)$ is transient. 
	In view of  \cite[Theorem 3.5.2]{CF}, by enlarging the Borel properly exceptional set $\sN_0$ if needed, we may and do assume that
	\begin{equation}\label{e:3a} 
	\bP^x  \left(\zeta =\infty \hbox{ and } \lim_{t\to \infty} X_t =\partial \right) = \bP^x  \left(\zeta = \infty \right) 
	\quad \hbox{for every } x\in \sX\setminus \sN_0.
	\end{equation}

		\begin{prop} \label{P:3.1}
	For any relatively compact open subset $D\subset \sX$, 
	$$
	\bP^x (\tau_D <\zeta)=1  \quad \hbox{for every } x\in D \setminus \sN_0.
	$$
	\end{prop} 
	
	\pf   	Recall that $\sX_\partial=\sX \cup \{\partial\}$ is the one-point compactification of $\sX$.  Since $D$ is a
	relatively compact open subset of $\sX$, $\sX \cup \{\partial \}\setminus \overline D$
	is an open neighborhood of $\partial$.
	We know from \eqref{e:2.8} that  
	$$
	\bP^x( X_{\zeta-} =\partial, \, \zeta <\infty) =\bP^x (\zeta <\infty)
	\quad \hbox{for } x\in \sX\setminus \sN_0.
	$$
	 This together with \eqref{e:3a} implies that 
	$$ \bP^x  \left(  \lim_{t\uparrow \zeta} X_t =\partial \right) =1
		\quad \hbox{ for any } x\in \sX\setminus \sN_0.
	$$ 
	Consequently, 
	we have that $\bP^x(\tau_D < \zeta)=1$ for every $x\in D\setminus \sN_0$. 
	\qed 
 
 The transience condition on $(\sE, \sF)$ can be dropped from Proposition \ref{P:3.1} 
 if we assume   $(\sE, \sF)$ is irreducible and $D^c$ is of positive $\sE$-capacity.  We emphasize that  the
 transience of $(\sE, \sF)$ is not assumed in the next Proposition.

\begin{prop}	\label{P:3.2}
		Suppose that  $(\sE, \sF)$ is a symmetric,  irreducible, strongly local regular Dirichlet form
		on $L^2(\sX; m)$ and that $D\subset \sX$ is a relatively compact open set of $\sX$ so that 
		$D^c$ is not $\sE$-polar. Then $\bP^x(\tau_D <\infty ) =1$ for  $\sE$-q.e. $x\in D$.
	\end{prop}
	
	\pf Since $(\sE, \sF)$ is an irreducible,  it is either transient or recurrent.
	The conclusion of the proposition follows readily from Proposition \ref{P:3.1} if
	$(\sE, \sF)$ is transient.
	When $(\sE, \sF)$ is recurrent, the desired conclusion holds as well since 
	$X_t$ in fact visits $\sX\setminus D$ infinitely often under $\bP^x$
	for  $\sE$-q.e. $x\in \sX$ by \cite[Theorem 3.5.6]{CF}. 
	\qed

	\medskip

	\begin{lemma}\label{L:1}
	 		$\bP^x (\tau_x >0)=0$ 
 	for every $x\in \sX\setminus \sN_0$,
	  where $\tau_x:= \inf\{t\geq 0: X_t\not= x\}$.  
	 \end{lemma}

	\pf 
	We have by   
	\eqref{e:2} and \eqref{e:3a} 
	that 
	$\tau_x<\zeta $ $\bP^x$-a.s. for every $x\in \sX \setminus \sN_0$.  
	Clearly  $X_{\tau_x}=x$ on $\{0<\tau_x<\zeta\}$ since $X$ is a diffusion. 
	Let  $A_x:=\{\tau_x>0\}$.
	Since $\mathbf{1}_{A_x} = \mathbf{1}_{A_x^c} \circ \theta_{\tau_x}$ on $\{0<\tau_x<\zeta \}$, we have by the strong Markov property of $X$ that
	for $x\in \sX \setminus \sN_0$, 
	$$ 
	\bP^x (A_x)=\bE^x \left[ \bP^{X_{\tau_x}} (A_x^c); 0<\tau_x<\zeta \right] = \bP^x (A_x^c) \,\bP^x(0<\tau_x<\zeta)
	=(1-\bP^x (A_x))\bP^x(A_x).
	$$
	It follows that $\bP^x(A_x)=0$. 
	\qed
	
	\bigskip

 Let $\sB_+(\sX)$ denote the set of non-negative Borel measurable functions on $\sX$.	Denote by  $\{P_t; t\geq 0\}$  the transition semigroup of  the process $X$; that is,  
	$$
	P_t f(x)=\bE^x [ f(X_t)],  \qquad x\in \sX \setminus \sN_0,  \, t>0 ,  \, f\in {\cal B}_+ (\sX),
	$$  
	with the convention that  $f(\partial ):=0$.  
	Define the Green operator $G$ by
	$$
	Gf(x):=\bE^x \int_0^\infty f(X_t) dt= \int_0^\infty \bE^x [f(X_t)] dt = \int_0^\infty P_t f (x) dt, \quad x\in \sX \setminus \sN_0, \, 
	f\in {\cal B}_+(\sX).
	$$

\begin{lemma} \label{L:g0}
By enlarging the Borel properly exceptional set $\sN_0$ if necessary,
	there is 
 	a function $g_0 \in  L^1(\sX; m)$   that takes values in $(0, 1]$ on $\sX$ 
	such that 
	\begin{equation}\label{e:4a}
	G g_0(x) \leq 1   \hbox{  for } x\in \sX \setminus \sN_0,  
 	\quad  Gg_0\in \sF_e  \  \  \hbox{ and }  \ \ \sE (Gg_0, Gg_0) \leq 1.
	\end{equation}
	 \end{lemma} 

	\proof
 By \cite[Theorems 2.1.5(i), Lemma 2.1.4(ii) and Theorem 2.1.12(i)]{CF},  there is an $L^1(\sX; m)$ function $g_1$ bounded by 1,
	strictly positive  on $\sX$, such that 
	$Gg_1<\infty$ $m$-a.e. on $\sX$ and  
	 	$Gg_1\in \sF_e$ with $\sE (Gg_1, Gg_1)\leq 1$. 
	Since $Gg_1$ is excessive and hence finely continuous \cite[Theorem A.2.2]{CF},  by enlarging the properly exceptional set 
	$\sN_0$ if necessary, we may and do assume
	that $Gg_1(x)<\infty$ for every $x\in \sX \setminus \sN_0$, thanks to \cite[Theorem A.2.13(v)]{CF}. 
	
 	Let $k \ge 1$ and 
	set $f_k = g_1 1_{\{Gg_1\leq k\}}$ and $T_k = \inf\{ t \ge 0: Gg_1(X_t) \le k \}$.
	Since $Gg_1$ is finely continuous we have $G f_k(X_{T_k})\le Gg_1(X_{T_k})\le k$  $\bP^x$-a.s.~on ${\{ T_k<\infty \}}$ for every 
	$x \in \sX \setminus \sN_0$.
	Hence by the strong Markov property,  
	$$ Gf_k(x)  = \bE^x 1_{\{ T_k<\infty \}}  \int_{T_k}^\infty f_k(X_s)ds = \bE^x 1_{\{ T_k<\infty \}} Gf_k(X_{T_k}) \le k
	\quad \hbox{for } x \in \sX \setminus \sN_0 . 
	$$
 	Let $  g_0=\sum_{k=1}^\infty k^{-1}2^{-k} f_k + 4^{-1} 1_{\sN_0} $. 
	Then $\  g_0$ is $(0,1]$-valued on $\sX$,  and 
	$$
	G  g_0 (x) \leq  1 \quad \hbox{for every } x\in \sX \setminus \sN_0, 
	$$  
	$G   g_0\in \sF_e$ and $\sE (G  g_0, G  g_0)\leq \sE (G g_1, G g_1)\leq 1$  
	in view of \cite[Theorem 2.1.12(i)]{CF}.  \qed

	It follows from \eqref{e:4a} that for every $x\in \sX \setminus \sN_0$,
	$G(x, dy)$, defined by $G(x, A)=G1_A (x)$, is a $\sigma$-finite  Borel measure on $\sX$.
	By the symmetry of the process $X$, each $P_t$ is a symmetric operator in $L^2(\sX; m)$. Hence  
	\begin{equation}\label{e:3} 
	\int_\sX g(x) Gf(x) m(dx)= \int_\sX f(x) G g(x) m(dx) \qquad \hbox{for } f, g\in {\cal B}_+(\sX).
	\end{equation}

	\medskip

	We introduce a function  $d_{\sX}$  on a locally compact separable metric space
	 	$(\sX, d)$.   Define for every $x\in \sX$, with the convention   $\inf \emptyset := \infty$,
\begin{equation} \label{e:3.5}
 	d_{\sX} (x)= \sup \{r>0: B(x, r) \hbox{ is relatively compact in } (\sX, d) \}.
 	\end{equation}

 	\begin{remark}
	\begin{enumerate}  
	\item[\rm (i)]   
	   Let $\{D_k; k\geq 1\}$ be an increasing sequence of relatively compact open subsets with $\cup_{k=1}^\infty D_k=\sX$. 
	It is easy to see that  for $x\in \sX$, 
	\begin{equation} \label{e:4.1}
 	d_{\sX} (x)=  \lim_{k\to \infty}  \inf_{y\in \sX \setminus D_k} d(x, y).  
 	\end{equation}

	\item[\rm (ii)]  Either $d_{\sX}$ is identically infinite on $\sX$ or $d_{\sX}(x)<\infty$ for every $x\in \sX$ and 
	$|d_{\sX}(x) - d_{\sX}(y)| \leq d(x, y)$ for every $x, y\in \sX$.

	 \item[\rm (iii)]    $d_\sX \equiv \infty$ on $\sX$ if and only if $B(x,r)$ is relatively compact for every $x \in \sX$ and every $r>0$.  
	 \end{enumerate} 
	\end{remark}

	 	We call $d_{\sX}$ the \emph{distance to the boundary function} for the metric space $(\sX, d)$. 
	  For an open subset $D\subset \sX$, we can similarly define $d_D$ on $(D, d|_{D\times D})$ with the locally compact separable 
	metric subspace $(D, d|_{D\times D})$ in place of $(\sX, d)$; that is, 
	\begin{equation} \label{e:4.3}
 	d_D (x)= \sup \{r>0: B(x, r) \hbox{ is relatively compact in } (D, d_{D\times D}) \}. 
 	\end{equation}
	 Clearly, we have
	\begin{equation}
	d_D (x) \leq d_{\sX}(x) \quad \hbox{for any } x\in D.
	\end{equation}

	\begin{definition} {\rm
		We say condition {\rm (HC)} holds if there is an $\sE$-nest $\{F_n; n\geq 1\}$
		consisting of an increasing sequence of compact subsets with $\sN_0\subset \sX \setminus \cup_n F_n$
		such that for all  $x_0 \in \sX$ 
		  there exists $r_{x_0} \in (0,   d_{\sX} (x) \wedge 1)$, 
		such that if $r \in (0, r_{x_0})$} and $f \in \mathcal{B}_+(\sX)$ 
		has compact support in $B(x_0, 2r)^c$, and satisfies 
	    $0\le f \leq cg_0$ for some $c>0$, then 
		$Gf(x)$ is  continuous  in $ B(x_0, r) \cap F_n$ for every $n\geq 1$.   
	\end{definition}

 	Condition (HC) is a weaker condition than the condition that bounded harmonic functions
	in $B(x_0, 2r)$ are continuous in $B(x_0, r)$. We show in Example \ref{E:8.3} an example
	of a MMD space $(\sX, d, \sE, \sF, m)$ 
	such that (HC) holds but there is a bounded harmonic function in $B(x_0, 2r)$ 
	that is not continuous at $x_0$  and in Example \ref{E:8.cylinder} we 
	use the infinite product spaces studied in \cite{BSC} to give an example of a space which does not satisfy (HC). 	 
	
	 	Note that it follows from the   definition of $\sE$-nest in Section \ref{S:2}, 
		if $\{F_n; n\geq 1\}$ is an $\sE$-nest, then so is $\{ K_n; n\geq 1\}$, 
	where $K_n={\rm supp} [1_{F_n}m]$.
	Thus without loss of generality, in this paper we always assume that the $\sE$-nest in {\rm (HC)} has the property
	that $F_n= {\rm supp} [1_{F_n}m]$ for every $n\geq 1$.  For an $\sE$-nest $\{F_n\}$, 
	$\sN=\sX \setminus \cup_n F_n$ is $\sE$-polar  and, in particular,   has zero $m$-measure. 
	
	\begin{thm} \label{T:1} 
		Assume that condition {\rm (HC)} holds with an $\sE$-nest $\{F_n\}$.
	 		Let $\sN$ be a Borel properly exceptional set that contains $\sX \setminus \cup_n F_n\supset \sN_0$. 
		Then for every $x\in \sX \setminus \sN$, $G(x, dy)$ is absolutely continuous with respect to $m$.
		Consequently,  for every $x\in \sX \setminus \sN$ and $t>0$,
		$P_t(x, dy):=\bP^x (X_t \in dy)$ is absolutely continuous with respect to $m$. 
	\end{thm}
	
	\pf     It follows from \eqref{e:3} that $G(x, A)=0$ $m$-a.e. on $\sX$ for every 
	$A\subset \sX$ with $m(A)=0$ (by taking $f=1_A$ and $g=1$). 
	Let $g_0$ be the strictly positive function from Lemma \ref{L:g0}. 
	 	Fix $x_0\in \sX\setminus \sN$ 
	and   $r\in (0,r_{x_0})$. 
	For $j\geq 1$, let $E_j=\{x\in \sX: 2^{-j}<  g_0(x) \leq 2^{1-j}\}$. 
	Then the $E_j$ form a partition of $\sX$. 
	Let  $A\subset B(x_0, 2r)^c$ 
	be a Borel set with compact closure in $B(x_0,2r)^c$ and satisfy $m(A)=0$. Since $1_{A\cap E_j} \leq 2^j g_0$, 
	we have  by condition {\rm (HC)} that for each $k \ge 1$   
	the function $x\mapsto G(x, A  \cap E_j)$ is continuous 
	and therefore zero on $B(x_0, r)\cap F_k$.
	Thus $G(x, A  \cap E_j)=0$ for every $x\in B(x_0, r)\setminus \sN$.  
	Consequently, $G(x, A)=\sum_{j=1}^\infty G(x, A\cap E_j)=0$ for every $x\in B(x_0, r)\setminus \sN$. 
	In particular, this shows that  for every $x_0\in \sX \setminus \sN$,
	\begin{equation}\label{ e:4}
	G(x_0, dy) \hbox{  is absolutely continuous with respect to } m(dy) \hbox{ on }
	\sX \setminus \{x_0\}. 
	\end{equation}
	
	 	We claim that $G(x_0, dy)$  is absolutely continuous with respect to  $ m(dy)$  on $\sX$.
	This is clearly true if  $m(\{x_0\})>0$.  We thus assume $m(\{x_0\})=0$ and for $t \ge 0$ define
	$$  h(x,t):= \bE^x \int_t^\zeta 1_{\{x_0\}} (X_s) \,ds.$$
We set $h(x)=h(x,0) =  (G 1_{\{x_0\}})(x)$, and need to  prove that $h(x_0)=0$.

The function $h$ is   harmonic on $\sX \setminus \{x_0\}$	and since 
$m(\{ x_ 0\})=0$, we have by \eqref{e:3} that  $h=0$  $m$-a.e. on $\sX$.
Further, by condition {\rm (HC)},  $h(x)=0$ on  $\sX \setminus (\sN \cup \{x_0\})$.
	Thus if $A=\{  y: h(y)>0\}$ then $A \subset \sN \cup \{x_0\}$. 
Let $T= \inf\{ s \in \bQ \cap [0, \infty): X_s \not\in A\}$; as $\sN$ is properly exceptional and 
(by Lemma \ref{L:1}) $X$ leaves $x_0$ immediately, we have $\bP^{x_0}(T=0)=1$.

 	Let $ \{  {\cal F}_t  \}_{t\geq 0}$ be the  minimum augmented admissible filtration generated 
	by $\{ X_t \}_{t \geq  0}$. 
	Let $t>0$, and set $M_s = h(X_s, t-s)$ for $s \in [0,t]$.   
By the Markov property of $X$, for $s\in [0, t]$, 
$$
	\bE^{x_0}   \left[ \int_t^\zeta 1 _{\{x_0\} } (X_r ) dr   \big| {\cal F}_s \right]
	=\bE^{X_s} \int_{t-s}^\zeta 1_{\{x_0\} } (X_r ) d r = h(X_s, t-s) =M_s  
	\quad   \bP^{x_0} \hbox{-a.s.}.
	$$
  Thus $\{  M_s; s\in [0, t]\}$ is a  non-negative $\bP^{x_0}$-martingale that has  a right continuity modification (cf. \cite[Theorem 1.3.13]{KS}). 
In particular, there is $\Omega_0\subset \Omega$  with $\bP^{x_0}  (\Omega_0)=1$ so that 
$  M_s (\omega) $  is a right continuous in $s\in [0, t] \cap \bQ$ for every $\omega \in \Omega_0$.
For every  $\omega \in \{ T=0\}\cap \Omega_0$,  
   there exists a sequence $s_n \downarrow 0, s_n \in   [0, t] \cap \bQ$
such that $X_{s_n}(\omega) \not\in A$, and thus 
$$
0\leq M_{s_n}(\omega)=h(X_{s_n}(\omega), t-s_n) \leq h(X_{s_n}(\omega)) =0.
$$
 Consequently we have $M_0(\omega)=0$ $\bP^{x_0}$-a.s. 
Thus we have $h(x_0,t) =\bE^{x_0} M_0=0$.
It follows then   $h(x_0)=\lim_{t\to 0} h(x_0, t)= 0$.
   	This together with \eqref{ e:4} shows that  $G(x, dy)$ is absolutely continuous with respect 
	to $m(dy)$ on $\sX$ for every $x\in \sX \setminus \sN$. 
	That  $\bP^x (X_t \in dy)$ is absolutely continuous with respect to $m$ for every $x\in \sX \setminus \sN$ and $t>0$
	follows immediately from  \cite[Proposition 3.1.11]{CF} or \cite[Theorem  4.2.4]{FOT}.
	\qed

	 	With Theorem \ref{T:1} at hand, we can deduce the following. 
 	
	\begin{thm}\label{T:2}
		Assume that  condition {\rm (HC)} holds with  an
		$\sE$-nest $\{F_n\}$. 
		Let $\sN$ be a Borel properly exceptional set that contains 
		$\sX \setminus \cup_n F_n\supset\sN_0$. Then
		there exists  a non-negative
		jointly   ${\cal B}(0, \infty)\times {\cal B}(\sX \times \sX)$-measurable
		function $p(t, x, y)$ on $(0, \infty) \times ( \sX \setminus \sN) \times (\sX \setminus \sN)$ such that 
		\begin{description}
			\item{\rm (i)} for every $f\in {\cal B}_+ (\sX)$,  $x\in \sX \setminus \sN$ and $t>0$,
			$\displaystyle \bE^x f(X_t)=\int_\sX p(t, x, y)  f(y) m(dy)$;
			
			\item{\rm (ii)} $p(t, x, y)=p(t, y, x)$ for every $x, y\in \sX \setminus \sN$ and $t>0$; 
			
			\item{\rm (iii)} For every $t, s>0$ and $x, y\in \sX \setminus \sN$, $\displaystyle p(t+s, x, y) = \int_\sX p(t, x, z)p(s, z, y) m(dz)$.
		\end{description}
		Consequently, $g(x, y):=\int_0^\infty p(t, x, y) dt$, $x, y\in \sX \setminus \sN$, 
		is  a non-negative
		jointly   ${\cal B}(\sX \times \sX)$-measurable
		function   on  $( \sX \setminus \sN) \times (\sX \setminus \sN)$ such that 
		\begin{description}
			\item{\rm (iv)} $\displaystyle Gf(x)=\int_\sX g(x, y) f(y) m(dy)$ for every $x\in \sX \setminus \sN$ and $f\in {\cal B}_+ (\sX)$;
			
			\item{\rm (v)} $g(x, y)=g(y, x)$ for every $x, y\in \sX \setminus \sN$, and $x\mapsto g(x, y)$ is
			$X|_{\sX\setminus \sN}$-excessive and finite  $\sE$-q.e.~on $\sX$ for  every $y\in  \sX \setminus \sN$.
			
			\item{\rm (vi)} For every $y_0\in  \sX \setminus \sN$, $x\mapsto g(x, y_0)$    is harmonic in 
		 			$ \sX \setminus \{y_0\}$. 
 In fact, for any open subset $U$ of $\sX$ with $y_0 \notin \overline{U}$, the function
					$x\mapsto g(x, y_0)$ is regular harmonic in $U$.   
 		\end{description}
	\end{thm}
	
	\pf   We first show that for each $x\in \sX \setminus \sN$ and $t>0$, $X$ has a pointwisely defined  transition density function $p(t, x, y)$.
	This part is almost the same as that for  \cite[Theorem 3.1]{BBCK}. 
	For the reader's convenience, we spell out the details here.
	
	By Theorem \ref{T:1}, 
	for every $t>0$ and $x\in \sX \setminus \sN$ there is a $[0,\infty]$-valued integrable kernel $y\mapsto
	p_0(t, x, y)$ defined  
	 on $\sX$ such that
	\begin{equation}\label{eqn:heat2}
	\bE^x \left[ f(X_t)\right]=P_t f(x) = \int_{\sX} p_0(t, x, y) f(y)\, m(dy)
	\qquad \hbox{for every } f\in {\cal B}_b(\sX), 
	\end{equation}
	where ${\cal B}_b(\sX)$ denotes the set of bounded Borel measurable functions on $\sX$.
    We note that $p_0(t, x, y)$ can be chosen to be jointly Borel measurable on $(0, \infty) \times \sX \times \sX$ by an application of the martingale convergence theorem (see, e.g., \cite[Proposition 5.6]{GK} with $H_t(x,y) \equiv \infty$).
	From the semigroup property $P_{t+s}=P_tP_s$, we have for every $t,
	s>0$ and $x\in \sX \setminus \sN$,
	\begin{equation}\label{eqn:heat3}
	p_0(t+s, x, y)=\int_\sX p_0(t, x, z) p_0(s, z, y) m(dz) \qquad
	\hbox{for } m \hbox{-a.e. } y\in \sX.
	\end{equation}
	Note that since $P_t$ is symmetric, we have for each fixed
	$t>0$,
	\begin{equation}\label{e:7} 
	p_0(t, x, y)=p_0(t, y, x) \qquad \hbox{for $m$-a.e. }
	(x, y)\in \sX \times \sX .
	\end{equation}
	
	For every
	$t>0$ and $x, y \in \sX \setminus \sN $, let $s\in (0,  t/3)$ and define
	\begin{equation}\label{eqn:newp}
	p(t, x, y):= \int_{\sX} p_0(s, x, w)
	\left( \int_{\sX} p_0(t-2s, w, z)p_0(s, y, z) m(dz) \right) m(dw).
	\end{equation}
    Clearly, $p(t, x, y)$ is jointly Borel measurable on $(0, \infty) \times (\sX \setminus \sN)
	\times (\sX \setminus \sN)$.
	By \eqref{eqn:heat3} and \eqref{e:7},  the above definition is independent
	of the choice of $s\in   (0, t/3)$.
	Clearly by \eqref{e:7} with $t-2s$ in place of $t$ and $(w, z)$ in place of $(x, y)$, we see that 
	\begin{equation}\label{e:C7}
	p(t, x, y)=p(t, y, x) \quad  \hbox{for every } x, y\in \sX \setminus \sN.
	\end{equation}
	By the semigroup property  \eqref{eqn:heat3},  \eqref{eqn:heat2} and \eqref{e:7}, we have for any $\phi \geq
	0$ on $\sX$ and $x\in \sX \setminus \sN$,
	\begin{align}
	& \bE^x \left[ \phi (X_t)\right]  \nonumber \\
	&= \int_{\sX} \left(
	\int_{\sX} p_0(s, x, w) p_0(t-s, w,  y) m( dw )\right)
	\phi (y) m(dy)  \nonumber \\
	&= \int_{\sX} \left(
	\int_{\sX} p_0(s, x, w) \left( \int_{\sX} p_0(t-2s, w, z) p_0(s, y, z) m(dz) \right)
	m(dw) \right) \phi(y) m(dy ) \nonumber \\ 
	&= \int_{\sX} p(t, x, y) \phi(y) m(dy). \label{e:C8}
	\end{align}
	Thus for each $x\in \sX \setminus \sN$,  $p(t, x, y)$ coincides with  $p_0(t, x, y)$ for $m$-a.e. $y\in \sX$. 
	For $t, s>0$ and $x, y\in \sX \setminus \sN$, take $s_0\in   (0,
	\,  (t\wedge s)/3)$. We have by (\ref{eqn:heat3})-(\ref{eqn:newp})
	\begin{align}
	& p(t+s, x, y)  \nonumber \\
	&= \int_{\sX} p_0(s_0, x, w) \left( \int_{\sX}
	p_0(t+s-2s_0, w, z) p_0(s_0, y, z) m(dz) \right) m(dw)   \nonumber \\
	&=\int_{\sX^5} p_0(s_0, x, w)  p_0(t-2s_0, w, u_1) p_0(s_0, u_1, u_2)
	p_0(s_0, u_2, v) p_0(s-2s_0, v, z)   \nonumber \\
	& \q \q  \q p_0(s_0, y, z) m(dw) m(du_1) m(du_2) m(dz) m(dv)   \nonumber \\
	&= \int_\sX p(t, x, v) p(s, v, y) m(dv) . \label{e:C9}
	\end{align}

	Define $g(x, y):=\int_0^\infty p(t, x, y) dt$ for $x, y\in \sX \setminus \sN$.
    Note that $g: (\sX \setminus \sN) \times (\sX \setminus \sN) \to [0, \infty]$
	is jointly Borel measurable. 
	It follows from \eqref{e:C8} and Fubini's theorem that for every $f\in {\cal B}_+(\sX)$,
	$$ 
	G f(x):= \bE^x \int_0^\infty f(X_s) ds = \int_\sX g(x, y) f(y) m(dy) \quad \hbox{for every } x\in \sX \setminus \sN.
	$$
	Clearly by \eqref{e:C7}, $g(x, y)=g(y, x)$ for every $x, y\in \sX \setminus \sN$.
	Note that for each fixed $y\in \sX \setminus \sN$ and $t>0$, by Fubini's theorem and \eqref{e:C9}, 
	$$
	P_t  g(\cdot , y) (x) =  \int_\sX p(t, x,  z) g(z, y) \,m(dz)=\int_t^\infty p ( s, x, y)\, ds \leq g(x, y)
	\quad \hbox{for every } x\in \sX  \setminus \sN,
	$$
 and hence $t \mapsto P_t  g(\cdot , y) (x) $ is non-increasing and $\lim_{t\downarrow 0} P_t g(\cdot , y) (x) =g(x, y)$. This shows that for each fixed $y\in \sX \setminus \sN$, $x\mapsto g(x, y)$ is an excessive
	function of $X|_{\sX\setminus \sN}$. Moreover, for each $x \in \sX \setminus \sN$, since $\int_{\sX} g(x,y) g_0(y) \,m(dy) = Gg_0(x) \le 1$, $g(x,\cdot)<\infty$ $m$-a.e.  and hence $g(x,\cdot)<\infty$ $\sE$-q.e. on $\sX$  by  \cite[Theorem A.2.13]{CF}.
	
	\smallskip
	
	The proof of  (vi) is similar to that for \cite[Proposition 6.2]{KW}.
	Let $y_0\in \sX \setminus \sN$ . Let $U$ be an open subset of $\sX$ with $y_0 \notin \overline{U}$, and take $t_1>0$ so that $\overline{U} \cap B(y_0,r_1) = \emptyset$.
		For any  non-negative $f\in\sB_+(X)$, by Fubini's theorem and 
	the strong Markov property of $X$,
	for each $x\in U \setminus \sN$ with $\restr{f}{B(y_0,r_1)^c}\equiv 0$, 
	\begin{align} \label{e:n1}
		\int_{B(y_0, r_1)} \bE^x \left[ g(X_{\tau_{U}}, y) \right] f (y) m(dy)
		&= \bE^x  \left[(Gf)(X_{\tau_{U}}) \right]  \cr
		&=Gf(x) 
		 =  \int_{B(y_0, r_1)} g(x, y) f(y) m(dy),
	\end{align}
	and this is finite if $0 \le f \le c g_0$ on $\sX$ for some $c>0$.
	Hence for each $x\in U\setminus \sN$,  
	$$
	\bE^x \left[ g(X_{\tau_{U}}, y) \right] = g(x, y)  \quad \hbox{for } m \hbox{-a.e.  } y\in B(y_0, r_1) .
	$$
	Since $y\mapsto g(z, y)$ is $X|_{\sX\setminus \sN}$-excessive, it follows from the  monotone convergence theorem, Fubini theorem,
	the fine continuity of $y\mapsto g(z, y)$ (see \cite[Theorem A.2.2]{CF} or \cite[Theorem A.2.5]{FOT}) and Fatou's lemma that
	\begin{eqnarray*}
		\bE^x \left[ g(X_{\tau_{U}}, y_0) \right] 
		&=& \lim_{t\downarrow 0}   P_t  \left(\bE^x   g(X_{\tau_{U}}, \cdot) \right)(y_0) \\
		&\geq& \limsup_{t\downarrow 0}    P_t \left(  g(x, \cdot) 1_{B(y_0, r_1)} \right) (y_0) \\
		&\geq & \bE^{y_0} \left[ \liminf_{t\downarrow 0} g(x, X_t)  1_{B(y_0, r_1)} (X_t) \right] \\ 
		&=& g(x, y_0).
	\end{eqnarray*}
	On the other hand, the optional sampling theorem (see, e.g., \cite[Theorem 1.3.22]{KS}) implies that $\bE^x \left[ g(X_{\tau_{U}}, y_0) \right] \leq g(x, y_0)$ as $x\mapsto g(x, y_0)$
	is excessive for $X|_{\sX\setminus \sN}$ and hence $\{g(X_t,y_0)\}_{t \ge 0}$ is a right-continuous $\bP^x$-supermartingale for any $x \in \sX \setminus \sN$ with $g(x,y_0)<\infty$ by \cite[Proof of Theorem A.2.2]{CF}.  Thus we have $  \bE^x \left[ g(X_{\tau_{U}}, y_0) \right] = g(x, y_0)$ for every $x\in U \setminus \sN$; that is, $x\mapsto g(x, y_0)$ is regular harmonic in $U$.
	This in particular proves that $x\mapsto g(x, y_0)$ is harmonic in $\sX \setminus \{y_0\}$. 
	\qed

   We call the function $g(x, y)$ in Theorem \ref{T:2} the \emph{Green function} of $X$.

\begin{remark} {\rm
There are gaps in the proofs of  the existence of a Green function
in \cite{BBKu} and \cite[Lemma 5.2]{GH}.
For details of the gap in \cite{GH},  see \cite[Remark 4.19]{BM2}.
The gap in [BBK] is that it is not proven that the Green's function is the  integral kernel of the 
Green operator  (cf. Theorem \ref{T:2}\,(iv)). 
}\end{remark}

We next give a sufficient condition for {\rm (HC)}. 
  The following definition is valid for any regular Dirichlet form on $\sX$, that is, for that definition, 
we do not  assume  $(\sE, \sF)$ is transient.  

   	\begin{defn} \label{D:3.8} \rm  
  We say that the  {\em (non-scale-invariant)  elliptic Harnack inequality} {\rm (Ha)} holds on $\sX$ if 
  for 	any ball $B=B(x_0, r)$ in $\sX$  with $r\in (0, d_{\sX}(x_0) \wedge 1)$, 
    there are constants $C_B>1$ and $\delta_B \in (0, 1)$ such that  
		 			for any non-negative 
			bounded $u \in \sF_e$ that is regular harmonic 
			in $B(x_0, r)$, 
	\be  \label{e:3.14}
	 {\rm esssup}_{B(x_0, \delta_B  r)} u \leq C_B \, {\rm essinf}_{B(x_0, \delta_B r)} u. 
	 \ee
	\end{defn}

\ms

	\begin{remark}\label{R:3.10} \rm
	 	\begin{enumerate} [\rm (i)]
 \item   If (Ha) holds, then \eqref{e:3.14}  holds for any non-negative
			  $u\in \sF_e$ that is regular harmonic in $B(x_0, r)$.  This is because 
				$u(x) = \bE^x \left[ u (X_{\tau_{B(x_0, r)}}) \right]$  
				is the increasing limit of $u_n(x):= \bE^x \left[ (u\wedge n) (X_{\tau_{B(x_0, r)}}) \right]$ 
				as $n\to \infty$  $\sE$-q.e. on $\sX$, and by  \cite[Theorem 3.4.8]{CF}, 
				$u_n\geq 0$ is a bounded function in $\sF_e$ that is regular harmonic in $B(x_0, r)$; see also the proof of Proposition \ref{P:2.10}(i).

   \item		  Note that if    property  \eqref{e:3.14} 
     holds for a ball $B=B(x_0,r)$ with constants $C_B$ and  $\delta_B$ then it holds
			 for any larger ball $B(x_0,R)$ with constants $C_B$ and $ \delta_B r/R$.
	\end{enumerate} 
	 \end{remark}

\ms	
	
	\begin{thm}\label{T:4}
Assume that 	{\rm (Ha)}  holds  and that 
		$$
		\lambda_\sX:= \inf \left\{ \sE (f, f): f\in \sF \hbox{ with } \| f \|_{L^2(\sX; m)}=1 \right\}>0 .
		$$
		Then {\rm (HC)} holds.   
	 	Moreover,  the $\sE$-nest $\{F_n\}$ in {\rm (HC)} can be taken so that for any compact 
		subset $K\subset \sX$ and any $x_0 \in \sX \setminus K$, there exist some $r_0>0$ with
		$B(x_0, r_0) \subset K^c$ and a constant $C=C(K, x_0, r_0)>0$ such that for any 
		$f\in L^1(K; m)$, $Gf$ is continuous on $K^c \cap F_n$ for any $n \ge 1$ and  
		\begin{equation}\label{e:15} 
		\sup_{x\in B(x_0, r_0) \setminus \sN}  \abs{Gf(x)} \leq 
		C \int_K \abs{f(x)} m(dx),
		\end{equation} 
	where $\sN:= \sX \setminus \cup_{n=1}^\infty F_n$.
	\end{thm}
	
	\pf  Since $\lambda_\sX >0$, by \cite[Lemma 2.1.4(ii) and Theorem 2.1.12(i)]{CF} we have $Gf\in \sF_e= \sF \subset  L^2(\sX; m)$  and $\|Gf\|_{L^2}\leq \lambda_\sX^{-1}\|f\|_{L^2}$ for every $f\in L^2(\sX; m)$. 
	
	Under $\lambda_\sX>0$ and {\rm (Ha)},
	for every $x_0\in \sX$,    $r\in (0, d_{\sX}(x_0))$,  
	and any ball $B(y_0, R)\subset  \sX \setminus B(x_0, r )$, 
	  by Remark \ref{R:3.10}(iii) above and 
	the same argument as that for \cite[(5.10)]{GH}, 
	we have  for any   $f\in L^1(\sX; m)$ with $f=0$ on $B(y_0, \delta_{B(y_0, R)} R)^c$, 
	\begin{equation}\label{e:13}
	{\rm esssup}_{B(x_0, \delta r)} |G f| \leq \frac{C_{B(x_0, r)} C_{B(y_0, R) }}{\lambda_\sX \sqrt{m(B(x_0,  \delta_{B(x_0, r)} r)) 
	 m(B(y_0, \delta_{B(y_0, R)}  R))}} \| f\|_{L^1(\sX; m)}.
	\end{equation}
 
Observe that, in the same way as \eqref{e:n1},
 by the strong Markov property of $X$, for such $f$, $Gf$ is regular harmonic 
in $B(x_0, r)$; note that $G \abs{f}<\infty$  $\sE$-q.e.~by \cite[Proposition 2.1.3(i) and Theorem A.2.13(v)]{CF}. 
 
	Let $\{x_k; k\geq 1\}\subset \sX$ be a  sequence of points in $\sX$, 
   $\{r_k; k \ge 1\}=\bQ \cap (0,\infty)$,  
	and 
\begin{eqnarray*}
 \Lambda  &:=& \big\{\eta =(x_i, x_j, r_k, r_l): i, j, k, l\geq 1, 
  r_k \in (0, d_\sX(x_i) \wedge 1) ,  \,  r_l \in (0, d_\sX (x_j)\wedge 1) ,    \\
&&\hskip 1.3 truein  \hbox{ and } 	 B(x_i, r_k )\cap B(x_j, r_l)=\emptyset \big\}.
\end{eqnarray*}
Let $\delta_{i,k}\in (0, 1)$ be one half of the largest positive constant $\delta_B$ in {\rm (Ha)} for the validity of the Harnack inequality in the ball $B(x_i, r_k)$ and $C_{i,k}>0$ be the corresponding comparison constant. We select $\delta_{i,k}$ in this way to ensure the uniformity of the constants
$\delta_B$ when the centers and radii of the balls are close to each other. This uniformity is needed
when we do the finite covering of such balls over each compact set $K$ in the last step of this proof.
 	Note that $\Lambda$ is a countable set. For each $\eta = (x_i, x_j, r_k, r_l)\in \Lambda$, 
  	 let $\{f_p, p\geq 1\}$ be a dense sequence in $C_c(B(x_j, \delta_{j, l} r_l)) $ with respect to the 
	supremum norm chosen so that $\{f_p\}_{p \ge 1} \cap C_c(U_n)$ is uniformly dense in $C_c(U_n)$ for any $n \ge 1$, where $\{U_n\}_{n \ge 1}$ is an increasing sequence of relatively compact open subsets of $B(x_j, \delta_{j,l}r_l)$ with $\bigcup_{n=1}^\infty U_n= B(x_j, \delta_{j,l}r_l)$.
   Since $f_k\in L^2(\sX; m)$, $Gf_k\in \sF$ and  it is quasi-continuous by   \cite[Theorem 4.2.6]{FOT}.
	Thus  there is an ${\sE}$-nest $\{F_n^{(\eta)}, n \geq 1\}$ consisting of an
	increasing sequence of compact sets such that $F_n^{(\eta)}=\on{supp}\left[\one_{F_n^{(\eta)}}m \right]$ for any $n \ge 1$, $\sN_0\subset \sX \setminus \cup_{n=1}^\infty
	F^{(\eta)}_n$ and  	$Gf_p$ is finite and continuous
	on each $F_n^{(\eta)}$ for every integer $p\geq 1$; see \cite[Lemma 1.3.1]{CF}.
	Let ${\cal N}_\eta:=\sX \setminus \cup_{n=1}^\infty
	F^{(\eta)}_n$, which is ${\cal E}$-polar and in particular has zero
	$m$-measure, and  
	$C_\eta:=\frac{C_{i, k} C_{j, l} }{\lambda_\sX \sqrt{m(B(x_i,   \delta_{i,k} r_k)) m(B(x_j,  \delta_{j, l} r_l ))}} $.
	Inequality \eqref{e:13} yields that, with $f_0:=0$, for any $k_1,k_2 \ge 0$, 
 	\begin{equation}\label{e:3.17} 
	\sup_{x\in B(x_i, \delta_{i,k} r_k)\setminus {\cal N}_\eta} |G f_{k_1} (x) -G f_{k_2} (x)| \leq  C_\eta
	\| f_{k_1} -f_{k_2}\|_{L^1(B(x_j,  \delta_{j, l}  r_l); m)} .
	\end{equation}
	Since $\{f_p,  p\geq 1\}$ is uniformly dense in $C_c(B(x_j, \delta_{j,l} r_l) $ and $m$ is a Radon measure
	on $\sX$, 
	it follows that for every $ f\in   C_c(B(x_j, \delta_j r_j)$,
	 $G  f$ is continuous on   $B(x_i, \delta_{i,k} r_k)\cap F_n^{(\eta)}$ for any $n \ge 1$ and
	\begin{equation}\label{e:3.18} 
	\sup_{x\in B(x_i, \delta_{i,k} r_k) \setminus {\cal N}_\eta} |G f(x)| \leq C_\eta   \| f \|_{L^1(B(x_j,  \delta_{j, l}  r_l); m)}.
	\end{equation} 
	Let $D$ be an open subset of $B(x_j, \delta_{j,l} r_j)$.
	Since $1_D$ can be approximated pointwise by an increasing sequence of non-negative continuous functions
	with compact support in $B(x_j, \delta_{j,l} r_j)$, we have from \eqref{e:3.18}
	and the monotone convergence theorem that 
	\begin{equation}\label{e:3.19} 
	\sup_{x\in B(x_i, \delta_{i,k} r_k) \setminus {\cal N}_\eta}  G (x, D)   \leq C_\eta   m(D). 
	\end{equation} 
	For any Borel subset $A\subset B(x_j, \delta_{j,l} r_j)$, let $\{D_{k_1}; k_1\geq 1\}$ be a decreasing
	sequence of open subsets of $B(x_j, \delta_{j,l} r_j)$ so that $A \subset \cap_{k_1 \ge 1} D_{k_1}$ and
	$\lim_{k_1 \to \infty} m( D_{k_1}\setminus A)=0$. By \eqref{e:3.19} and the dominated convergence theorem,
	$$
	\sup_{x\in B(x_i, \delta_{i,k} r_k) \setminus {\cal N}_\eta}  G (x, A)   \leq 
	\lim_{k_1\to \infty} \sup_{x\in B(x_i, \delta_{i,k} r_k) \setminus {\cal N}_\eta}  G (x, D_{k_1}) 
	\leq \lim_{k_1 \to \infty} C_\eta   m(D_{k_1})= C_\eta   m(A).
	$$
	The above in particular implies that for each $x \in B(x_i, \delta_{i,k} r_k) \setminus {\cal N}_\eta$, 
	\begin{equation}\label{e:3.20}  
	\begin{aligned}
	&G (x, dy)   \hbox{ is absolutely continuous with $[0,C_\eta]$-valued density} \\
	&\hskip 0.8truein  \hbox{with respect to $m(dy)$ on } B(x_j, \delta_{j,l} r_l).
	\end{aligned}
	\end{equation} 
	 	Since $\{f_p,  p\geq 1\}\subset C_c(B(x_j, \delta_{j,l} r_j))$ is dense in $L^1(B(x_j,  \delta_{j, l}  r_l); m)$, 
	it follows from \eqref{e:3.17} and \eqref{e:3.20}
	that $G  f$ is continuous on $B(x_i, \delta_{i,k} r_k)\cap F_n^{(\eta)}$ for any $n \ge 1$ and \eqref{e:3.18} holds
	for every $ f\in   L^1(\sX; m) $  with $f=0$ on $B(x_j, \delta_{j, l} r_l )^c $.

	By \cite[Lemma 1.3.1]{CF} and its proof, by taking suitable intersections of $F^{(\eta)}_{n_k}$'s,  
	there is an ${\cal E}$-nest $\{F_n , n \geq 1\}$ consisting of an
	increasing sequence of compact subsets of $\sX$ such that 
	$\sN_0\subset \sX \setminus \cup_{n=1}^\infty F_n$ and that 
	for every $\eta=(x_i,x_j,r_k,r_l) \in \Lambda$, 
	$G  f$ is continuous on   $B(x_i, \delta_{i,k} r_k)\cap F_n $ for any $n \ge 1$ and
	\begin{equation}\label{e:14} 
	\sup_{x\in B(x_i, \delta_{i,k} r_k) \setminus {\cal N} } |G f(x)| \leq
	C_\eta   \| f \|_{L^1(B(x_j,  \delta_{j, l}  r_l); m)} 
	\end{equation}
	for every $ f\in   L^1(\sX; m) $  with $f=0$ on $B(x_j, \delta_{j,l} r_l)^c $, where $\sN:=\sX \setminus \cup_n F_n$ which is $\sE$-polar.

	Let $K$ be a compact subset of $\sX$.
		As $\{x_i\}$ is dense in $\sX$ one can deduce from \eqref{e:14} by finite covering  
	that for any $f\in L^1(\sX; m)$ that vanishes outside $K$, 
	 	$Gf$ is continuous on   $K^c\cap F_n $ for any $n \ge 1$, and for any $x_0\in K^c$, there is some $r_0>0$ 
		with $B(x_0, r_0)  \subset K^c$ so that  
	$$
	\sup_{x\in B(x_0, r_0)  \setminus {\cal N} } |G f(x)| \le	C(K, B(x_0,r_0))    \| f \|_{L^1(K; m)}.
 $$
   This in particular proves that (HC) holds with the above $\sE$-nest $\{F_n\}$ of compact sets. 
	\qed
 	 	
	Under the assumption of Theorem \ref{T:4}, we have by \eqref{e:15} that 
	for every compact subset $K\subset \sX$, $G(x, K)<\infty$ for $x\in (\sX \setminus \sN)\setminus   K$. In other words, $\restr{G(x,\cdot)}{\sX\setminus \{x\}}$ is a Radon measure on $\sX \setminus \{x\}$ for any $x \in \sX\setminus \sN$.

	\medskip

	Using a time change argument, we can remove  the assumption of $\lambda_\sX>0$ in
 	Theorem \ref{T:4}.
	
	\begin{thm}\label{T:5}
		 		Assume that  {\rm (Ha)}  holds.
	 	  Then   the conclusions of Theorem \ref{T:4}  hold.   
		Consequently,  the conclusions of Theorems \ref{T:1} and  \ref{T:2}  hold. 
	\end{thm}

	\pf 	Recall that our running assumption is that  the Dirichlet form $(\sE, \sF)$ on $L^2(\sX; m)$
	(or equivalently, its associated Hunt process $X$)  is transient.
	Let $g_0$ be as in Lemma \ref{L:g0}, and $\mu (dx)=g_0(x) m(dx)$.
	We now make a time change of $X$ via the inverse of the positive continuous additive functional
	$A_t:=\int_0^t g_0(X_s)ds$. That is, let $Y_t=X_{\tau_t}$, where $\tau_t:=\inf\{s>0: A_s>t\}$.
	Then $Y$ is $\mu$-symmetric and transient, and its extended Dirichet form is the same as that of $X$; see \cite[Theorem 6.2.1 and Corollary 6.2.1]{FOT} and \cite[Corollary 5.2.12]{CF} (since $\mu$ and $m$ are mutually absolutely continuous). 
	So the Dirichlet form of $Y$ is $(\sE, \sF_e \cap L^2(\sX; \mu))$ on $L^2(\sX; \mu)$. 
	 	Since $Y$ and $X$ share the same family of regular harmonic functions,  {\rm (Ha)} 
	holds for $Y$.  We claim that 
	\begin{equation}\label{e:16}
	\lambda_\sX^Y:=\inf\{\sE (f, f):  f\in \sF_e \cap L^2(\sX; \mu) \hbox{ with } \| f \|_{L^2(\sX;\mu)}=1\} \geq   1.
	\end{equation}
	Denote by $\wt G$ the Green potential of $Y$, that is, for $f\geq 0$ on $\sX$,
	$$ 
	\wt Gf(x):= \bE^x \int_0^\infty f(Y_t) dt =\bE^x \int_0^\infty f(X_{\tau_t}) dt.
	$$
	Using the time change, we see that $\wt Gf(x)= \bE^x \int_0^\infty (fg_0) (X_t) dt = G(fg_0)(x)$.
	In particular, we have $\wt G1 =Gg_0\leq 1$. Thus  for $u\in L^2(\sX; \mu)$,  by Cauchy-Schwarz 
	and the symmetry of $\wt G$ with respect to $\mu$, 
	\begin{eqnarray}
	\int_\sX (\wt Gu )^2  (x)\mu (dx) & \leq &\int_\sX \wt  G (u^2) (x)  \wt G1 (x) \mu (dx) \leq  \int_\sX \wt  G (u^2) (x)   \mu (dx)  \nonumber\\
	&\leq & \int_\sX  u(x)^2 \wt G1(x) \mu (dx) \leq \int_\sX  u(x)^2   \mu (dx).  \label{e:17}
	\end{eqnarray} 
 	Therefore, for any $u\in L^2(\sX; \mu)$,  $\int_\sX u (x)\wt Gu (x) \mu (dx)\leq \int_\sX u(x)^2 \mu (dx) <\infty$ by \eqref{e:17}.
	It follows  (cf. \cite[Theorem 2.1.12]{CF} or  \cite[Theorem 1.5.4]{FOT}) that $\wt Gu \in \sF_e \cap L^2(\sX; \mu)$  
	with $\sE (\wt Gu, \wt Gu)=\int_\sX  u  (x)\wt Gu  (x)\mu (dx)$.
	Hence for   $u\in \sF_e\cap L^2(\sX; \mu)$,  we have by \eqref{e:17} and the Cauchy-Schwarz,
	$$
	\int_\sX u^2(x) \mu (dx)  =  \sE (\wt Gu, u) \leq \sE (u, u)^{1/2} \sE(\wt Gu, \wt Gu )^{1/2}
	\leq   \sE (u, u)^{1/2} \, \|u\|_{L^2(\sX; \mu)}.
	$$ 
	Consequently,  
\[  \|u\|_{L^2(\sX; \mu)} \leq \sE (u, u)^{1/2} \quad \hbox{for every } u\in \sF_e\cap L^2(\sX; \mu).
\]
	This proves the claim \eqref{e:16}. 
	
	For the process $Y$, we can take $g^Y_0=1$   in  the role of $g_0$ for $X$ in \eqref{e:4a} as $\wt G1= Gg_0\leq 1$.
	By Theorem \ref{T:4}, {\rm (HC)} holds for the process $Y$.  Since $\wt G f= G(fg_0$) and   
    an increasing sequence of  compact sets $\{F_n\}$ of $\sX$ is an $\sE$-nest
	for the process $X$ if and only if it is an $\sE$-nest for the process $Y$ in view of  \cite[Theorem 5.2.11]{CF}, 
	we conclude that {\rm (HC)} as well as the rest of the conclusions of Theorem \ref{T:4} hold for the process $X$. 
	\qed

	\medskip

 	  \begin{remark}\label{R:3.13} \rm
	One can see from its proof  that 
 the properly exceptional set $\sN$ in the conclusion of Theorem \ref{T:2}  in fact depends only 
 on the properly exceptional set  $\sN_1\supset \sN_0$ of $X$ for which  
 \begin{equation} \label{e:3.24}
    \bP^x(X_t\in dy) \hbox{ is absolutely continuous with respect to } m  \hbox{ for every } x\in \sX \setminus \sN_1
    \hbox{ and } t>0.
    \end{equation}
 By Theorem \ref{T:1}, under condition {\rm (HC)}, $\sN_1$ depends only on the $\sE$-nest $\{F_n\}$ in the definition of {\rm (HC)}. 
 In view of the proof of Theorem \ref{T:5}, under condition  {\rm (Ha)},   the conclusion of Theorem \ref{T:2} holds 
for a properly exceptional set $\sN$ that only depends on $\sN_1$ in \eqref{e:3.24} and the properly exceptional set   in the property 
that $Gg_0\leq 1$ $\sE$-q.e. on $\sX$. 
\end{remark}

\medskip
	 
  \begin{remark}\label{R:3.14} \rm
 In Example \ref{E:8.3}, we will give an example of an irreducible  strongly local regular Dirichlet form for which 
(Ha) fails but (HC) holds. In fact, for this example,
there exists a discontinuous positive harmonic function.  
	\end{remark}

	\section{ Green functions } \label{S:4}

	We now drop the hypothesis that $(\sE,\sF)$ is transient.   Recall that $d_{\sX}(x)$ is the distance to the boundary function on $\sX$ defined in \eqref{e:3.5}.
	  	
	\medskip
	
	\begin{definition}\label{D:ehi}  \rm   For a MMD space $(\sX,d,m,\sE, \sF )$, we say
			 
			\begin{enumerate} [\rm (i)]
			
		\item	the {\em (scale invariant) elliptic Harnack inequality} (EHI)  holds if there exist
			constants  $\delta_H \in (0, 1)$ and $C_H \in (1,\infty)$ so that 
			 for any $x \in \sX$,  	
 			 $R\in (0, d_{\sX}(x))$,   
			and for any nonnegative bounded harmonic function $h$  on a ball $B(x,R)$,  one has
			\be \label{e:ehi2} 
			\esssup_{B(x,\delta_H R)} h \le C_H  \essinf_{B(x, \delta_H R)} h;
			\ee

\item  the {\em (scale invariant) local elliptic Harnack inequality} ${\rm EHI}_{\rm loc}$  if 
   	  \eqref{e:ehi2}  holds for  any nonnegative bounded harmonic function on
 			 balls $B(x,  R)$ with  $0<R< d_\sX(x)  \wedge 1$. 
			   \end{enumerate} 
			  \end{definition}
	 	
	\begin{remark} \label{R:4.2} \rm 
	 	\begin{enumerate} [(i)]
	 	\item  Clearly,   the EHI implies 
 	 the ${\rm EHI}_{\rm loc}$, 
	and the ${\rm EHI}_{\rm loc}$ implies  implies (Ha).
	 	
	\item   If $(\sX,d)$ is a geodesic metric space and    
 	 ${\rm EHI} $  	  or  ${\rm EHI}_{\rm loc}$ 
	 holds for 
	some value of $\delta_H$, then the same holds for any other $\delta' \in (0,1)$ with a constant $C_H(\delta')$.
	 
	\item  Suppose ${\rm EHI}_{\rm loc}$ holds  and 
		 $u$ is a nonnegative  harmonic function on a 
			 ball $B(x,  R)$ with $0<R< d_\sX(x) \wedge 1$.  
			 By Proposition \ref{P:2.10}(ii), 
			$u$ is the increasing limit in $B(x, R/2)$ of a sequence of   functions
			$\{u_n; n\geq 1\}$ that is non-negative, bounded and harmonic in $B(x, R/2)$. 
	    It follows that \eqref{e:ehi2} holds for $u$ on the ball $B(x, \delta_H'R)$, 
			 where $\delta_H'=\delta_H/2$. 
		In other words, if ${\rm EHI}_{\rm loc}$ holds, it holds for any non-negative (possibly unbounded) harmonic function.

	\item  If ${\rm EHI}_{\rm loc}$ holds, in view of (iii) above,
	iterating the condition \eqref{e:ehi2} gives a.e.~H\"older continuity of harmonic functions, and it follows that any  locally bounded harmonic function has a continuous modification \cite[Lemma 5.2]{GT12}.
   	\end{enumerate} 
	\end{remark}

	 	  Let $D$ be an open set of $\sX$.
	Note that if $(\sX,d,m,\sE,  \sF )$ satisfies  the ${\rm EHI}_{\rm loc}$, then so does $(D,d,m|_D, \sE,  \sF^D )$,
	where $(\sE, \sF^D)$ is the Dirichlet form for the part process $X^D$ of $X$ killed upon leaving $D$
	(see \eqref{e:FD}). 
	Let $D_{\diag}$ denote the diagonal in $D \times D$. 
	 	For a subset $A\subset \sX$, we use $\ol A$ to denote its closure and $\partial A$ its boundary.

	\medskip

	 	\begin{definition}  \label{D:goodgreen} \rm 
	 		\begin{enumerate}
	 			\item[(a)]
			Let $D$ be a  non-empty open subset of $\sX$ such that $D^c$ is not $\sE$-polar. 
			We say that $(\sE,\sF)$ has a {\em regular Green function}     
			on $D$ if  there exists a 
		 non-negative $\sB (D\times D)$-measurable function 
			$g_D(x,y)$   on $D \times D$
			with the following properties:
			\begin{enumerate}[(i)]
				\item (Symmetry) $g_D(x,y)= g_D(y,x)  $ for all $(x,y) \in D \times D$;
				
				\item (Continuity) $g_D(x,y)$ is $[0,\infty)$-valued and jointly continuous in $(x,y) \in D \times D \setminus D_{\diag}$;
				
				\item (Occupation density) 
       There is a Borel properly exceptional set $\sN_D$ of $X^D$ such that
				\begin{equation}\label{e:4.2}
				\bE^x \int_0^{\tau_D} f(X_s) ds = \int_D g_D(x, y) f(y)\, m(dy)
				\quad \hbox{for every }  x\in D\setminus \sN_D, 
				\end{equation}  
				 for any   $f\in \sB_+ (D)$;
				
				\item (Excessiveness) For each 
				$y\in D$, $x\mapsto g_D(x, y)$
				is $X^D|_{D\setminus \sN_D}$-excessive, where $\sN_D$ is the set given in (iii);

				\item (Harmonicity) For any fixed $y \in D$, the function $x \mapsto g_D(x,y)$ is 
				in $\sF_{\rm loc}^{D \setminus \{y\}}$ and  
			 for any open subset $U$ of $D$ with $y \notin \overline{U}$, 
					$x\mapsto g_D(x, y)$ is regular harmonic in $U$ with respect to $X^D$.

				\item (Maximum principles) If $x_0 \in U \Subset D$, then
				\begin{equation}\label{e:max}
				\inf_{\overline{U} \setminus \set{x_0}} g_D(x_0,\cdot) = \inf_{\partial U} g_D(x_0,\cdot),
				\qquad \sup_{D \setminus U} g_D(x_0, \cdot) = \sup_{\partial U}g_D(x_0,\cdot).
				\end{equation}
			\end{enumerate}
 We call   $g_D(x, y)$ the regular Green function of $(\sE, \sF)$ in $D$.
			
			\medskip
			
		\item[(b)]	We say that $(\sE,\sF)$ has {\em regular Green functions} if 
			for any non-empty open set $D \subset \sX$ whose complement $D^c$ is not $\sE$-polar, 
			$(\sE, \sF)$ has a regular
			Green function $g_D(x, y)$ on $D$, 
 where the properly exceptional set $\sN_D$ in (iii),(iv) can be taken to be $\sN_D=D \cap \sN$ for a Borel properly exceptional set $\sN$ of $X$ that is  independent of $D$ as long as $D \subset B(x,R)$ for some $x \in \sX$ and $R \in (0,\diam(\sX,d)/2)$.  
\end{enumerate}
	\end{definition}

	\begin{theorem} \label{T:reggreen}
 	Suppose that the MMD space  $(\sX, d, m , \sE, \sF)$  and let  $D$ be a non-empty
 	open subset of $\sX$  such that  the part process $X^D$ is transient.
		\begin{description}
			\item{\rm (i)} 
			   Assume that  
			$(D, d, m|_D, \sE, \sF^D)$ satisfies {\rm (Ha)}.
			 			Then $(\sE,\sF^D)$ has a Green function $g_D(x, y)$ in the sense that 
			\begin{enumerate} 			
			\item[{\rm (i.a)}] $g_D(x, y) $  is  a non-negative  jointly   ${\cal B}(D \times D)$-measurable
			function  and there is   a Borel properly exceptional
			set  $\sN_D$ of $X^D$ such that 
			$$
			\bE^x \int_0^{\tau_D} f(X_s) ds = \int_D g_D(x, y) f(y) \, m(dy) , \quad x\in D\setminus \sN_D, 
			$$
			for any  $f\in \sB_+ (D)$;
			
			\item[{\rm (i.b)}] $g_D(x, y)=g_D (y, x)$ for every $x, y\in D \setminus \sN_D$, 
			and 
			$x\mapsto g_D(x, y)$ is $X^D|_{D\setminus \sN_D}$-excessive and bounded  $\sE$-q.e.~on any $U \Subset D \setminus \{y\}$ for every $y\in D\setminus \sN_D$;
			
			\item[{\rm (i.c)}] For every $y_0\in  D \setminus \sN_D$, $x\mapsto g(x, y_0)$ is harmonic in $ D  \setminus \{y_0\}$.  Moreover, for any open subset $U$ of $D$ with $y_0 \notin \overline{U}$ is regular harmonic
			in $U$ with respect to $X^D$. 
			\end{enumerate} 
			
			\item{\rm (ii)} If $(D, d, m|_D, \sE, \sF^D)$ satisfies  the ${\rm EHI}_{\rm loc}$,   
			then  $(\sE, \sF)$ has a regular Green function on $D$. 
		\end{description}
		In particular, the above properties hold if   $(\sX, d, m , \sE, \sF)$ is irreducible and $D$ is a non-empty
		open subset of $\sX$  such that $D^c$ is not $\sE$-polar
	\end{theorem}
	
	\pf   	(i) 
	The conclusion of this part   follows directly from Theorem \ref{T:5} by replacing $\sX$ and $X$ by $D$ and $X^D$, 
	 and from Proposition \ref{P:2.9}.   Denote the corresponding Borel properly exceptional set in Theorem 
	\ref{T:5} for $X^D$ by $\sN_D$. The local   $\sE$-q.e.  boundedness  of  $x \mapsto g_D(x,y)$ on $D \setminus \{y\}$ for every $y \in D \setminus \sN_D$ follows from      (Ha),  \eqref{e:15} of Theorem \ref{T:5}, \cite[Theorem A.2.7]{FOT} and \cite[Theorem A.2.13(iv)]{CF}.
 Note that this $\sN_D$ has the property that for every $x\in D\setminus \sN_D$,
	the law of $X^D_t$ under $\mathbb P^x$ is absolutely continuous with respect to $m|_D$
	for every $t>0$. This property will be used in \eqref{e:4.7} 
	as well as at the end of this proof when establishing
	the excessiveness of $x\mapsto g_D (x, y)$ for every $y\in D$.

	\smallskip
	
	(ii)  Suppose now that $(D, d, m|_D, \sE, \sF^D)$ satisfies the 
 	 ${\rm EHI}_{\rm loc}$ in $(D, d)$. Let $x_0 \in D$ and $r \in (0, (d_D(x_0) \wedge 1)/2)$.  	
	 Let $u$ be a bounded harmonic function in $B(x_0, 2r) \subset D$. Iterating the condition \eqref{e:ehi2} yields that there are constants $c_0 >0$ and $\beta \in (0, 1)$ that depend only  $\delta_H (D)$ and $C_H(D)$ in \eqref{e:ehi2} (with $D$ in place of $\sX$) such that 
	\begin{equation}\label{e:3.6}
	|u(x)-u(y)| \leq c_0 \| u\|_{L^\infty (B(x_0, 2r)} \left( \frac{d(x,y)}{r}\right)^\beta \quad \hbox{for a.e. } x, y \in B(x_0, r).
	\end{equation}
	Since {\rm (Ha)} holds on $D$, by (i) there is a Green function  $g_D(x, y)$ in $D$.
	For each fixed $y_0\in D\setminus \sN_D$,  $x\mapsto g_D(x, y_0)<\infty$ $m$-a.e. and is harmonic in $D \setminus \{y_0\}$. 
	It follows from   \eqref{e:15} of Theorem \ref{T:5} (of ${\rm EHI}_{\rm loc}$ 
   and Remark \ref{R:4.2}(iii))
 	that $x\mapsto g_D(x, y_0) $   is  (essentially) locally bounded in $D \setminus \{y_0\}$. 
The H\"older estimate \eqref{e:3.6} implies that there is a locally H\"older continuous function $\wt g_D(\cdot, y_0)$ on $D\setminus \{y_0\}$
such that $\wt g_D(x, y_0)= g_D(x, y_0)$ for $m$-a.e.  $x\in D$. 
Since $g_D(\cdot, y_0) $ is $X^D|_{D\setminus \sN_D}$-excessive,
$t\mapsto g_D(X^D_t, y_0) \in [0, \infty ]$ is right continuous on $[0, \infty)$ 
$\bP^x$-a.s. for every $x\in D \setminus \sN_D$; see, e.g.,  \cite[Theorem A.2.2]{CF}.
Since the law of $X^D_t$ under $\bP^x$ is absolutely continuous with respect to $\restr{m}{D}$ 
for every $x\in D\setminus \sN_D$, we have for every 
$x\in D\setminus (\sN_D\cup \{y_0\})$, $\bP^x$-a.s., 
\begin{equation}\label{e:4.7}
g_D(x, y_0) = \lim_{\bQ\ni t\to 0} g_D (X^D_t, y_0) = \lim_{\bQ  \ni   t\to 0} \wt g_D (X^D_t, y_0) = \wt g_D(x, y_0).
\end{equation}
 	Thus  $g_D(x, y_0) = \wt g_D(x, y_0)$ for every $x\in D\setminus (\sN_D\cup \{y_0\})$.
 For each $y_0\in D\setminus \sN_D$, we define $\wt g_D(y_0, y_0) := g_D(y_0, y_0)$.  The above shows that there is a function $\wt g_D(x, y)$ defined on $D\times (D\setminus \sN_D)$ so that for each $y_0\in D\setminus \sN_D$, $x\mapsto \wt g_D(x, y)$ is locally H\"older continuous on
		$D\setminus \{y_0\}$ and $\wt g_D(x, y_0) =g_D(x, y_0)$ for every $x\in D \setminus \sN_D$. 

 Since $g_D(x, y)$ is symmetric on $(D \setminus \sN_D) \times (D\setminus \sN_D)$ by (i.b), we have $\wt g_D(x, y)=\wt g_D(y, x)$ for every $x,  y\in D \setminus \sN_D$. 
	We extend the definition of $\wt g_D(x, y)$ on $ D\times (D  \setminus \sN_D )$
	to $(D\times D ) \setminus  (\sN_D \times \sN_D)$ by setting 
	 $\wt g_D(x, y) = \wt g(y, x)$ for $x\in D\setminus \sN_D$ and $y\in \sN_D$.
 	Note that $y\mapsto \wt g_D(x_0, y)$ is continuous in $D\setminus \{x_0\}$.
 	We next show that such defined $\wt g_D(x, y)$ 
 on $(D\times D ) \setminus  (\sN_D \times \sN_D)$
	is locally jointly H\"older continuous off the diagonal and thus  its definition  can be continuously extended to $(\sN_D\times \sN_D)\setminus  \{ (x, x): x\in \sN_D\} $. 
	
	Let $x_0, y_0\in D\setminus \sN_D$ with $x_0\not= y_0$. There is $r \in (0, (\rho_D(x_0)\wedge \rho_D(y_0))/2)$ such that $B(x_0, 2r)\cap B(y_0, 2r) = \emptyset$. 
	By the ${\rm EHI}_{\rm loc}$ in $D$,  for every $x\in B(x_0,  2 \delta_H(D)r)\setminus \sN_D$ and $y\in B(y_0,  2 \delta_H(D) r)\setminus \sN_D$, 
	$$
	\wt g_D(x, y ) \leq C_H(D) \wt g_D(x , y_0) \leq  C_H(D)^2 \wt g_D(x_0 , y_0).
	$$
	It follows from \eqref{e:3.6} that for $x_1, x_2 \in B(x_0,  \delta_H(D)r)\setminus \sN_D$ and $y_1, y_2\in B(y_0,  \delta_H(D)r)\setminus \sN_D$,
	\begin{eqnarray*}
		| \wt g_D(x_1, y_1)-\wt g_D(x_2, y_2)|
		&\leq & | \wt g_D(x_1, y_1)-\wt g_D(x_2, y_1)|+ | \wt g_D(x_2, y_1)-\wt g_D(x_2, y_2)| \\
		&\leq & C_H(D)^2 \wt   g_D(x_0 , y_0) c_0 \left(\delta_H(D) r\right)^{-\beta}\left( d(x_1,x_2)^\beta + d(y_1,y_2)^\beta\right).
	\end{eqnarray*}
	Consequently $\wt g_D(x, y)$ can be extended continuously to  
	 $B(x_0, \delta  r)\times B(y_0, \delta r)$  and hence to $D\times D\setminus D_{\diag}$
	as a locally H\"older continuous function. Clearly, $\wt g_D (x, y)=\wt g_D (y, x)$ for $x, y\in D$ with $x\not=y$.
 For $w\in \sN_D$, we define $\wt g_D(w, w) :=\limsup_{x\not=y \in D, x\to w, y\to w}
		\wt g_D(x, y)$.
	In summary, we now have a symmetric Borel measurable function $\wt g_D(x, y):
	D\times D\to [0, \infty]$ defined on $D\times D $ so that $\wt g_D(x, y)$ is locally
	H\"older continuous on $ (D\times D)\setminus D_{\diag}$, and $\wt g_D(x, y)= g_D(x, y)$ for $x, y\in D\setminus \sN_D$.
	From now on, we   take this refined  version $\wt g_D(x, y)$ for the Green function $g_D(x, y)$ 	in $D$ and drop the tilde from $\wt g_D(x, y)$.

 	 We next show that $ g_D(x, y)$ is a regular Green function on $D$.
	We already have the symmetry and continuity properties
	of $g_D(x, y)$. The occupation density holds
	for this refined $g_D(x, y)$ as well since $\sN_D$ is $\sE$-polar so $m(\sN_D)=0$.
	Thus  it remains to show the excessiveness, regular harmonicity  and maximum principle for $g_D(x, y)$  (see Definition \ref{D:goodgreen}).
	Suppose $U$ is a relatively compact open subset  of $D$ and $x_0\in U $. 
	Let $r_0>0$ be such  that $B(x_0, r_0) \subset U$. For  every $y\in D\setminus ( U\cup \sN_D)$ and for every $r\in (0, r_0)$, 
	we have by the symmetry of $g_D(x, y)$ and the strong Markov property of $X$ that  
	\begin{eqnarray*}
		\int_{B(x_0, r)}  g_D(x, y)  m(dx) 
		&=& \bE^{y} \int_0^{\tau_D} 1_{B(x_0, r)} (X_s) ds 
		= \bE^{y} \bE^{X^D_{\sigma_U}} \int_0^{\tau_D} 1_{B(x_0, r)} (X_s) ds \\
		&=&\int_{B(x_0, r)}  \bE^y    g_D(x, X^D_{\sigma_U}) m(dx).
	\end{eqnarray*}
	Dividing both sides by $m(B(x_0, r))$ and then taking $r\to 0$ yields that 
	\begin{equation}\label{e:4.5} 
		g_D(x_0, y)=  \bE^y \left[  g_D(x_0, X^D_{\sigma_U}) \right].
		\end{equation} 
 This together with the continuity of $g_D$ on $D\times   D \setminus D_{\diag}$ shows that 
	\begin{equation}\label{e:max1}
	\sup_{y\in D \setminus U}   g_D(x_0,  y) = \sup_{y\in \partial U}   g_D(x_0, y).
	\end{equation} 
	Identity \eqref{e:4.5} holds with $\tau_V$ in place of $\sigma_U$ 
	for any open subset $V$ of $D$ with $x_0 \notin \overline{V}$ and any $y \in V \setminus \sN_D$ by exactly the same argument. For this, we note that \eqref{e:max1} along with $B(x_0,r_0) \Subset D \setminus \overline{V}$ is used to apply dominated convergence theorem.
	Thus the regular harmonicity in Definition \ref{D:goodgreen}(v) holds for $g_D(x, y)$. In particular, $g_D(\cdot,y) \in \sF_{\on{loc}}^{D \setminus \{y\}}$ for any $y \in D$ by Proposition \ref{P:2.9}(iii).
	
	For each fixed $x\in U\setminus \sN_D $ and  any Borel measurable function $f \geq 0$ on $D$,  by
	the strong Markov property of $X$,
	$$
	\int_D g_D(x, y) f(y) m(dy) = \bE^x \int_0^{\tau_D} f(X_s) ds \geq  \int_D \bE^x [ g_D (X^D_{\tau_U}, y) ] f(y) m (dy).
	$$
	Hence $g_D(x, y) \geq \bE^x [ g_D (X^D_{\tau_U}, y) ] $ for $m$-a.e. $y\in D$.
	Since for every $z\in D$, 
	$y\mapsto g_D(z, y)$ is continuous on $D\setminus \{z\}$,  we have by Fatou's lemma that 
	$$ 
	g_D(x, z) \geq \bE^x [ g_D (X^D_{\tau_U}, z) ] \quad \hbox{for every } z\in U\setminus \{x\}.
	$$
	Taking $z=x_0$ and by the symmetry of $g_D(x, z)$,  we get  for every $ x \in U \setminus (\sN_D \cup\{x_0\}) $, 
	$$
	g_D(x_0, x) \geq \bE^x [ g_D (x_0, X^D_{\tau_U}) ] \geq \inf_{y\in \partial U} g_D (x_0, y).  
	$$
	In the last inequality, we used the fact that $\bP^x(\tau_U<\tau_D)=1$
 due to Proposition \ref{P:3.1}. 
	By the continuity of  $x\mapsto g_D (x_0, x)$ on $D\setminus \{x_0\}$,  
	we get  
	$$
	\inf_{y\in \ol{U}\setminus \{x_0\}} g_D(x_0, y) = \inf_{y\in \partial U} g_D (x_0, y).
	$$ 
	This together with \eqref{e:max1} establishes the maximum principle for $g_D(x, y)$. 
	
 We now show that  for every $y\in D$, $x\mapsto g_D(x, y)$ is $X^D|_{D\setminus \sN_D}$-excessive. Note that $\sN_D$ is a properly exceptional set 
   for $X^D$.
	By (i.b), the above property holds for every $y\in D\setminus \sN_D$. 
	For $y\in \sN_D$, let $y_n\in D\setminus \sN_D$ so that $y_n\to y$ as $n\to \infty$. 
	Let $x\in D\setminus \sN_D$.
	Observe that since  $y$ is $\sE$-polar, $m(\{y\})=0$. 
	As mentioned earlier, the law of $X^D_t$ under $\bP^x$ is absolutely continuous with respect to $m$ for every $t>0$.  Thus by the local H\"older continuity of 
	$g_D$ on $D\times D \setminus D_{\diag}$ and Fatou's lemma, 
	for every $t>0$, 
	\begin{equation}\label{e:exces1}
	\bE^x g(X^D_t, y) \leq \liminf_{n\to \infty} \bE^x g_D (X^D_t, y_n)
	\leq \liminf_{n\to \infty} g_D  (x, y_n) = g_D (x, y). 
	\end{equation}
	 On the other hand, since $\bP^x ( \lim_{t\to 0} X^D_t=x)=1$ and $x\neq y$,
	 we have by Fatou's lemma again,
	 $$
	 g_D(x, y) =\bE^x \left[ \lim_{t\to 0} g_D(X^D_t, y) \right]
	 \leq \liminf_{t\to 0} \bE^x g_D(X^D_t, y) .
	 	$$ 
	 This together with \eqref{e:exces1} proves that $x\mapsto g_D(x, y)$ is $X^D|_{D\setminus \sN_D}$-excessive for every $y\in \sN_D$ and hence for every $y\in D$.   
	 This completes the proof that $g_D(x, y)$  
	is a regular Green function on $D$. This concludes the proof of (i) and (ii).
	
	 Finally, if $(\sE, \sF)$ is irreducible and $D^c$ is not $\sE$-polar,  the regular 
	Dirichlet form $(\sE, \sF^D)$ on $L^2(D; m|_D)$ 
	is transient by  Proposition \ref{P:2.1}.  	 	\qed

   The following  property on the consistency   of exceptional sets follows readily from Remark \ref{R:3.13}. 
This observation will be used to choose a common properly exception set as required   in Definition \ref{D:goodgreen}.
 
	\begin{lem} \label{l:monoton}
	Suppose that  $(\sX,d,m,\sE,\sF )$ is irreducible and  satisfies the ${\rm EHI}_{\rm loc}$. Let $D_1,D_2$ be open subsets of $\sX$ such that $\sX \setminus D_2$ is non-$\sE$-polar and $D_1 \subset D_2$, and let $\sN_{D_2}$ be a Borel properly exceptional set of $X^{D_2}$  in Theorem \ref{T:reggreen}(i). Then we can choose $\sN_{D_2} \cap D_1$ as the Borel properly exceptional set $\sN_{D_1}$   in Theorem \ref{T:reggreen}(i) for $D=D_1$.
	\end{lem}

	\begin{theorem} \label{T:Rgreen}
	Suppose that  $(\sX,d,m,\sE,\sF )$ is irreducible and  satisfies the ${\rm EHI}_{\rm loc}$. 
		Then $(\sE,\sF)$ has regular Green functions. Moreover, the Borel properly exceptional set $\sN$ of $X$ in the definition of regular Green functions for $(\sE, \sF)$ 
			in Definition 
	\ref{D:goodgreen}(b) has the property that for each $x\in \sX \setminus \sN$ and
	every $t>0$,
$\bP^x(X_t \in dy)$ is absolutely continuous with respect to $m(dy)$. 
	 	\end{theorem}
	
	\pf   Since $(\sX, d, m, \sE,  \sF )$ satisfies the ${\rm EHI}_{\rm loc}$, 
	 	every locally bounded harmonic function is locally H\"older continuous, and $(D, d, m|_D, \sE, \sF^D)$ satisfies 
		the ${\rm EHI}_{\rm loc}$
	for every  non-empty open subset $D$ of $\sX$ whose complement $D^c$ is not $\sE$-polar. 
	The existence of regular Green functions  follow directly from 
	 	Proposition \ref{P:2.1} and Theorem \ref{T:reggreen}(ii) if we can show that there is
		a common properly exceptional set $\sN$ of $X$, independent of $D$,  so that  
		\eqref{e:4.2} of Definition \ref{D:goodgreen} for the regular Green function
		$g_D(x, y)$ on $D$ hold for all $x\in D\setminus \sN$.  
			 
		We first  show
		 the existence of a common properly exceptional set $\sN$. 	 Let $\sX_0$ be a dense countable subset of $\sX$, let $\{(x_n,r_n)\}_{n=1}^\infty$ be an enumeration of   $\sX_0 \times  		   \left((0,\diam(\sX,d)/2) \cap \bQ\right)$, 
		  and set $D_n=B(x_n,r_n)$ for each $n \ge 1$. Then for any $n \ge 1$,  
		 since $\sX\setminus \overline D_n \not=\emptyset$ and hence is non-$\sE$-polar,   
		 we can take a Borel properly exceptional set $\sN'_{D_n}$ of $X^{D_n}$ as in Theorem \ref{T:reggreen}(i) for $D=D_n$. Next, choose a Borel properly exceptional set $\sN'$ of $X$ so that $\bigcup_{n=1}^\infty \sN'_{D_n} \subset \sN'$. Then   for any $n \ge 1$, we can take $\sN' \cap D_n$ as a properly exceptional set $\sN_{D_n}$ of $X^{D_n}$ as in Theorem \ref{T:reggreen}(i) for $D=D_n$. Let $D$ be an open subset of $\sX$ such that $D \subset B(x,R)$ for some $x \in \sX$ and 
		  $R \in (0, \diam(\sX,d)/2)$. Then $D \subset D_n$ for some $n=n_D \ge 1$ and therefore by Lemma \ref{l:monoton} and the definition of $\sN'$ above, we can take $(\sN' \cap D_n) \cap D= \sN'  \cap D$ as the Borel properly exceptional set $\sN_D$ of $X^D$   in Theorem \ref{T:reggreen}(i) for $D$. This concludes the proof that $(\sE,\sF)$ has regular Green functions.
		
		It remains to verify the absolute continuity of $\bP^x(X_t \in dy)$ with respect to $m$ for all $x$ in the complement of an exceptional set. 
		To this end,  
      fix some $D_1=B(x_0, r)$ with $r=\diam(\sX,d)/2$ and let 
		$D_2= {\overline B(x_0, r/2)}^c$, which are non-empty open subsets of $\sX$
		with $D_1\cup D_2 = \sX$.  
		 Define   $\tau_0=0$, $\tau_1:=\tau_{D_1}=\inf\{t>0: X_t\notin D_1\}$,  $\tau_2:= \tau_{D_2}\circ \theta_{\tau_1}=\inf\{t>\tau_1:
		X_t\notin D_2\}$, and for $k\geq 1$,  
		$$
		\tau_{2k+1} :=\tau_{D_1} \circ \theta_{\tau_{2k}} \quad \hbox{ and } \quad \tau_{2k+2} :=\tau_{D_2} \circ \theta_{\tau_{2k+1}}.
		$$ 
Then for every $\alpha>0$, $x\in \sX\setminus \sN'$ and Borel measurable function $f\geq 0$ on $\sX$,
\begin{eqnarray*}
  G_\alpha f(x)  &:=&  \bE^x \int_0^\infty e^{-\alpha t} f(X_t) dt \\
& =& \sum_{k=0}^\infty \bE^x \int_{\tau_{2k}}^{\tau_{2k+1}}  e^{-\alpha t} f(X_t) dt 
+  \sum_{k=0}^\infty   \bE^x \int_{\tau_{2k+1}}^{\tau_{2k+2}}  e^{-\alpha t} f(X_t) dt 
\\
&=&  \sum_{k=0}^\infty \bE^x  \left[ e^{-\alpha \tau_{2k}} \bE_{X_{\tau_{2k}}}    \int_0^{\infty} e^{-\alpha t} f(X^{D_1}_s) ds ; \tau_{2k}<\zeta \right] \\
&& + \sum_{k=0}^\infty \bE^x  \left[ e^{-\alpha \tau_{2k+1}} \bE_{X_{\tau_{2k+1}}}    \int_0^{\infty} e^{-\alpha t} f(X^{D_2}_s) ds ;
\tau_{2k+1}<\zeta  \right].
\end{eqnarray*}
By Theorem \ref{T:5} applied to part processes $X^{D_1}$ and $X^{D_2}$,  respectively, 
there is a properly exceptional set $\sN \supset \sN'$ of $X$ so that 
$G_\alpha (x, dy) \ll m(dx)$ for every $x\in \sX \setminus \sN$.  We then conclude by \cite[Proposition 3.1.11]{CF} that  
		$\bP^x(X_t \in dy) \ll m(dy) $ for every $x\in \sX \setminus \sN$. \qed

	 	We next   give a sufficient condition for a strongly local MMD space  $(\sX,d,m,\sE,\sF )$ to be irreducible.
	First we present a characterization of irreducibility for such a Dirichlet form, which in fact holds also for any
	strongly local quasi-regular Dirichlet forms by using quasi-homeomorphism. 
	See \cite[Theorem 5.2.16]{CF} for a   characterization of irreducible recurrent Dirichlet forms. 
 
	\begin{theorem}\label{T:3.5}
	 Let  $(\sE,\sF )$ be a strongly local regular Dirichlet form on $L^2(\sX; m)$. Then the following are equivalent.
	\begin{description}
	\item{\rm (i)} $(\sE, \sF)$ is irreducible;
	
	\item{\rm (ii)} If $u\in \sF_{\rm loc}$ having $\sE (u, u)=0$, then $u$ is constant  $\sE$-q.e. on $\sX$.

	\item{\rm (iii)}   If $u\in \sF_{\rm loc} \cap  L^\infty (\sX; m)$ having $\sE (u, u)=0$, then $u$ is constant  $\sE$-q.e. on $\sX$.

	\end{description}  
	\end{theorem}
	
	\pf  (i) $\Rightarrow$ (ii):    Suppose $u\in \sF_{\rm loc}$ and $\sE (u, u)=0$. 
	Let $\{U_k; k\geq 1\}$ be an increasing sequence of relatively compact open subsets whose union is $\sX$.
	Then for each $k\geq 1$, there is some $u_k\in \sF$ so that $u_k=u$ $m$-a.e. on $U_k$.
           By Fukushima's decomposition \eqref{e:Fu}, 
           $$
           u_k(X_t)- u_k(X_0) =  M^{u_k}_t+N^{u_k}_t, \quad t\geq 0, \mbox{$\bP^x$-a.s.~for  $\sE$-q.e.~$x\in \sX$,}
	$$
	where $M^{u_k}$ is a martingale additive functional of $X$ having finite energy and $N^{u_k}$ is a continuous additive
	functional of $X$ having zero energy.  Since $\mu_{\<u_k\>}(U_k) =\mu_{\<u \>}(U_k) =0$,  we have
	$M^{u_k}_t =0$ for every $t\in [0, \tau_{U_k})$ $\bP^x$-a.s.~for $\sE$-q.e.~$x \in \sX$ by \cite[Theorems 5.2.3, 5.1.5 and (5.1.22)]{FOT} (and the optional sampling theorem \cite[Problem 1.3.24(i)]{KS} applied to $(M^{u_k})^2- \langle M^{u_k}\rangle$) and 
	$$ 
	\sE (u_k, \varphi ) =0 \quad \hbox{for  every } \varphi \in \sF \cap C_c (U_k).
	$$
	The last display implies by \cite[Theorem 5.4.1]{FOT} that $N^{u_k}_t=0$ for every $t\in [0, \tau_{U_k})$ $\bP^x$-a.s.~for  $\sE$-q.e.~$x \in \sX$. 
	Consequently, we have for each $k\geq 1$ that for  $\sE$-q.e.~$x\in \sX$, $\bP^x$-a.s.~
	$$ 
	 u (X_t)- u (X_0) = u_k(X_t)- u_k(X_0)  = 0 \quad \hbox{for every} t\in [0, \tau_{U_k}).
	 $$
	 As $\lim_{k\to \infty} \tau_{U_k} =\zeta$, we have  for quasi-every $x\in \sX$, $\bP^x$-a.s., 
	 \begin{equation}\label{e:3.7}
	 u (X_t)= u (X_0) \quad \hbox{for every } t\in [0, \zeta).
	 \end{equation}   
	 For $a\in \bR$,  define $A_a=\{x\in \sX: u(x)>a\}$. In view of \eqref{e:3.7}, 
	 $P_t 1_{A_a}  \leq 1_{A_a}$  $m$-a.e. on $\sX$. Hence by the irreducibility of $(\sE, \sF)$, 
	 either  $m(A_a)=0$ or $m(\sX\setminus A_a)=0$. This proves that    $u$ is constant $m$-a.e.  and hence
	 $\sE$-q.e. on $\sX$ by \cite[Theorem 1.3.7]{CF}. 
	
	\medskip
 
	(ii) $\Rightarrow$ (iii) is trivial.  
	
	\medskip
	
	 (iii)   $\Rightarrow$ (i): Were  the Dirichlet form $(\sE, \sF)$ on $L^2(\sX; m)$ not irreducible, there
	 would exist a Borel set $A$ with $m(A)>0$ and $m(A^c)>0$ such that 
	 $1_Au\in \sF$ for any $u\in \sF$ and  \eqref{e:2.1} holds.  
	  In particular,  both $1_A$ and $1_{A^c}$ are in $\sF_{\rm loc}$, and \eqref{e:2.1} with $\one_Au, \one_{A^c}v$ in place of $u,v$ yields 
	 $$
	 \sE (1_Au,  1_{A^c} v) =0 \quad \hbox{for every }  u, v\in \sF.
	 $$ 
	 This together with  \eqref{e:2.7} gives that for any bounded $u, v\in \sF$,
	 \begin{eqnarray*}
	 \int_{\sX} v(x) \mu_{\< 1_A u\>} (dx) 
	 &=& \ 2\sE (1_A u, 1_A u  v) - \sE (1_A u^2,  v) \\
	  &=& \ 2\sE ( u,    (1_A v) u ) - \sE ( u^2, 1_A v) \\
	  &=& \int_{\sX}  (1_A v) (x) \mu_{\<   u\>} (dx) .
	 \end{eqnarray*}  
	 This yields  
	 \begin{equation}\label{e:3.8}
	  \mu_{\< 1_A u\>} (dx) = 1_A (x) \mu_{\<   u\>} (dx).
	 \end{equation}
	 Let $\{U_k; k\geq 1\}$ be an increasing sequence of relatively compact open subsets whose union is $\sX$
	 and $u_k \in \sF\cap C_c (\sX)$ be such that  $u_k=1$ on $U_k$.  
            We have by \eqref{e:3.8} and Proposition \ref{P:2.3}(ii),(i) that 
            $$
            \mu_{\< 1_A  \>} (U_k) = \mu_{\< 1_A u_k \>} (U_k) =  \mu_{\<   u_k\>} (U_k\cap A)=0
            \quad \hbox{for each } k\geq 1,
            $$	 
	 and so 
	 $$
	 \sE (1_A, 1_A) = \frac12  \mu_{\< 1_A  \>} (\sX) =0.
	 $$
	 Since $1_A$, which is in $\sF_{\rm loc}$ and bounded,   
	  is not constant $m$-a.e. on $\sX$ this is a contradiction, and 
	proves that 
	  the Dirichlet form $(\sE, \sF)$ on $L^2(\sX; m)$ is  irreducible.
	 \qed

	 Theorem \ref{T:3.5} in particular implies that irreducibility is invariant under 
	 form-bounded perturbations  in the following sense.
	If   $(\sX, d, m, \sE,\sF )$ and $(\sX, d, \mu, \sE',\sF' )$ are two strongly local regular MMD spaces such that 
	  the Radon measure  $\mu$ does not change $\sE$-polar sets and has full  quasi support on $\sX$,  $\sF \cap C_c(\sX)=\sF'\cap C_c(\sX)$
	   and there is a constant $C\geq 1$ so that 
	\begin{equation}\label{e:4.12} 
	C^{-1} \sE (u, u) \leq \sE' (u, u) \leq C \sE (u, u) \quad \hbox{for } u\in \sF \cap C_c(\sX),
\end{equation}
	then $(\sE, \sF)$ on $L^2(\sX, m)$ is irreducible if and only if so is $(\sE', \sF')$ on $L^2(\sX, \mu)$. 
	Note that, by \cite[Corollary 5.2.12]{CF},
	  the condition that $\mu$ is a smooth measure of $(\sX, d, m, \sE,\sF )$  with full  quasi support on $\sX$
	  ensures that two extended Dirichlet spaces coincide. 
	By \cite[Proposition 1.5.5(b)]{LJ}, \eqref{e:4.12} implies that 
	 	\begin{equation}\label{e:4.13} 
	 	 C^{-1}\mu_{\langle u \rangle} \le \mu'_{\langle u \rangle} \le C \mu_{\langle u \rangle} \quad \hbox{on $\sX$  \,  for any }u \in \sF_e =\sF'_e.
 	\end{equation}

\medskip	
	
	\begin{theorem} \label{T:3.6}
	  Suppose that $(\sX, d)$ is connected. If  a MMD space
	  $(\sX,d,m,\sE,\sF )$ satisfies   (Ha)   and   any bounded function that is harmonic in 
	    each ball $B(x_0, r)$ with $r\in (0, d_{\sX} (x_0) )$  has a continuous modification there, 
	then $(\sX,d,m,\sE,\sF )$ is irreducible.
	\end{theorem} 
	
	\pf    (i) We first show that  for any relatively  compact open subset $U\subset   \sX$ with $\overline{U} \neq \sX$,
	 \begin{equation}\label{e:4.13a}
	 \bP^x(\tau_U< \infty) =1  \quad \hbox{for q.e. } x \in U.  
	 \end{equation} 
	 Indeed, since $\partial U$ is compact, there is some  $\varphi \in \sF\cap C_c(\sX)$ so that $\varphi =1$ in a neighborhood of $\overline U$.
	By the strong locality of $(\sE, \sF)$, $\varphi$ is $\sE$-harmonic in $U$. 
	Define $h(x) :=  \bE^x [ \varphi (X_{\tau_U})]$.
	By \cite[Theorems A.2.6(i), 4.1.3 and 4.2.1(ii)]{FOT}, $h=\varphi$ $\sE$-q.e. on $U^c$. In particular,  
	there is a properly exceptional set  $\sN$ so that for every $x\in U\setminus \sN$, 
	\begin{equation}\label{e:4.14} 
	h(X_{\tau_U}) =\varphi (X_{\tau_U})= 1 \quad \bP^x \hbox{-a.s. on } \{\tau_U<\infty\} . 
	\end{equation}
	Moreover,   $h\in \sF_e$ is $\sE$-quasi-continuous on $\sX$ and $\sE$-harmonic in $U$  by \cite[Theorem 3.4.8]{CF} or \cite[Theorem 4.6.5]{FOT}.
	Consequently, $ f :=\varphi - h$ is a bounded $\sE$-quasi continuous element in $\sF_e$  that is $\sE$-harmonic in $U$ and vanishes $\sE$-q.e. on $U^c$.
	Therefore $\sE(f, f)=0$. By \cite[Lemma 2.2]{C}, there is a properly exceptional set $\sN_1\supset \sN$ of $X$ so that 
\be \label{e:4.15}
\bP^x(f(X_t) = f(X_0) \mbox{ for every $t \ge 0$})=1 \quad \hbox{for  every  } x\in \sX\setminus \sN_1. 
\ee
Note that for $x\in U\setminus \sN_1$, $h(x)=  \bE^x [ \varphi (X_{\tau_U})] = \bP^x (\tau_U<\infty)$. 
For $  x\in U \setminus \sN_1 $ having $h(x)=\bP^x (\tau_U < \infty)>0$,     we have by \eqref{e:4.14} and 
\eqref{e:4.15} that  
$$
h(x)=h(X_0)=   h(X_{\tau_U}) =1 \quad \bP^x \hbox{-a.s. on } \{\tau_U< \infty\}. 
$$	
This proves that $h$ takes values either 0 or 1  $\sE$-q.e. everywhere on $U$.
Consequently, $f=\varphi -h$ takes values  either 0 or 1 $\sE$-q.e. everywhere on $U$ and hence on $\sX$. 
Since $f$ is   bounded  $\sE$-quasi-continuous and $\sE$-harmonic on $\sX$ due to  $\sE(f, f)=0$, 
it is harmonic on $\sX$ on Proposition \ref{P:2.9}(i). Hence by the assumption of the theorem,
 it has a continuous $m$-modification $\wt f$ which takes values either 0 or 1. 
 Since $(\sX, d)$ is connected and $\overline U \not= \sX$,  $\wt f =0$ on $\sX$. 
 It follows that $f=0$ and hence $h=1$ $\sE$-q.e. on $U$. This establishes the claim \eqref{e:4.13a}. 
  
 \medskip
 
	(ii) Suppose $u\in \sF_{\rm loc} \cap L^\infty(\sX; m)$ and $\sE (u, u)=0$.  
	By the proof of (i) $\Rightarrow$ (ii) part of Theorem \ref{T:3.5},
	we know that \eqref{e:3.7} holds $\bP^x$-a.s.~for  $\sE$-q.e.~$x \in \sX$, and in particular that $u$ is   harmonic  on $\sX$.
  By the assumption of the theorem, $ u$ has a continuous version, which we still denote it by $u$. 	 	 
  Since $(\sE, \sF)$ is strongly local,  $1\in \sF_{\rm loc}$ and $\sE (1, 1)=0$. 
	 Let  $x_0$ be an arbitrary point in $\sX$ and  denote $u(x_0)$ by $a_0$. Then $u-a_0 \in \sF_{\rm loc}$ and   
	 $\mu_{\<u-a_0\>}=\mu_{\<u\>}$ by Proposition \ref{P:2.3}(ii).
	 Thus 
	 $$\sE(  u-a_0, u-a_0)=\frac12 \mu_{\<u-a_0\>}(\sX)= \frac12 \mu_{\<u\>}(\sX)= \sE(u, u)=0. 
	 $$
	 Let $v=|u-a_0|$.
	 By \cite[Theorem 4.3.10]{CF},   
	$v\in \sF_{\rm loc}$ and $\sE (v, v)=0$.   
	  By \eqref{e:4.13a}  and Proposition \ref{P:2.9}(ii), $v$ is regular harmonic in $B(x_0,r)$ for any $r \in (0, d_{\sX}(x_0))$.
  	   Since $v(x_0)=0$, by  (Ha) and the continuity of $v$,  we have
            $v(x)=0$ on $B(x_0, r)$ for some $r\in (0, d_\sX (x_0))$;  that is, 
           $u(x)=u(x_0)$ on $B(x_0, r)$ for some $r\in (0, d_\sX (x_0))$. This shows that for any constant $a\in \bR$, 
                       both $A_a:=\{x\in \sX: u(x)>a\}$ and its complement $A_a^c=\{x\in \sX: u(x) \leq a\}$ are  open
           subsets of $\sX$.  If $u$ is not a constant, then there is a constant $a$ so that neither $A_a$ nor $A_a^c$ is empty.
           This would contradict  the assumption that $(\sX, d)$ is connected. So $u$ must be constant.
           This establishes the irreducibility of  $(\sX,d,m,\sE,\sF )$  by Theorem \ref{T:3.5}. \qed

	Combining Theorem \ref{T:3.6} with Theorem \ref{T:Rgreen} shows that if 
	  $(\sX, d)$ is connected and $(\sX,d,m,\sE,\sF )$ satisfies the ${\rm EHI}_{\rm loc}$,
	  then  $(\sE,\sF)$ has regular Green functions.

 \section{Implications of EHI} 	\label{S:5}

Recall the definition of  
  metric doubling from Definition \ref{D:MD}. 
 	We introduce a few related properties.   For a metric space $(\sX,d)$ and $A \subset \sX$, we set $\operatorname{diam}(A)= \sup_{x,y \in A} d(x,y)$.
 	For $x \in \sX, R>0$, we denote by $\overline{B}(x,R):= \{y \in \sX: d(x,y) \le R\}$ the closed ball centered at $x$ of radius $R$.
	The following definition of relative $K$ ball connectedness is adapted from \cite[Definition 5.5]{GH}.
	
	 	\begin{definition} \label{d:metric} \rm 
	\begin{description}
	\item{(i)}
 		Let $K>1$.
			A metric space $(\sX,d)$ is {\em relatively $K$ ball connected} if 
			for each $\eps \in (0,1)$  there exists an integer
			$N=N_\sX(\eps) \ge 1$ such that if $x_0 \in \sX$, $R>0$ and 
			$x, y \in \ol {B}(x_0,R)$ then there exists a chain of balls
			$B(z_i, \eps R )$ for $i=0, \dots, N$ such that $z_0=x$, $z_N =y$, $B(z_i, \eps R) \subset B(x_0, KR)$ 
			for each $i=0,\ldots,N$ and $d(z_{i-1}, z_i) <  \eps R $ for $1\le i \le N$.
			We write $N_\sX$ for the integer $N_\sX(\eps)$ with $\eps = 1/4$. 
			We say also that  $(\sX,d)$ satisfies the property RBC$(K)$.
			We say that  $(\sX,d)$ is {\em relatively ball connected} if 
	there exists $K>1$ such that $(\sX,d)$ is relatively $K$ ball connected.

	\item{(ii)} A metric space $(\sX,d)$ is said to be \emph{uniformly perfect}, if there exists $C>1$ such that if $x \in \sX$,  $r>0$ and $B(x,r)^c \neq \emptyset$	then $B(x,r) \setminus B(x,r/C) \neq \emptyset$.
			
	\item{(iii)}		A metric space $(\sX,d)$ is said be \emph{$L$-linearly connected} (for some $L >1$), 
	 if for all $x,y \in \sX$, there exists a connected compact set $J$ such that $x,y \in J$ and $\diam(J) \le L d(x,y)$.
			
	\item{(iv)}		A \emph{distortion function} is a homeomorphism of $[0,\infty)$ onto itself. 
			Let $\eta$ be a distortion function. 
			A map $f:(\sX_1,d_1) \to (\sX_2,d_2)$ between metric spaces is said to be
			\emph{$\eta$-quasisymmetric} or an \emph{$\eta$-quasisymmetry}, if $f$ is a homeomorphism and
			\[
			\frac{d_2(f(x),f(a))}{d_2(f(x),f(b))} \le \eta\left(\frac{d_1(x,a)}{d_1(x,b)}\right)  
			\]
			for all triples of points $x,a,b \in \sX_1$, $x \neq b$. We say $f$ is a \emph{quasisymmetry} 
			if it is $\eta$-quasisymmetric for some distortion function $\eta$.
			We say that metric spaces $(\sX_1,d_1)$ and $(\sX_2,d_2)$ are \emph{quasisymmetric}, 
			if there exists a quasisymmetry $f:(\sX_1,d_1) \to (\sX_2,d_2)$.
			We say that  metrics $d_1$ and $d_2$  on $\sX$ are \emph{quasisymmetric} (or, $d_1$ is 
			\emph{quasisymmetric} to $d_2$), if the identity map $\operatorname{Id}:(\sX,d_1) \to (\sX,d_2)$ is a quasisymmetry.

	\item{(v)}		 We say a metric space $(\sX,d)$ is \emph{quasi-arc connected}, if there exists a distortion function $\eta:[0,\infty) \to [0,\infty)$ such that for all  pairs of distinct points $x,y \in \sX$, there exists a subset
$J \subset \sX$ and  an $\eta$-quasisymmetry $\gam:[0,1] \to J$ such that 
 $\gam(0)=x$ and $\gam(1)=y$. 
Here $J$ is endowed with the metric $d$ and $[0,1]$ has the Euclidean metric.
\end{description}
	 \end{definition}
	
	The following lemma clarifies some relationships between these conditions.
	
		\begin{lem} \label{l:metric} \rm{
			Let $(\sX,d)$ be a complete metric space.
			\begin{enumerate}[(a)]
				\item 	Assume that $(\sX,d)$ is relatively $K$ ball connected.  Then there exists $L >1$, such that for all $x,y \in \sX$, there exists a continuous map $\gam:[0,1] \to \sX$ such that $\gam(0)=x$, $\gam(1)=y$, 
			 				and  $\diam(\gamma([0, 1])) \le L d(x,y)$.
				In particular, $(\sX,d)$ is $L$-linearly connected.
				
				\item If $(\sX,d)$ is  connected  then it is uniformly perfect. 	
				
				\item  If $(\sX,d)$ is relatively ball connected and satisfies metric doubling, then $(\sX,d)$ is quasi-arc connected.
				
				\item If $(\sX,d)$ is quasi-arc connected, then $(\sX,d)$ is relatively ball connected. 
				
				\item Assume that $(\sX,d)$ is relatively $K$ ball connected. If $\rho$ is a metric on $\sX$ quasisymmetric to $d$, then $(\sX,\rho)$ is also relatively ball connected. In other words, the property of being relatively ball connected is a quasisymmetry invariant.
				
			\end{enumerate}
	}\end{lem}
	
	\proof
	(a) Fix $\eps \in (0,1)$ and let $K, N=N_\sX(\eps)$ be the constants of relative ball connectivity.
	
	Let $x,y \in \sX$ be a pair of distinct points. For each $k \in \bN$, we define $\gam_k:[0,1] \to \sX$ as follows. Let $z^{(1)}_0, z^{(1)}_1,\ldots, z^{(1)}_N$ 
	be a sequence of points in $B(x, Kd(x,y))$ such that $d(z^{(1)}_i,z^{(1)}_{i+1})< \eps d(x,y),$ with $z^{(1)}_0=x, z^{(1)}_N=y$. Let $\gam_1:[0,1] \to \sX$ be a piecewise constant function on intervals defined by
	\[
	\gam_1(t)= z^{(1)}_i \quad \mbox{for all $i=0,1,\ldots, N-1$ and for all $i/N\le t < (i+1)/N$}
	\]
	and $\gam_1(1)= y$. Similarly, for all $i=0,\ldots, N-1$, we choose
	$z^{(2)}_j \in B(z_i^{(1)},K\epsilon d(x,y)), j=iN, iN+1,\ldots, iN+N$ such that 
	$z^{(2)}_{iN}=z^{(1)}_i, z^{(2)}_{iN+N}=z^{(1)}_{i+1}, d(z^{(2)}_j, z^{(2)}_{j+1}) < \eps^2 d(x,y)$ and set 
	\[
	\gam_2(t)= z^{(2)}_j   \quad  \mbox{for all $i=0,1,\ldots,N^2-1$ and for all $i/N^2\le t < (i+1)/N^2$},
	\]
	with $\gam_2(1)=y$. We similarly define $\gam_k:[0,1] \to \sX$ a piecewise constant function on intervals $[j/N^k, (j+1)/N^k)$, $j=0,1,\ldots, N^k -1$. 
	Since for all $t \in [0,1]$, $d(\gam_k(t), \gam_{k+1}(t)) < K \eps^k d(x,y)$, the sequence $\set{\gam_k(t), k \in \bN}$ is Cauchy, and hence converges to say $\gam(t) \in \sX$. This defines a function $\gam:[0,1] \to \sX$.
	Note that
	\[
	d(x,\gam(t)) \le \sum_{k=0}^\infty K \eps^k d(x,y)= K d(x,y)/(1-\eps).
	\]
	If $ \abs{t_1 - t_2} \le \frac{1}{N^k}$ for some $k \in \bN$, we have
	\begin{align*}
	d(\gam(t_1),\gam(t_2)) &\le d(\gam_k(t_1), \gam(t_1)) + d(\gam_k(t_2), \gam(t_2))+ d(\gam_k(t_1),\gam_k(t_2)) \\
	& \le 2 \left( \sum_{l=k}^\infty K \eps^l d(x,y) \right) +  \eps^k d(x,y)\\
	&\le  (2K(1-\epsilon)^{-1}+1)\eps^k d(x,y),
	\end{align*}
	which implies the continuity of $\gam$.
	
	This shows that the image $J=\gam([0,1])$ is a compact, connected set  and $x,y \in J$ 
	with $\diam(J) \le L d(x,y)$, where $L=2K/(1-\eps)$. Therefore $(\sX,d)$ is $L$-linearly connected. 
	
	\noindent (b)
	Let $B(x,r)$ be a ball such that $B(x,r)^c \neq \emptyset$  and let $C>1$. Then $B(x,r)^c$ and $\overline{B(x,r/C)}$ are non-empty disjoint closed sets.
	 Since $(\sX,d)$ is connected, $B(x,r) \setminus B(x,r/C) \neq \emptyset$.
	 
 \noindent	(c) 
	By part (a), $(\sX,d)$ is linearly connected. By Tukia's theorem (\cite[Corollary 1.2]{Mac} and \cite[Theorem 4.9]{TV}), 
	$(\sX,d)$ is quasi-arc connected. 
	
\noindent	(d) Let $\eta$ be the distortion function corresponding to quasi-arc connectedness. Define $K= 1+ \eta(1)$. 
	Let $x,y \in \ol{B}(x_0,R), x \neq y$ and $\gam:[0,1] \to J$ be an $\eta$-quasisymmetry such that $\gam(0)=x, \gam(1)=y$. For all $t \in [0,1]$, 
	\[
	d(x,\gam(t)) \le \eta(t) d(x,y) \le \eta(1) d(x,y).
	\]
 	
	Let $\eps \in (0,1)$ be arbitrary. Let $N \in \bN$ be such that $2 \eta(2/N) \eta(1) < \eps$ and define $z_i= \gam(i/N)$ for $i=0,1,\ldots,N$.
	By $\eta$-quasisymmetry, we have 
	\[
	d(z_i,z_{i+1}) \le \eta(2/N) d(z_i,w) \le \eta(2/N) \eta(1) d(x,y) \le 2  \eta(2/N) \eta(1)R < \eps R, 
	\]
	where $w=x$ if $i \ge N/2$ and $w=y$ otherwise.
	This implies that $(\sX,d)$ is relatively $K$ ball connected, where $K=2(1+\eta(1))$.	\\
\noindent	(e) Let $\on{Id}:(\sX,d) \to (\sX,\rho)$ be an $\eta$-quasisymmetry, where $(\sX,d)$ is relatively $K$ ball connected.
	Let $\eps \in (0,1)$ and let $x,y \in \sX, x \neq y$ be arbitrary. Choose $\eps' \in (0,1)$ such that
	\[
	\eta(2\eps') \left( \eta(K)+1 \right) < \eps.
	\]
	Choose points $z_0,z_1,\ldots,z_N$ such that $z_0=x, z_N=y, B(z_i, \eps' d(x,y)) \subset B(x,K d(x,y))$ and $d(z_i,z_{i+1})< \eps' d(x,y)$, where $N=N_{(\sX,d)}(\eps')$ is the constant associated with the relative ball connected property of $(\sX,d)$. For any $i=0,1,\ldots, N-1$, let $w \in \set{x,y}$ be such that $d(z_i,w) = \max(d(z_i,x), d(z_i,y))$. Since $d(x,y)/2 \le d(z_i,w)$, we obtain
	\begin{align*}
	\rho(z_i,z_{i+1}) &\le  \eta(d(z_i,z_{i+1})/d(z_i,w)) \rho(z_i,w)  
	\le \eta(2 \eps') (\rho(x,z_i) + \rho(x,y)) \\
	& \le \eta(2\eps') \left( \eta(K)+1 \right)\rho(x,y) < \eps \rho(x,y).
	\end{align*}
	Since $\rho(x,z_i) \le \eta(K) \rho(x,y)$, $(\sX,\rho)$ is relatively $K_\rho$ ball connected, where $K_\rho =2+ 2\eta(K)$.
	\qed
	
	\begin{remark} 
		{\rm 
			See \cite[Definition 5.5]{GH} for the definition of {\em relatively $(\eps,K)$ ball connected}.  
			It is immediate that if $(\sX,d)$ is  relatively $K$ ball connected then it is 
			relatively $(\eps,K)$ ball connected for any $\eps \in (0,1)$. Conversely
			it is straightforward to show that if for some $\eps \in (0,1)$, $K>1$
			$(\sX,d)$ is  relatively $(\eps,K)$ ball connected then it is 
			relatively $K'$ ball connected with $K' = 1+ \frac{K}{1-\eps}$.
	}\end{remark}

The main result of this section is the following.

	\begin{thm} \label{T:ehitomd}
		Let $(\sX,d)$ be a metric space such that $\overline{B(x,r)}$ is compact for all $x \in \sX$ and $r>0$. Assume that   $(\sX,d,m,\sE,\sF)$ is a MMD space that satisfies the EHI. 
		The following are equivalent: \\
		{\rm (a)}  $(\sX,d)$ is relatively $K$ ball connected for some $K>1$.   \\
		{\rm (b)}   $(\sX,d)$ satisfies metric doubling.\\
		{\rm (c)}    $(\sX,d)$ is quasi-arc connected.
	\end{thm}
	
	\proof  (b) $\Rightarrow$  (a). This follows by the argument in \cite[Proposition 5.6]{GH}:
for any $K>1+\delta_H^{-1}$ we obtain relative $K$-ball connectedness. (The hypothesis
of volume doubling there is only used to obtain metric doubling). 
 Note that this implication requires EHI only for $[0,1]$-valued $h \in \sF \cap C_c(\sX)$.

   (c)  $\Rightarrow$  (a),  (a) $+$ (b)  $\Rightarrow$  (c) (and so   (b) $\Rightarrow$  (c)) 
     are proved in Lemma \ref{l:metric}.
 
   The proof of (a) $\Rightarrow$  (b) needs more preparation and will be given  after Lemma \ref{L:3.10}. 
\qed

\sms
 
For the proof of (a) $\Rightarrow$  (b) in Theorem \ref{T:ehitomd}, 
we follow \cite[Section 3]{BM1}; however it was assumed 
there that the metric $d$ was geodesic, and some changes are needed to handle
the case when we only have that $(\sX,d)$ is relatively $K$ ball connected.  We now outline these changes.

\begin{definition} \label{D:HG}
		{\rm We say that $(\sX,d,m,\sE, \sF)$ satisfies the condition (HG)
		if $(\sX,d,m,\sE, \sF)$ has regular Green functions and 
		 there exist constants 
			$C_G, K_G>0$ such that for any $x_0 \in \sX$, $R>0$ and open set $D$ in $\sX$ with
			$B(x_0,K_G R) \subset D$ and $D^c$ non-$\sE$-polar
			\be
			\sup_{y_2 \in D \setminus B(x_0,R)} g_D(x_0, y_2) \le C_G 
			 \inf_{y_1 \in \ol {B(x_0, R)} \setminus \{x_0\}}
			  g_D(x_0, y_1).
			\ee
		}\end{definition}
	
	\medskip

	\begin{assumption}  \label{A:main}
		{\rm 
			Throughout the remainder of this section
		 			except for Lemma \ref{L:regpath1}, 
			 we assume that $(\sX,d)$ is a metric space such that $\overline{B(x,r)}$ is compact for all $x \in \sX, r>0$ and 
			 that $(\sX,d)$ is relatively $K$ ball connected for some $K\ge 2$. 
			Furthermore  we assume that the MMD space $(\sX, d, m, \sE, \sF)$
			satisfies the (scale invariant) EHI with constants $C_H, \delta_H$. 
	} \end{assumption}
	 	
	\medskip
	
	 	Recall that by Lemma \ref{l:metric}(a),  a complete metric space  $(\sX, d)$ that is relatively $K$ ball connected  is connected. 
	Thus under Assumption \ref{A:main},  $(\sX, d, m, \sE, \sF)$ is irreducible by Theorem \ref{T:3.6}  and
	has regular Green functions by Theorem \ref{T:Rgreen}. 
	 	 	By the maximum principle \eqref{e:max} for the regular Green function 
			$g_D$ in Theorem \ref{T:reggreen}(ii),   
	 	we have for any $B=B(x_0, R)\Subset D$, 
	\be
	\sup_{D \setminus B} g_D(x_0, \cdot ) =   \sup_{\pd B} g_D(x_0, \cdot ),
	\qquad   
	 	\inf_{ \ol B \setminus \{x_0\}}
	g_D(x_0, \cdot )  =  \inf_{\pd B} g_D(x_0, \cdot ).
	\ee

	\begin{proposition} \label{P:HG}
		 		Let $(\sX,d,m,\sE,\sF)$ be a MMD space that satisfies Assumption \ref{A:main}.
	 Then (HG) holds
		with constants $C_G, K_G$ for $K_G=2K+1$, where $C_G$ depends only on $C_H, \delta_H, K_G, N_\sX$.
	\end{proposition}
	\proof
 	This follows from the proof of \cite[Lemma 5.7]{GH}. (The statement of the 
	result in \cite{GH} has stronger hypotheses, but these are only used to obtain the existence
	and regularity of the Green function, and to prove that $(\sX,d)$ is relatively $(\eps, K)$ ball connected
	for some $\eps \in (0,1)$ and $K>1$.)  
	\qed
	
	Under Assumption \ref{A:main}, 
	$(\sE,\sF)$ has regular Green functions by Theorems \ref{T:3.6} and \ref{T:Rgreen}, and 
	(HG) holds with constant $K_G=2K+1$ and $C_G>1$ by Proposition \ref{P:HG}.

	\begin{corollary}  \label{C:3.2}
		(See \cite[Corollary 3.2]{BM1}.) 
			Let $(\sX,d,m,\sE,\sF)$ be a MMD space that satisfies Assumption \ref{A:main}.
		Let $K_1  =K+1$. 
 For any $\delta \in (0, 1/2]$,  there exists a positive constant $C$ that depends only on $\delta$ and the constants in Assumption \ref{A:main}  such that the following holds:
		for  any open set $D$ in $\sX$ whose complement $D^c$ is non-$\sE$-polar and for any 
		$B(x_0, K_1 R )\subset D$, 
		$$
		 g_D(x_0,x) \le C g_D(x_0,y)  \quad \hbox {for } x,y \in B(x_0,R) \setminus B(x_0, \delta R) .
		$$
	\end{corollary}
	
	\proof
Let $x,y \in B(x_0,R) \setminus B(x_0, \delta R)$.  Let $\eps =  \delta \delta_H / ( 1+ \delta_H)$; we have $\eps <  \delta /2$.
We connect $y$ to $x_0$ by a chain of balls $B(z_i, \eps R)$, $i =0, 1, \dots N$,
with the properties given in  
Definition \ref{d:metric}(i) of  relatively $K$ ball connected. 
Let $i_0$ be the first integer such that $d(z_{i_0}, x_0) < \delta R$.  
With the definition of $\eps$ given above, $g_D(x_0, \cdot)$ is harmonic
on $B(z_i, \eps R/\delta_H)$ for $i=0, \dots, i_0-1$, and so we can use the 
EHI to deduce that $g_D(x_0, z_{i_0}) \le C_H^N g_D(x_0,y)$. 
Finally by (HG) we have  $g_D(x_0,x) \le C_G g_D(x_0,z_{i_0})$. \qed

	\begin{lemma} \label{L:3.3}
		(See \cite[Lemma 3.3]{BM1}.)
			Let $(\sX,d,m,\sE,\sF)$ be a MMD space that satisfies Assumption \ref{A:main}.
		Let $K_2 = 2K+3$. There exists  
				$C_0>1$ that depends only on the constants in Assumption \ref{A:main}
		such that the following holds:
		Let $x_0 \in \sX$, $R>0$ and let $B(x_0,K_2R)\subset D$, where $D$ is an open set in $\sX$ such that $D^c$ is non-$\sE$-polar. 
		Then 
		 		if $x_1,x_2, y_1, y_2 \in B(x_0, R)$ with
		$d(x_j,y_j) \ge R/4$,  then
		\be
		g_D(x_1, y_1) \le C_0 g_D(x_2,y_2). 
		\ee 
	\end{lemma}
	
	\proof 
 Note that for any four numbers $a_i\in [0, R)$, $1\leq i\leq 4$, 
$$ 
\left| (0, R)\setminus \cup_{i=1}^4 [a_i-(R/9), a_i+(R/9)] \right|\geq R-(8R/9)=R/9>0
$$
so there is some $a_0\in (0, R)$ so that  $|a_0-a_i|>R/9$ for all $1\leq i\leq 4$.
	As $d(x_0,\cdot)$ is continuous and by the RBC$(K)$ property and Lemma \ref{l:metric}(a),
$B(x_0,R)$ contains a connected path from $x_0$ to $B(x_0,R)^c$, we have
$\{ d(x_0,x): x \in B(x_0, R ) \} = [0,R)$.
Thus    there is $a_0\in (0, R)$ so that   
 $|a_0 - d(x_0,w)|>R/9$ for $w \in \{ x_1,x_2,y_1, y_2\}$. 
Let $z\in B(x_0, R)$ having $d(x_0, z)=a_0$. 
Now applying Corollary \ref{C:3.2} to the balls $B(x_1, 2R)$, $B(z, 2R)$ and $B(x_2, 2R)$ with $\delta = 1/18$
consecutively, we get  by the symmetry of the Green function $g_D(x, y)$ that 
$$
g_D(x_1, y_1) \leq C g_D(x_1, z) \leq C^2 g_D(x_2, z) \leq C^3 g_D(x_2, y_2).
$$
This establishes the lemma by taking $C_0=C^3$. 
  \qed

 	 As in \cite{BM1},  we define for an open set $D\subset \sX$ with non-$\sE$-polar complement:
	\begin{align*}
	g_D(x,r) &= \inf_{y \in \partial B(x,r)} g_D(x,y) \quad \mbox{ provided } \overline{B(x,r)} \subset D, \\
	\operatorname{Cap}_D(A) &= \inf\{ \sE(f,f): f \in (\sF_e)^D, \, f \ge 1 
	 	\   \sE \hbox{-q.e.  on } A\}, \quad A\subset D,
	 	\end{align*}
 	where $(\sF_e)^D= \{ u \in \sF_e: u=0~\sE-\mbox{q.e.~ on $\sX \setminus D$} \}$. By noting that $(\sF_e)^D$ is equal to the extended Dirichlet space of $(\sE,\sF^D)$ \cite[Theorem 3.4.9]{CF}, we set $\sF_e^D:= (\sF_e)^D$.
	 	We call $\operatorname{Cap}_D(A)$ the relative capacity of $A$ in $D$. 
	 	 	The maximum principle \eqref{e:max} 
	implies that $g_D(x,r)$ is non-increasing in $r$, and 
	an easy application of (HG)  gives that if
 $y \in \partial B(x,r)$ and $\overline{B(x,  r)}\cup B(y, K_Gr) \subset D$ then
	\be
	g_D(x,r) \le C_G g_D(y,r). 
	\ee
	
	Given Proposition \ref{P:HG} the proof of the next Lemma is the same as in  \cite[Lemma 3.5]{BM1}.

	\begin{lemma} \label{L:3.5}
	 			Let $(\sX,d,m,\sE,\sF)$ be a MMD space that satisfies Assumption \ref{A:main}.
		There is a constant $C_G>0$ depending only on the constants in Assumption \ref{A:main}
		such that for any    open set $D$ whose complement  $D^c$ is non-$\sE$-polar and for any  
		$B(x_0, K_G r) \subset D$ where  $K_G=2K+1$, 
		\be
		g_D(x_0,r) \le     \textstyle \operatorname{Cap}_D(B(x_0,r))^{-1}  \le  C_G  g_D(x_0,r).
		\ee 
	\end{lemma}
	
	\ms 
	
	\begin{remark} \label{r:nonpolar} 
\rm For any $x \in \sX$, $0<R<\diam(\sX,d)/2$, the complement of  ball $B(x,R)^c$ is non-$\sE$-polar.  This is because  by the triangle inequality, there exist  $z \in \sX$ and $0<r<\diam(\sX,d)/2 -R$ so that  $B(z,r) \subset B(x,R )^c$. Since $m$ has full support, $m\left(B(x,R)^c\right) \ge m(B(z,r)) >0$ and thus $B(x,R)^c$  has positive capacity.
	\end{remark}

 \begin{lem}(See \cite[Lemmas  7.1 and 7.4]{GH} and \cite[Lemma 2.5]{GNY})  \label{l:nashwill}
 	Let $(\sX,d,m,\sE,\sF)$ be a MMD space that satisfies Assumption \ref{A:main}.
 There exists $C$ such that for any sufficiently large $A>1$, and for any ball $B(x,r), n \in \bN$ with $A^nr<\diam(\sX,d)/A$, denoting 
 	$B_k=  B(x, A^kr)$, 
 	we have
 	\[
 	\sum_{i=0}^{n-1} \operatorname{Cap}_{B_{i+1}}(B_i)^{-1} \le \operatorname{Cap}_{B_n}(B_0)^{-1} \le C\sum_{i=0}^{n-1} \operatorname{Cap}_{B_{i+1}}(B_i)^{-1}.
 	\]
 \end{lem}
 
 \proof The upper bound is contained in   \cite[Lemmas  7.1 and 7.4]{GH}. The upper bound in \cite{GH} is stated under the additional volume doubling property assumption but the proof only uses the properties of the regular Green functions and (HG). Also, the constant $C$ can be chosen to be $C_G$ for any $A \ge K_G$, where $C_G, K_G$ are as given in Definition \ref{D:HG}.
 
 The lower bound is a general fact that does not require the EHI -- see \cite[Lemma 2.5]{GNY}.
 \qed

	\begin{lemma} \label{L:3.7}
			Let $(\sX,d,m,\sE,\sF)$ be a MMD space that satisfies Assumption \ref{A:main}. 
		Let $B=B(x_0,R) \subset \sX$, and  
		$B_1= B(x_1, R/(8K_3))$ with $x_1  \in B(x_0, R/(4K_3))$, where $K_3=K+2$ and $K$ is as given in Assumption \ref{A:main}.
		There exists $p_0 >0$ depending  only on the constants in Assumption \ref{A:main} such that
		\be \label{e:h11}
		\bP^y( \sigma_{B_1} < \tau_B ) \ge p_0> 0  \quad \hbox{for $\sE$-q.e.  } y \in B(x_0, 
	 		 R/ (2K_3) ).
		\ee
 	\end{lemma}
	
	\proof  We consider two cases.

   (i) Suppose $B(x_0,R)^c$ is non-$\sE$-polar.  
 	By the maximum principle it is enough to prove this for $y \in \pd             B(x_0, R/ (2K_3) )$.
The argument, which uses Corollary \ref{C:3.2}, is the same as in 
	 	\cite[Lemma 3.7]{BM1}. 
	
   (ii) Now suppose $B(x_0,R)^c$ is $\sE$-polar. 
 By Remark \ref{r:nonpolar}, $R \ge \diam(\sX,d)/2$ and $\diam(\sX,d)<\infty$. If 
    $R> 8K_3 \diam(\sX,d)$, 
 then $B_1=\sX$ and \eqref{e:h11} is obviously true. Therefore, it suffices to consider the case when 
   $R \le  8K_3 \diam(\sX,d) < \infty$.

Let  $x_0 \in \sX,x_1 \in B(x_0,R/(4K_3)),y \in  B(x_0,R/(2K_3))$.
 Let 
   $\eps=1/(130 K_3^2)$ and let $B(z_i,\frac{\eps R}{4K}), 0 \le i \le  N:=N_{\sX}(\eps)$ 
 be a chain of balls with $z_0=y$, $ z_N=x_1$,    $B(z_i,\frac{\eps R}{4K}) \subset
    B(x_0,R/4 )$ 
 for all each $i$ and 
 $d(z_{i-1}, z_i)<\eps R/(4K)$ for $1\leq i\leq N$ as in Definition \ref{d:metric}(i) with $R/(4K)$ in place of $R$ there. 
 Since $8 K_3 \eps R
   \leq 64 K_3^2 \eps \diam(\sX,d)
 < \diam(\sX,d)/2$, by Case 1 and Remark \ref{r:nonpolar}, we have
	$$ \bP^w ( \sigma_{ B(z_i, \eps R)} < \tau_{ B( z_{i-1}, 8 K_3 \eps R) } ) \ge p_0
\quad \hbox{for $\sE$-q.e.~}  w \in  B( z_{i-1}, 
 4\eps R). $$
Since $B(z_i,8K_3 \eps R) \subset B(x_0,
 R/2)$ 
for all $i$, using the strong Markov property, we conclude that
$$
	\bP^w( \sigma_{B_1} < \tau_B ) \ge p_0^N > 0  \quad \hbox{for $\sE$-q.e.~ } w \in B(y, 4 \epsilon R).
$$
Since $B(x_0,R/(2K_3))$ is covered by a countable family of balls $B(y,4\epsilon R), y \in B(x_0,R/(2K_3))$, by replacing $p_0$ by $p_0^N$, we obtain \eqref{e:h11} in the second case as well.
	\qed
	
	\begin{remark} \label{R:chainprob}  \rm \begin{enumerate}[\rm (i)]
	\item 		In \cite[Lemma 3.7]{BM1}, 
		the corresponding result held for $y \in B(x_0, 7R/8)$; we cannot expect that
			here, since such a point $y$ might not be connected to $B_1$ by a path inside $B$.

	\item Let $B_i= B(z_i, \eps R)$, $0 \le i \le n$ be a chain of balls as in Definition \ref{d:metric}(i). 
			Using this Lemma we have for each $i$
			$$ \bP^y( \sigma_{ B(z_i, \eps R) }< \tau_{ B( z_{i-1}, 8 K_3 \eps R) } ) \ge p_0
			\quad \hbox{ for $\sE$-q.e.~} y \in  B( z_{i-1}, 4\eps R). $$
			Thus if
			$$ D = \cup_{i=0}^n B( z_i, 8  K_3\eps R), $$
			then 
			\be \label{e:chainmove}
			\bP^y (  \sigma_{B_n} < \tau_D ) \ge p_0^n  \, \mbox{ for $\sE$-q.e.} \, y \in B(z_0, 4 \eps R).
			\ee
			\end{enumerate} 
	 \end{remark}
	
	\begin{corollary} (See \cite[Corollary 3.8]{BM1}). \label{C:3.8}
		Let $(\sX,d,m,\sE,\sF)$ be an MMD space that satisfies Assumption \ref{A:main}.
		Let $B(x,R) \subset D$, where $D$ is an open set in $\sX$ and $D^c$ is non-$\sE$-polar. 
 		There exist positive constants  $c$ and $\theta$ that depend only on the constants in Assumption \ref{A:main}
		such that if 
		$0< s < r < R/(K+1)$ then
		\be
		\frac{ g_D(x,r) }{ g_D(x,s) } \ge  c \Big(  \frac{ s}{ r } \Big)^\theta.
		\ee
 	\end{corollary}
	
	\proof This follows easily from Corollary \ref{C:3.2}. \qed

	The following Lemma is used to regularize chains of balls obtained 
	by using the RBC$(K)$ property.

	\begin{lemma}  \label{L:regpath1}  
		Suppose that $(\sX,d)$  is a metric space satisfying the RBC$(K)$ property.
		Let $d(x,y)<R$, $\eps \in (0,1)$ and $ \eps R <r < R$. 
		There exists a chain of balls  $B(z_i, \eps R)$, $0\le i \le n$ with the following  properties: \\
		{\rm (i)}  $z_0=x$, $z_n=y$ and $d(z_{i-1} , z_i ) < \eps R$ for $1\le i \le n$; \\
		{\rm (ii)}  $B(z_i, \eps R) \subset B(x,KR)$ for $0\le i \le n$;  \\
		{\rm (iii)} If $j = \max\{i: z_i \in B(x,r)\}$ then $B(z_i, \eps R) \subset B(x,Kr)$ for $0\le i \le j$;  \\
		{\rm (iv)}  $n \le N_\sX(\eps) + N_\sX(\eps R/r)$.
	\end{lemma}
	
	\proof By the RBC$(K)$ property there exists a chain of balls
	$B(w_i, \eps R)$, $0\le i\le m_1$ connecting $x$ and $y$ and satisfying the conditions of Definition \ref{d:metric}(i) with $x_0=x$. 
	Let $k= \max\{i: w_i \in B(x,r)\}$. By the RBC$(K)$ property   for $x$ and $w_k$, and with
	$\eps$ replaced by $\eps' = \eps R/r$, there exists a chain $B(w'_i, \eps R)$, $0\le i \le m_2$
	with $B(w'_i , \eps R) \subset B(x,Kr)$. Joining the paths 
	$w'_0, \dots, w'_{m_2}$ and $w_{k+1}, \dots w_{m_1}$ gives a path $(z_i)$ which
	satisfies the conditions (i)--(iv). \qed

	\begin{lemma}  \label{L:3.10} 
		 Let $(\sX,d,m,\sE,\sF)$ be a MMD space that satisfies Assumption \ref{A:main}.
		There exists an integer $N\geq 1$ that depends only on the constants in Assumption \ref{A:main}
		such that if $x_0 \in \sX$, $R>0$ and $B(z_i, R/8)$, $1\le i \le m$, 
		are disjoint balls with $z_i \in B(x_0, R) \setminus  B(x_0, R/2)$,
		 then $m \le N$.  
		 	\end{lemma}
	
	\proof 
 	This lemma corresponds to \cite[Lemma 3.10]{BM1}. 
	In \cite{BM1} the metric $d$ on $\sX$ is assumed to be a geodesic distance, and the 
	proof in \cite{BM1}	 uses this property quite strongly. The proof here
	is much longer since we only have the weaker property that $(\sX,d)$
	is relatively $K$ ball connected.
	
	Let $(z_k, 1\le k \le m)$ satisfy the hypotheses of the Lemma, and write $B_k = B(z_k, R/8)$. 
	Choose $\eps = 1/(720 K^2)$, and let $n= N_\sX(\eps) + N_\sX(24K\eps)$.
	For each $k$ we use Lemma \ref{L:regpath1} with $r= R/(24K)$ to find a 
	chain of balls $B(w_{ki}, \eps R)$ with $0\le i \le n$ connecting $z_k$ to $x_0$. 
	Note that by taking $w_{ki} =x_0$ for large $i$ if necessary, we can assume that
	all the chains have length $n$.
	We set $l_k = \max\{ i: w_{ki} \in B(z_k,R/(24K)) \}$, and write $z'_k = w_{k l_k}$.
	
	\sms
	We now find a subset $I$ of the balls $B_i$ 
	such that the chain $(w_{ik}, 0\le k \le n)$ associated with one ball does not hit any other ball 
	with index in $I$.
	
	For $1\le i,j \le m$ set $a_{ij} =1$ if $\{ w_{ik}, 1\le k \le n\} \cap B_j \neq \emptyset$, and let
	$a_{ij}=0$ otherwise. Let $b_j = \sum_i a_{ij}$. 
	Since each $w_{ik}$ is in at most one ball $B_j$, we have
	$\sum_j a_{ij} \le n$, and hence $\sum_j b_j = \sum_i \sum_j a_{ij} \le mn$.
	Thus if  $J=\{j: b_j \le 2n\}$ then $|J| \ge m/2$. 
	
	We now consider the collection of balls $(B_i, i \in J)$, and relabel them
	$B_1, \dots, B_{m_1}$ where $m_1 = |J| \ge m/2$.
	We now start with the ball $B_1$, and remove from the collection of balls
	$B_2, \dots , B_{m_1}$ any ball $B_j$ such that either $a_{1j}=1$ or $a_{j1}=1$.
	Since $1 \in J$ and $a_{11}=1$, less than  $3n$ balls are removed. Set $j_1=1$. 
	Let $j_2$ be the smallest label of a ball which has not been removed; we
	now repeat the procedure above with this ball, and remove any ball $B_i$ such that
	$i > j_2$ and $a_{j_2 i} + a_{i j_2} \ge 1$. We continue until there are no balls
	left, and write $I=\{ j_k, 1\le k \le m' \}$ for the set of balls which are retained. 
	Since at each step we remove at most $3n-1$ balls, we have
	$3 n m' \ge \half m$. 
	
	By the construction above we have that
	$$ w_{i k} \not\in \cup_{j \in I \setminus \{i\}} B_j \q \hbox{ for } i \in I,\,  k=0, \dots, n.$$
	For $i \in I$ set $B'_i =  B(z_i, \eps R)$, $A_i = B( z'_i, \eps R)$,
	and let
	$$ D = B(x_0, 2 K R ) \setminus \cup_{i \in I} B'_i. $$
	
	We now claim that  for $i \in I$
	\begin{align} 
	\label{e:hitB'}
	&\bP^y(  \sigma_{B'_i} < \tau_{B_i} ) \ge p_0^n \q \hbox{ for $\sE$-q.e.  } y \in 4A_i:=B(z_i',4 \eps R),  \\
	\label{e:hitA}
	&\bP^{y}( \sigma_{A_i} < \tau_D ) \ge p_0^n \q \hbox{ for $\sE$-q.e. }  y \in B(x_0,4 \eps R). 
	\end{align} 
	Both these inequalities follow by chaining the bound 
	in Lemma \ref{L:3.7}, as in Remark \ref{R:chainprob}(ii), along a sequence of balls.
	For \eqref{e:hitB'} we use the sequence $B(w_{ij},\eps R), 0\le j \le  l_i$ (starting at $j=l_i$ and ending at $j=0$), and for
	\eqref{e:hitA} we use  $B(w_{ij}, \eps R), l_i \le j \le n$.
	(We start at $j=n$ and end at $j=l_i$.)
	
	The remainder of the proof is as in \cite[Lemma 3.10]{BM1}.
	Let $F_i = \{ \sigma_{A_i} < \tau_D\}$, and 
	$$Y = \sum_{i \in I} 1_{F_i} $$
	be the number of distinct balls $A_i$ hit by
	$(X_t, 0 \le t \le \tau_D)$. The bound  \eqref{e:hitB'} implies that if $X$ hits
	$A_i$ then with probability at least $p_0^n$ it leaves $D$ before it hits
	any other ball $A_j$ with $j \neq i$. Thus $Y$ is stochastically dominated by a geometric r.v.
	with mean $p_0^{-n}$, and so
	$$ 
	\bE^{y} Y \le   p_0^{-n} 	\quad \hbox{ for $\sE$-q.e. }  y \in B(x_0, 4 \epsilon R).    
	$$
	However, by \eqref{e:hitA} we also have 
	$$ \bE^{y} Y = \sum_{i\in I} \bP^{y}(F_i) \ge |I| p_0^n = m' p_0^n
	\quad   \hbox{ for $\sE$-q.e. }  y \in  B(x_0, 4 \eps R).
	$$
  Since $m' \ge \frac{m}{6n}$ as given above, it follows that 
	$m \le 6n p_0^{-2n}$. \qed
	
	\bigskip
	
	 	Now we can finish the proof of Theorem  \ref{T:ehitomd} by giving the 
	
	\med{\em Proof of (a) $\Rightarrow$ (b) in Theorem \ref{T:ehitomd}.}
	 (i) Suppose  that a metric space $(\sX, d)$ has the property that there is an integer $N'\geq 1$, 
	so that any ball $B(x, R)$ contains at most $N'$ points
	that are at distance of at least $R/2$. 
	Given any ball $B(x, R)\subset \sX$, 
	take $z_1 \in B(x, R)$, $z_2\in B(x, R)\setminus B(z_1, R/2)$, 
	and for $k\geq 3$,  $z_{k } \in  B(x, R)\setminus \cup_{j=1}^{k-1} B(z_j, R/2)$ if the set is non-empty. 
	By the assumption, we can only proceed this procedure up to  some number $k_0$ no larger than $N'$.
	Clearly $\cup_{j=1}^{k_0} B(z_j, R/2) \supset B(x, R)$. Thus  $(\sX, d)$ is metric doubling.
	Conversely, suppose $(\sX, d)$ is metric doubling with positive integer $N\geq 2$ in Definition 
 	\ref{D:MD}. 
	For any ball $B(x, R)$, by applying the definition of (MD) to $B(x, R)$ and to balls with radius $R/2$, and then with radius $R/4$ (to guarantee $x_1,\ldots,x_{N^3} \in B(x,R)$)
	there are  $N^3$  points $x_1, \dots, x_{N^3}$ in $B(x, R)$ so that $\cup_{j=1}^{N^3} B(x_j, R/4)
	\supset B(x, R)$.  Suppose $\{z_1, \dots, z_n\} $  are $n$ points in $B(x, R)$
	that are at distance of at least $R/2$, then  $\abs{\{z_1,\ldots,z_n\} \cap B(x_k,R/4)} \le 1$ for any $1 \le k \le N^3$.
	Thus $n\leq N^3$. This proves that   a metric space $(\sX, d)$  is (MD) if and only
	if there is some constant $N'$ so that any ball $B(x, R)$ contains at most $N'$ points
	that are at distance of at least $R/2$ from each other. 
  
          (ii) Now let $N\geq 1$ be the integer in Lemma \ref{L:3.10}.
           Let $x_0 \in \sX$, $R>0$, and let $z_i \in B(x_0,R)$, $1\le i \le n$, with the property 
	that the balls $B(z_i, R/8)$ are disjoint. 
	By Lemma \ref{L:3.10} applied first to $B(x_0, R)$ and then
	to $B(x_0, R/2)$, there are at most $2N$ of the $z_i$ in $B(x_0, R) \setminus B(x_0, R/4)$.
	Using the relative $K$-ball connectivity of $\sX$, we can find $x_1$ such that
	$R/2 \le d(x_0, x_1) < 3R/4$ (here we assume without loss of generality that $B(x_0,R/2) \neq \sX$; otherwise we can cover $B(x_0,R)$ with $B(x_0,R/2)$). 
	Thus 
 	$B(x_0, R/4) \subset B(x_1, R) \setminus B(x_1, R/4)$.
	So  by Lemma  \ref{L:3.10} applied  to $B(x_1, R)$, there are at most $2N$ points $z_i$ in $B(x_0, R/4)$.
	Consequently,  $n \le 4N$.  This proves that $(\sX, d)$ is (MD) in view of its equivalent characterization given in (i).  \qed

	\ms
	We need to compare the Green functions in two concentric balls.
	
	\begin{lemma} \label{L:3.12}
		(See \cite[Lemma 3.12]{BM1}.)
			Let $(\sX,d,m,\sE,\sF)$ be a MMD space that satisfies Assumption \ref{A:main}.
 		There exists a constant $C_1$ that depends only on the constants in Assumption \ref{A:main}
		such that if 
		$B=B(x_0, R)$, $2B = B(x_0,2R)$ and $B(x_0,(2+ 1/(128K^2))R)^c$ is non-empty, then
		$$  g_{2B}(x,y) \le C_1 g_{B}(x,y) \quad \hbox{ for } x,y  \in 
 	           B(x_0, R/(8K)),  x \neq y. $$
 	\end{lemma}
	
	\proof 
	 	Let $a_1 = 1/(8K)$, $a_2 = 1/(4K)$, $a_3 = 1/(2K)$, $\eps = 1/( 128K^2)$ 
	and $B_i = B(x_0, a_i R)$.
	Let $p_1>0$. Suppose that there exist an open set $D$ with
	$\overline{B_2} \subset D \subset B$ and a Borel set $A \subset \sX$ such that 
	\begin{align}
	\label{e:xhitA}
	&\bP^{x}( X_{\tau_{D}} \in A) \ge p_1, \hbox{ for $\sE$-q.e.~} x \in  B_2, \\
	\label{e:Ato2B}
	&\bP^w( \tau_{2B} < \sigma_{B_2} ) \ge p_1 \hbox{ for $\sE$-q.e.~} w \in A. 
	\end{align} 
	 Let $y \in B_1$.
	Let $x_1 =x_1(y) \in \pd B_2$ be chosen to maximize $g_{2B}(x', y)$ for $x' \in \pd B_2$.
	Write $h(w) = \bP^w( \tau_{2B} < \sigma_{B_2})$. 
	Then by the strong Markov property at $\tau_D,\sigma_{B_2}$ and the occupation density formula  \eqref{e:4.2}, for $\sE$-q.e.~ $z \in D \setminus\{y\}$
	\begin{align*}
	g_{2B}(z,y) &= g_{D}(z,y) + \bE^{z} g_{2B}(X_{\tau_{D}},y)
	\le g_B(z,y)+ \bE^z \left[ \bE^{X_{\tau_D}} [\mathbf{1}_{\{\sigma_{B_2} < \tau_{2B}\}} g_{2B}(X_{\sigma_{B_2}},y)]\right] \\
	&\le g_{B}(z,y) + \bE^{z} ( 1- h(X_{\tau_{D}}) ) g_{2B}(x_1,y).
	\end{align*}
	Using \eqref{e:xhitA} and \eqref{e:Ato2B}, for $\sE$-q.e.~$z \in B_2 \setminus \{y\}$ we have
	$$  g_{B}(z,y) +g_{2B}(x_1,y)-g_{2B}(z,y)\ge  g_{2B}(x_1,y)  \bE^{z} h(X_{\tau_{D}}) 
	\ge   g_{2B}(x_1,y) p_1^2. $$ 	Letting $z \to x_1$ yields
	$$g_{B}(x_1,y) \ge g_{2B}(x_1,y) p_1^2.$$
	
	Then if $x \in B_1 \setminus \{y\}$,
	\begin{align*}
	g_{2B}(x,y) & =g_{B_2}(x,y) + \bE^x g_{2B}(X_{\tau_{2B}},y) \le g_{B}(x,y) + g_{2B}(x_1,y)\\
	&\le  g_{B}(x,y) +p_1^{-2} g_B(x_1,y). 
	\end{align*} 
 The first equality above only holds for $\sE$-q.e. $x \in B_1 \setminus \{y\}$ but $g_{2B}(x,y)\le  g_{B}(x,y) +p_1^{-2} g_B(x_1,y)$ follows for any $x \in B_1\setminus \{y\}$ by continuity.

	Let $x'_1$ be the point in $\pd B_2$ which minimizes $g_B(x',y)$.
	 	By the maximum principle \eqref{e:max}, 
	$g_B(x,y) \ge g_B(x'_1,y)$.
	We now apply Corollary \ref{C:3.2} to the ball $B(y,a_1 R +a_2R )$ to deduce that
	$g_B(x_1,y) \le c g_B(x'_1,y)$. Combining this with the inequalities above we obtain
	the bound $g_{2B}(x,y) \le C g_{B}(x,y)$. 
	(Note that the constant $C$ only depends on $p_1$ and the constants in 
	Corollary \ref{C:3.2}; it does not depend on $y$.)

	It remains to find $p_1>0$ such that there exist sets
	$D$ and $A$ satisfying  \eqref{e:xhitA} and \eqref{e:Ato2B}.
	Let $y_0 \in \partial B(x_0, (2+\eps)R)$. 
	By Lemma \ref{L:regpath1} there exists a sequence $x_0=z_0, \dots , z_n = y_0$
	 with $d(z_{i-1},z_i)< \epsilon R$ for $1 \le i \le n$ such that if $j = \max\{ i: z_i \in B_3 \}$ then
	$B(z_i, \eps R) \subset B(x_0, K a_3 R)$ for $0\le i \le j$.
	Write $B'_i = B(z_i, \eps R)$. Now let $D = B \setminus \ol B'_j$, and $A = \ol {B'_j}$.

	 We will use  Remark \ref{R:chainprob}(ii) repeatedly to obtain \eqref{e:Ato2B}.  If $i  \ge j$ then $B(z_i, 8 K_3 \eps R ) \cap B_2 = \emptyset$. So we can chain
	along the sequence of balls $B_j, \dots B_n$ to obtain \eqref{e:Ato2B}
	with $p_1 = p_0^n$.
	
	If $0\le i\le j$ then $d(x_0, z_i) \le K a_3 R$ and so
	$B(z_i, 8 \eps K_3R) \subset B(x_0, K a_3 R + 8 \eps K_3 R) \subset B$.
	Hence, chaining along this sequence we obtain
	$$ \bP^x( X_{\tau_D} \in A)  \ge p_0^j \hbox{ for $\sE$-q.e.~} x \in \ol{B'_0}. $$
	To complete the proof of \eqref{e:xhitA} we need to extend this estimate
	to $x \in B_2$. 
	
	Let $x_2 \in B_2$. Then there exists a chain of balls
	$B(w_j, \eps R)$,
 	$0\le j\le k$ with $w_0=x_2$, $w_k = x_0$,  $d(w_{j-1},w_j)<\eps R$ for $1 \le j \le k$,
	and with $B(w_j , \eps R) \subset B(x_0, K a_2 R)$. 
	Since  $B(w_j, 8 \eps K_3 R) \subset B$, we deduce that
	$$ \bP^x ( \sigma_{B'_0} < \tau_B ) \ge p_0^k, \quad \mbox{for $\sE$-q.e.~ $x \in B(x_2,\eps R)$.}$$
 By letting $x_2$ run over a countable dense subset of $B_2$, it follows that
	$$   \bP^x( X_{\tau_D} \in A) \ge p_0^{k+n}\q\mbox{for $\sE$-q.e.~$x \in B_2$}. $$
	Since $k$ and $n$ only depend on the constants 
 	$N_\sX(\eps)$, 
	this completes the proof of \eqref{e:xhitA}. 
	\qed

 	The following corollary is a direct consequence of   Lemmas \ref{L:3.5} and \ref{L:3.12} and Remark \ref{r:nonpolar}.
	
	\begin{corollary} \label{C:3.13}
		(See \cite[Corollary 3.13]{BM1}.)
			Let $(\sX,d,m,\sE,\sF)$ be a MMD space that satisfies Assumption \ref{A:main}.
		There exists $C_1$ that depends only on the constants in Assumption \ref{A:main} such that for all $A > 8K$ and for all $0<r< \diam(\sX,d)/(6A)$, $x \in \sX$,
		\be
		\operatorname{Cap}_{B(x,2Ar)} (B(x,r)) \le \operatorname{Cap}_{B(x,Ar)} (B(x,r)) \le C_1 \operatorname{Cap}_{B(x,2Ar)} (B(x,r)).
		\ee
	\end{corollary}
	
	 \medskip
 
 	In the following,  notations $f \asymp g $,  $f\lesssim g$ 
	and  $f \gtrsim g$ mean  that there are  positive constants $c_1, c_2$
	so that $c_1g \leq f\leq c_2 g$,  $f\leq c_2 g$ and $f\geq c_1 g$, respectively,
	on the common domain of definition of $f$ and $g$.
	
	\ms

	\begin{lemma}  \label{L:3.14} (See \cite[Lemma 3.14]{BM1}.) 
			Let $(\sX,d,m,\sE,\sF)$ be a MMD space that satisfies Assumption \ref{A:main}. 
		\begin{description}
		\item{\rm 	(a)} Let $D$ be an open set in $\sX$ such that $D^c$ is non-$\sE$-polar. 
		Let $x \in \sX$ and $r>0$ be such that $B(x, C_0 r) \subset D$, where
		$C_0=2K+3$.  There exists a constant $C_1>0$ such that
		$$   C_1^{-1} \operatorname{Cap}_{D}(B(y,r)) \le \operatorname{Cap}_{D}(B(x,r)) \le C_1  \operatorname{Cap}_{D}(B(y,r)) \quad \hbox{for } y \in B(x,r). $$
		
		\item{\rm (b)} Let $A >8K$. 
		There exists a constant $C_2 >0$ such that 
		\be \nonumber
		\operatorname{Cap}_{B(x,Ar)}(B(x,r)) \le C_2  \operatorname{Cap}_{B(y,Ar)}(B(y,r))  
		\ee
	 for  $x \in \sX, y \in B(x,A_1 r), 
		\, 0< r < \diam(\sX,d)/(6A)$.
		\item{\rm (c)} Let $A> 8 K$ and $ A_1 >0$. There exists a constant $C_3>0$ such that 
		\be \nonumber
		\operatorname{Cap}_{B(x,Ar)}(B(x,r)) \le C_3  \operatorname{Cap}_{B(y,Ar)}(B(y,r)) 
		\ee
		 for  $x \in \sX, y \in B(x,A_1 r), 
		\, 0< r < \diam(\sX,d)/(6A)$.
		Here the constants $C_1,C_2,C_3$ depend only on $A_1$ and  the constants in Assumption \ref{A:main}.
	\end{description} 
	\end{lemma}

	\proof
	(a) 
	As in the proof of Lemma \ref{L:3.3}, choose $z \in B(x,r)$  such that 
	 $\min(d(z,x), d(z,y)) \ge r/4$.
	By Corollary \ref{C:3.2} and Lemma \ref{L:3.5}, $ \operatorname{Cap}_{D}(B(x,r)) \asymp g_D(x,z)^{-1}$ and $\operatorname{Cap}_D (B(y,r)) \asymp g_D(y,z)^{-1}$.
	The conclusion now follows from Lemma \ref{L:3.3}.\\
	(b) By Corollary \ref{C:3.13} and part(a), we have \[ \operatorname{Cap}_{B(x,Ar)}(B(x,r))  \lesssim  \operatorname{Cap}_{B(x,2Ar)}(B(x,r))\asymp \operatorname{Cap}_{B(x,2Ar)}(B(y,r)).\] Since $B(y,Ar) \subset B(x,2 Ar)$, we have $\operatorname{Cap}_{B(x,2Ar)}(B(y,r)) \le \operatorname{Cap}_{B(y,Ar)}(B(y,r))$. \\
	(c) The case $A_1 \le 1$ follows from (b).  For $A_1  > 1$, by the RBC($K$) condition there exists $N$ such that  $x,y \in \sX$ with $d(x,y) < A_1 r$, can be connected by a sequence of points $x_0=x,x_1,\ldots,x_N=y$ with $d(x_{i-1},x_i)<r$ for $1\le i \le N$, where $N$ depends only on $A_1$ and the constants in RBC($K$) condition.
	By applying (b) repeatedly, we obtain (c) with $C_3= C_2^N$, where $C_2$ is the constant in (b).
	\qed

	\begin{prop} (See \cite[Proposition 3.15]{BM1}) \label{P:3.15}
			Let $(\sX,d,m,\sE,\sF)$ be a MMD space that satisfies Assumption \ref{A:main}.
		Let $D \subset \sX$ be an open set  such that $D^c$ is non-$\sE$-polar and let   
		$B(x_0, 2K R) \subset D$. Let  $b\geq 24$.
		Let $F \subset B(x_0,R)$, and suppose there exist disjoint Borel subsets 
		$\{ Q_i, 1 \le i \le n\}$ of $\sX$ with $n \ge 2$ such that 
		\[
		F = \cup_{i=1}^n Q_i
		\]
		and for each $i$, there exists $z_i \in \sX$ so that 
 		$B(z_i,R/b) \subset Q_i$. 
		Then there exists $\delta= \delta(\delta_H, b ,C_H,K) >0$ such that 
		\[
		\operatorname{Cap}_D(F)\le (1-\delta) \sum_{i=1}^n \operatorname{Cap}_D (Q_i).
		\]
	\end{prop}
	
	\proof
	The proof is similar to that of  \cite[Proposition 3.15]{BM1}.
	The only difference is that we use RBC($K$) condition and a chaining argument using the EHI along with 
	Lemma \ref{L:3.7} to obtain the lower bound on the equilibrium potentials $h_i$ for $\operatorname{Cap}_D (Q_i)$.
	
 The proof in  \cite[Proposition 3.15]{BM1} uses the fact that the $0$-order equilibrium measure $\nu_B^D$ of any $B \subset D$ with $\operatorname{Cap}_D(B)<\infty$ for the part Dirichlet form $(\sE,\sF^D)$ on $D$ satisfies 
	$\nu_B^D(D)=\operatorname{Cap}_D(B)$. 
	Since this is mentioned in \cite{CF,FOT} under the additional assumption that $B$ is compact, we provide further details on how to verify this equality for an arbitrary set $B$. Let $e^D_B$ denote the $0$-order equilibrium potential of $B$ for $(\sE,\sF^D)$, so that $\sE(e^D_B,u)= \int_D u \,d\nu_B^D$ for any $u \in \sF^D_e$ by \cite[Theorem 2.2.5]{FOT}. Then by \cite[the $0$-order version of Theorem 2.1.5]{FOT}, $\operatorname{Cap}_D(B)=\sE(e_B^D,e_B^D)= \int_D e_B^D \,d\nu^D_B, \int_D \varphi(1-e^D_B)\,d \nu^D_B= \sE(e_B^D,\varphi(1-e^D_B))=0$ for any $\varphi \in \sF \cap C_c(\sX)$ with $\restr{\varphi}{\sX \setminus D}=0$, hence $\int_D (1-e^D_B)\, d\nu_B^D=0$, namely $e_B^D=1$ $\nu^D_B$-a.e., and thus $\operatorname{Cap}_D(B)=\int_D e^D_B\, d\nu^D_B=\nu_B^D(D)$.
	\qed

	The following lemma is an extension of Corollary \ref{C:3.13}.

	\begin{lem} \label{l:ratiocap}
			Let $(\sX,d,m,\sE,\sF)$ be a MMD space that satisfies Assumption \ref{A:main}.
		 		 Let $1 < A_1 \le A_2< \infty$. 
 		 There exists a positive constant $C_2$ that depends only on $A_1,A_2$  and the constants in Assumption \ref{A:main} such that 
		for all $x \in \sX, 0<r< \diam(\sX,d)/(6 (A_2 \vee (9K)))$,
	\bes
	\operatorname{Cap}_{B(x,A_2r)} (B(x,r)) \le \operatorname{Cap}_{B(x,A_1 r)} (B(x,r)) 
	\le C_2 \operatorname{Cap}_{B(x,A_2 r)} (B(x,r)).  
	\ees
	\end{lem}

\proof    The estimate $\operatorname{Cap}_{B(x,A_2r)} (B(x,r)) \le \operatorname{Cap}_{B(x,A_1 r)} (B(x,r))$ follows from domain monotonicity. For the other estimate, by domain monotonicity we may assume 
 $A_2>8K$. 

Choose $A_3>8K$ so that $A_2/A_3< A_1 -1$. Then 
$B(y,A_2 r/A_3) \subset B(x,A_1r)$ for all $y \in B(x,r)$.  By the metric doubling property, there exists $N \in \bN$ (depending only on $A_3$ and  the constant associated with metric doubling) such that $y_1,\ldots,y_N \in B(x,r)$ and $\cup_{i=1}^N B(y_i,r/A_3) \supset B(x,r)$. By considering the function $e= \max_{1 \le i \le N} e_i$ where $e_i$ is the equilibrium potential corresponding to $\operatorname{Cap}_{B(y_i, A_2r/A_3)}(B(y_i,r/A_3))$, we obtain
\[
 \operatorname{Cap}_{B(x,A_1 r)} (B(x,r)) \le  \sum_{i=1}^N \operatorname{Cap}_{B(y_i, A_2r/A_3)}(B(y_i,r/A_3)).
\]
As $y_i \in B(x,r)= B(x,A_3(r/A_3))$, by  Lemma \ref{L:3.14}(c), we obtain
\[
 \operatorname{Cap}_{B(y_i, A_2r/A_3)}(B(y_i,r/A_3)) \asymp  \operatorname{Cap}_{B(x, A_2r/A_3)}(B(x,r/A_3)),
\]
for all $x \in \sX$,  $r < \diam(\sX,d)/(6 A_2/A_3)$, and $ i=1,\ldots,N$.
By Corollary \ref{C:3.13} and domain monotonicity, we have   
\[\operatorname{Cap}_{B(x, A_2r/A_3)}(B(x,r/A_3)) \asymp  \operatorname{Cap}_{B(x, A_2r)}(B(x,r/A_3))\le  \operatorname{Cap}_{B(x, A_2r)}(B(x,r)),\]
for all $x \in \sX, r <\diam(\sX,d)/(6A_2)$. We obtain the desired bound
$$
 \operatorname{Cap}_{B(x,A_1 r)} (B(x,r)) \lesssim \operatorname{Cap}_{B(x,A_2 r)} (B(x,r))
$$
 by combining the above three estimates.
\qed

The following lemma is used to compare capacities at different scales.

\begin{lem} \label{l:capcomp}
		Let $(\sX,d,m,\sE,\sF)$ be a MMD space that satisfies Assumption \ref{A:main}.
	Let $A>1$.  
	There exist constants $C_2>1$ and $\gamma>0$ that depend only on $A \wedge (8K)$ and the constants in Assumption \ref{A:main}
	such that for all $x \in \sX$ and 
	$0 < s \le r < \diam(\sX,d)/(6 (A \vee (9K)))$, we have
	\[
	C_2^{-1} \left(\frac{r }{s}\right)^{-\gamma} \operatorname{Cap}_{B(x,A s)}(B(x,s))		 \le	\operatorname{Cap}_{B(x,A r)}(B(x,r)) \le C_2 \left(\frac{r }{s}\right)^\gamma \operatorname{Cap}_{B(x,A s)}(B(x,s)).	
	\]	
 \end{lem}	

\proof 
By Lemma \ref{l:ratiocap}, we may assume without loss of generality that $A>8K$. By Remark \ref{r:nonpolar}, Corollary \ref{C:3.8}, Lemma \ref{L:3.5} and domain monotonicity, we have
\[
\operatorname{Cap}_{B(x,Ar)}(B(x,r)) \asymp g_{B(x,Ar)}(x,r)^{-1} \lesssim \left(\frac{r}{s}\right)^\theta g_{B(x,Ar)}(x,s)^{-1} \lesssim  \left(\frac{r}{s}\right)^\theta\operatorname{Cap}_{B(x,As)}(B(x,s))
\]
for all $x \in \sX, 0< s \le r  < \frac{\diam(\sX,d)}{2A}$, where $\theta>0$ is as given in Corollary \ref{C:3.8}. 

For the reverse inequality, we use Corollary \ref{C:3.13} repeatedly and domain monotonicity to obtain
\[
\operatorname{Cap}_{B(x,As)}(B(x,s)) \lesssim \left(\frac{r}{s}\right)^{\theta_1}\operatorname{Cap}_{B(x,Ar)}(B(x,s)) \le \left(\frac{r}{s}\right)^{\theta_1}\operatorname{Cap}_{B(x,Ar)}(B(x,r)) ,
\]
for all $x \in \sX, 0< s \le r < \frac{\diam(\sX,d)}{6A}$, where $\theta_1= \log_2 C_1>0$, where $C_1$ is as given in Corollary \ref{C:3.13}. Setting $\gamma=\max(\theta,\theta_1)$, we obtain the desired conclusion.
\qed

	\section{ Good doubling measures}\label{S:6} 

As in \cite[Section 4]{BM1} we now use the argument of Volberg and Konyagin \cite{VK}
to construct a new measure $\mu$ such that $(\sX, d,\mu)$ satisfies  {\it volume 
 doubling}; that is, there is a constant $c> 1$ so that $\mu (B(x, 2r) ) \leq c \, \mu (B(x, r))$
for all $x\in \sX$ and $r>0$.
 We need further that  $\mu$ relates well with capacities -- see 
Definition \ref{d:goodmeas} below. 
 One key difference from \cite{BM1}  is that we do not assume bounded geometry condition on 
the original MMD space  $(\sX, d, m, \sE, \sF)$.
Another difference 
 from \cite{BM1} is that we do not have any
cutoff at small length scales. This means that $\mu$ need not be absolutely continuous 
with respect to $m$, and it is not a priori clear that
   $\mu$ is a smooth measure having full  quasi support on $\sX$. 
This property is 
 established in Proposition \ref{p:fullquasisupport} of this section. 
	The key inputs from the previous section are inequalities controlling 
 	capacities 
	Corollary \ref{C:3.13}, Lemma \ref{L:3.14}, Proposition \ref{P:3.15}, and Lemma \ref{l:capcomp}.

\subsection{Construction of a doubling measure}

 In this section, we often make the following assumption.
 \begin{assumption}  \label{A:main2}
 	{\rm 
 	We assume that $(\sX,d)$ is a  metric space such that $\overline{B(x,r)}$ is compact for 
	all $x \in \sX$ and $r>0$, and such that $(\sX,d)$ satisfies one (and hence all) of the equivalent conditions 
	in Theorem \ref{T:ehitomd}.
		Furthermore, we assume that the MMD space $(\sX, d, m, \sE, \sF)$
 		satisfies (scale invariant) EHI with constants $C_H, \delta_H$. 
		 } \end{assumption}	
	
	The following definition is a simplification of \cite[Definition 4.1]{BM1}: 
	we do not require absolute continuity with respect to the reference measure $m$. 
	We do not require the volume doubling property for the measure $\nu$ either -- this will follow from Lemma \ref{l:good-VD}.

	\begin{definition}\label{d:goodmeas} 
{\rm 
 Let $(\sX,d,m,\sE,\sF)$ be a MMD space that satisfies Assumption \ref{A:main2}.
Let $D$ be  a Borel subset of $\sX$.
 Let $C_0\in (1, \infty)$, $0< \beta_1 \le \beta_2$ and let $A \in [2^{1/3},\infty)$. 
	Let $I=(0,A^{-4}\diam(D))$.
 We say a Borel measure $\nu$  on $D$   is {\em $(C_0, A, \beta_1, \beta_2)$-capacity good} if 
 for all $x \in D$ and $s_1, s_2 \in I$ with $s_1<s_2$, $0<\nu(B(x,s_1))\le \nu(B(x,s_2))<\infty$ and
\begin{equation} \label{e:m02}
			C_0^{-1} \left( \frac{s_2}{s_1}\right)^{\beta_1} \le
			\frac{\nu(B(x,s_2)) \operatorname{Cap}_{B(x,As_1)} (B(x,s_1))} {\nu(B(x,s_1)) \operatorname{Cap}_{B(x,As_2)}(B(x,s_2))}
			\le C_0 \left( \frac{s_2}{s_1}\right)^{\beta_2}.
\end{equation}
Since $\nu$ is locally finite, any capacity good measure $\nu$ is a Radon measure on $D$, if $D$ is open in $\sX$.
 }\end{definition}

 	Under Assumption \ref{A:main2},  we observe by Corollary \ref{C:3.13} that 
	the second inequality in \eqref{e:m02} of Definition \ref{d:goodmeas} 
	implies the volume doubling property for $\nu$.
	
	\begin{lem} \label{l:good-VD}
			Let $(\sX,d,m,\sE,\sF)$ be a MMD space that satisfies Assumption \ref{A:main2}.
		Let $\nu$ be a $(C_0,A,\beta_1,\beta_2)$-capacity good measure on $\sX$. Then it satisfies the volume doubling property. 
	\end{lem}
	\proof
	If $\diam(\sX,d)=\infty$, then the volume doubling property follows from Lemma \ref{l:capcomp} and domain monotonicity of capacity, since $\operatorname{Cap}_{B(x,As)}(B(x,s))$ and $\operatorname{Cap}_{B(x,2As)}(B(x,2s))$ are comparable.
	
	In the case $\diam(\sX,d)<\infty$, we view $\sX$ as the closure of the ball 
 	$ \ol{B(x_0,2\diam(\sX,d))} $ 
	and use  Lemma \ref{l:capcomp} to obtain the volume doubling property for balls $B(x,s)$ with  $s \lesssim \diam(\sX,d)$. The volume doubling property for larger balls follows from a covering argument, the metric doubling property and the fact that $\inf_{x \in \sX} \nu(B(x,s)) \gtrsim \nu(B(x_0,s))$ for $s=\frac{1}{3} A^{-4} \diam(\sX,d)$ and any $x_0 \in \sX$ by  RBC$(K)$.
	\qed

	The following  is the main result of this section.
	\begin{theorem}[Construction of a doubling measure] \label{T:meas}
		Let $(\sX,d,m,\sE,\sF)$ be a MMD space that satisfies Assumption \ref{A:main2}.
		Then there exist constants $C_0>1$, $A > 1$,  $0< \beta_1 \le \beta_2$
		and a Borel measure $\mu$ on $\sX$ which is $(C_0, A, \beta_1, \beta_2)$-capacity good.
	\end{theorem}
	
 	The proof of Theorem \ref{T:meas} requires a preparation of a few results. 
	 	We begin by adapting the argument in \cite{VK} to construct
	 	a measure with the desired properties   on a family of compact sets. 
	We then follow \cite{LuS} and 
	obtain $\mu$ as a weak$^*$ limit of measures defined on an
	increasing family of compact sets.

	The proof uses a family of generalized dyadic cubes, which 
	provide a family of nested partitions of a space. Such a decomposition of space was introduced by 
	Christ \cite[Theorem 11]{Chr}. The following is a slight modification of the construction in 
 	\cite[Theorem 2.1]{KRS}.
	 Since the requirements (g) and (h) are new, we provide some details.
	
	\begin{lemma}\label{L:gdc} 
 		Let $(\sX,d)$ be a complete metric space satisfying metric doubling and RBC$(K)$ property.
		Let $x_0 \in \sX$ and  $A \ge 8$.
	Then there exists a collection $\Sett{Q_{k,i}}{k \in \bZ, i \in I_k \subset \bZ_+}$  of Borel sets satisfying the following properties:
		\begin{enumerate}[\rm (a)]
	\item  $\sX = \cup_{i\in I_k} Q_{k,i}$ for all $k \in \bZ$,  and $Q_{k,i} \cap Q_{k,j}=\emptyset$ for all $k \in \bZ$ and $i,j \in I_k$ with $i\neq j$.
	
			\item If $m \le n$ and $i \in I_n$, $j \in I_m$,  then either 
			$Q_{n,i} \cap Q_{m,j} = \emptyset$ or $Q_{n,i} \subset Q_{m,j}$. 
			
			\item For every $k  \in \bZ$ and $i \in I_k$, there exists $x_{k,i}\in  Q_{k,i}$ such that
			\[
			B(x_{k,i},c_A A^{-k}) \subset Q_{k,i} \subset \overline {B(x_{k,i},  C_A A^{-k}) } ,
			\]
			where $c_A=  \frac{1}{2} - \frac{1}{A-1}$ and $C_A= \frac{A}{A- 1}$. 
			
  \item The sets $N_k= \Sett{x_{k,i}}{i \in I_k}$, where $x_{k,i}$ are as defined in  {\rm (c)},  are increasing in $k$ and $x_0 \in N_k$ for all $k \in \bZ$; that is $N_k \subset N_{k+1}$ for all $k \in \bZ$   and $x_0\in \cap_{k\in \bZ} N_k$.  
  
			\item  Properties  {\rm (a), \rm (b) and \rm{(c)}} define a partial order $\prec$ on $\mathcal{I}=\Sett{(k,i)}{k \in \bZ, i \in I_k}$ by inclusion, where
			$(k,i) \prec (m,j)$ if and only if $k \ge m$ and  $Q_{k,i} \subset Q_{m,j}$.
			
			\item There exists $C_M=C_M(A)>0$ such that, for all $k \in \bZ$  
 		and for all $x_{k,i} \in N_k$, the  `successors'
			\[
			S_{k}(x_{k,i})= \Sett{x_{k+1,j}}{(k+1,j) \prec (k,i)}
			\]
			satisfy
			\be \label{e:up}
	 			  \abs{S_k(x_{k,i})}  \leq  C_M 
			\q \mbox{for all $k \in \bZ, i \in I_k$.}
			\ee
			Furthermore, we have $d(x_{k,i},y) < A^{-k}$ for all $y \in S_k(x_{k,i})$.
\item 			
			Let  
			\be \label{e:defk0} 
			k_0= \inf \set {k \in \bZ: \abs{I_k} > 1},
			\ee where $\abs{I_k}$  denotes the cardinality of $I_k$. 
			Then $k_0 \in \bZ\cup \set{-\infty}$ satisfies
			\be \label{e:estk0}
			c_A A^{-k_0} \le \diam(\sX,d)\le 2 C_A A^{1-k_0}.
			\ee
			For all $k \ge k_0$,  $k \in \bZ$ and $i \in I_k$, we have $\abs{S_k(x_{k,i})} \ge 2$. 
			
		\item  For all $k \in \bZ$, $0 \in I_k, Q_{k,0}$ is compact and $x_{k,0}=x_0$.
		\end{enumerate}
	\end{lemma}
	\proof
	The sets $Q_{k,j}, k \in \bZ, j \in I_k$ are referred to as `generalized dyadic cubes'.
	We follow the construction in \cite{KRS} with a minor modification so as to ensure the property (h). 
	
	We choose $N_0 \subset \sX$ such that $x_0 \in N_0$ and $N_0= \set{x_{0,i}: i \in I_0}$ is a maximal subset of $\sX$ such that $d(x_{0,i},x_{0,j}) \ge 1$ for all $i \neq j$ with $i,j \in I_0$. For $k > 0$, we define  $N_k= \set{x_{k,i}: i \in I_k}$
as a maximal subset of $\sX$ such that $N_{k-1} \subset N_k$ and $d(x_{k,i},x_{k,j}) \ge A^{-k}$ for 
all distinct $x_{k,i},x_{k,j} \in N_k$. For $k < 0$, we define $N_k= \set{x_{k,i}: i \in I_k}$
as a maximal set such that $x_0 \in N_{k} \subset N_{k+1}$, and $d(x_{k,i},x_{k,j}) \ge A^{-k}$ for 
all distinct $x_{k,i},x_{k,j} \in N_k$. 

 We label the indices $I_k$ such that $0 \in I_k$ and $x_{k,0}=x_0$ for all $k \in \bZ$.
For each $(k,i) \in \bZ \times \mathbb{Z}_+$ with $i \in I_k$, we pick an element $(k-1,j) \in \bZ \times \mathbb{Z}_+$  with $j \in I_{k-1}$ such that \[d(x_{k,i},x_{k-1,j})= \min_{l \in I_{k-1}} d(x_{k,i},x_{k-1,l}).\]
We define $\prec$ as the smallest partial order that contains the relations $(k,i) \prec (k-1,j)$ for all $(k,i) \in \bZ \times  \mathbb{Z}_+$ with $i \in I_k$, where $(k-1,j) \in \bZ \times  \mathbb{Z}_+$ with $j \in I_{k-1}$ is chosen as above.

We relabel the indices $I_0$ of $N_0$  such that $0 \in I_0$ remains unchanged  and  
\be \label{e:extra}
l_1 < l_2 \q \mbox{for all $k < 0, l_1 \in \set{i \in I_0: (0,i) \prec (k,0)}$ and $l_2 \in I_0 \setminus \set{i \in I_0: (0,i) \prec (k,0)}$.} 
\ee
This relabeling exists since $\set{i \in I_0: (0,i) \prec (k,0)}$) is finite for all $k<0$ (by the doubling property) and $(k,0)\prec(k-1,0)$ for all $k \le 0$.
 
 Define the sets $Q_{0,i}$ as
 \[
 Q_{0,i}= \overline{\set{x_{l,k}: (l,k)\prec (0,i)}} \setminus \bigcup_{j<i, j \in I_0} Q_{0,j}.
 \]
	For $k<0$, we define the sets $Q_{k,i}$ inductively as
	\[
	Q_{k,i}= \bigcup_{(k+1,j)\prec (k,i)} Q_{k+1,j},
	\]
	whereas for $k>0$, we define
	\[
	Q_{k,i} = Q_{k-1,i' }\cap \overline{\set{x_{l,j}:(l,j) \prec (k,i)}} \setminus \bigcup_{j<i, j \in I_k} Q_{k,j}, \q \mbox{where $(k,i)\prec (k-1,i')$.}
	\]
	Properties (a)-(e) are contained in  \cite[Theorem 2.1]{KRS};  (f) is immediate from the above construction and metric doubling.
	
\noindent	(g) The estimate $\abs{S_k(x_{k,i})}  \ge 2$ relies on the following consequence of RBC($K$) (see Lemma \ref{l:metric}(a)): $r \le \diam(B(x,r)) \le 2r$ for all $B(x,r) \neq \sX$. 	Since
	$2 C_A/ c_A=  4A/ (A-3) < A$ for all $A \ge 8$,  we have
	\[
	 \diam(Q_{k,i}) \ge c_A A^{-k} > 2 C_A A^{-k-1} \ge \diam(Q_{k+1,j}) \q \mbox{for all $k \ge k_0, k \in \bZ$.}
	\]
	Hence $Q_{k,i} \neq Q_{k+1,j}$ for all $k \ge k_0, i \in I_k, j \in I_{k+1}$, and therefore $\abs{S_k(x_{k,i})} \ge 2$ for all $k \ge k_0$. 
	
Clearly by (c), $\diam(\sX,d)=\infty$ if and only if $k_0=-\infty$.
 If $k_0 \in \bZ$, the estimate \eqref{e:estk0} follows from $ B(x_0,c_A A^{-k_0}) \subset Q_{k_0,0} \subsetneq \sX = Q_{k_0-1,0} \subset  \overline {B(x_0,C_A A^{-k_0+1})}$.
	
\noindent (h) 	By \eqref{e:extra}, $Q_{k,0}$ is closed for all $k \in \bZ$, since $Q_{k,0}= \overline{ \{x_{l,j}: (l,j) \prec (k,0) \} }$. By (c) and \BC, $Q_{k,0}$ is compact for all $k \in \bZ$.
	\qed
	
	\sms

	We fix a family  
	\[ \Sett{Q_{k,i}}{k \in \bZ, i \in I_k \subset \bZ_+},\] 
	of generalized dyadic cubes as given by Lemma \ref{L:gdc},
	and define the nets $N_k$ and successors $S_k(x)$ as in
	the lemma. 
	
		\begin{definition} \label{d:pred}
		{\rm
		We define
			the {\it predecessor } $P_k(x)$ of $x\in N_k$ to be the unique 
			element of $N_{k-1}$ such that $x \in S_{k-1}(P_k(x))$. 
			Note that for $x \in N_k$, $S_k(x) \subset N_{k+1}$ whereas $P_k(x) \in N_{k-1}$. 
			For $x \in \sX$, we denote by $Q_k(x)$ the unique  $Q_{k,i}$ such that $x \in Q_{k,i}$.
	}\end{definition}

 Let $k_0 \in \bZ \cup \set{-\infty}$ be as defined in \eqref{e:defk0}.
 For $ k \in \bZ$ with $k \ge k_0+2$,   denote by $c_k (x) $ the relative capacity
	\be \label{e:defck}
	c_k(x)= \operatorname{Cap}_{B(x, A^{-k+1})} (Q_{k}(x)).
	\ee
	The following lemma provides useful estimates on $c_k$.
	Note that if $k \ge k_0+2$, then 
	$$
	A^{-k+1} \le A^{-k_0-1} \le  c_A^{-1} A^{-1} \diam(\sX,d)= \frac{2(A-1)}{(A-3)A}\diam(\sX,d).
	$$
	
	\begin{lemma}  [Relative capacity  estimates for generalized dyadic cubes] \label{l-ced}
			Let $(\sX,d,m,\sE,\sF)$ be a MMD space that satisfies Assumption \ref{A:main2}.
		There exists  $A_0 \ge 8$ such that the following hold.
		\begin{description}
	\item{\rm (a)}  For all $A \ge A_0$, there exists $C_1>0$ such that for all $k \ge k_0+2, x,y \in \sX$ with 
		 $d(x,y) \le 4 A^{-k}$, we have
		\begin{equation}\label{e:ce1}
		C_1^{-1} c_k(y) \le c_k(x) \le C_1 c_k(y).
		\end{equation}
		
		\item{\rm (b)} For all $A \ge A_0$, there exists $C_1=C_1(A)>0$ such that  for all $k \ge k_0+ 2, x \in N_k$ and $ y \in S_k(x)$, we have
		\begin{equation}\label{e:ce2}
		C_1^{-1} c_k(x) \le c_{k+1}(y) \le C_1 c_k(x).
		\end{equation}
		
		\item{\rm (c)}  For all $A \ge A_0$, there exists $C_1=C_1(A)>0$  such that  for all 
		$x \in \sX$ and $s < \diam(\sX,d)/A^4$,
		\begin{equation} \label{e:ce3}
		C_1^{-1}c_k(x) \le	 \operatorname{Cap}_{B(x,As)}(B(x,s))  \le C_1 c_k(x)
		\end{equation}
		where $k \in \bZ$ is the unique integer such that $A^{-k} \le s < A^{-k+1}$.
	\end{description} 
	\end{lemma}
	
	\proof
	We use domain monotonicity of capacity along with Corollary \ref{C:3.13}, Lemma \ref{L:3.14}(c)  and Lemma \ref{l:capcomp} to show first \eqref{e:ce3} with $s=A^{-k}$ for all $k \ge k_0 +2$ and $x \in \sX$ (with $C_1$ independent of $A$) and then use Lemma \ref{L:3.14}(c) and Lemma \ref{l:capcomp} to obtain the above estimates.
	For (c), note that $A^{-k} \le s < \diam(\sX,d)/A^4 \le 2 C_A A^{-k_0-3}< A^{-k_0-2}$ implies $k \ge k_0+3$.
	\qed
	
	We record one more estimate regarding the subadditivity of $c_k$,
	which will play an essential role in  ensuring \eqref{e:m02} and follows from Proposition \ref{P:3.15} and domain monotonicity of capacity.
	\begin{lemma} (\cite[Lemma 4.6]{BM1}) \label{L:ESP}
		Let $(\sX,d,m,\sE,\sF)$ be a MMD space that satisfies Assumption \ref{A:main2}.
		There exists $A_0 \ge 8$ such that the following holds: for all $A \ge A_0$ there exists $\delta =\delta(A) \in (0,1)$ such that for all $k \in \bZ$, $k \ge k_0 + 2$ and  for all $x \in N_k$, we have
		\[
		c_k(x) \le (1-\delta) \sum_{y \in S_k(x)} c_{k+1}(y).
 		\]
	\end{lemma}
	Henceforth, we fix an $A \ge 8$ large enough  such that the conclusions of Lemmas \ref{l-ced} and \ref{L:ESP} hold.
	
	\medskip
	
 	We need the following  modification of    	\cite[Lemma, p. 631]{VK}, which was stated in \cite[Lemma 4.7]{BM1} without a proof. For the reader's convenience,
	we provide its full proof below.

	\begin{lemma}\label{L:ind}
	 		Let $(\sX,d,m,\sE,\sF)$ be a MMD space that satisfies Assumption \ref{A:main2}.
		Let $c_k(\cdot)$, $k \ge k_0+2$, $k \in \bZ$ denote the capacities of the corresponding generalized dyadic cubes as defined in \eqref{e:defck}.
		There exists $C > 1$ satisfying the following.
	  Let $k \ge k_0+2$, $k \in \bZ$, and let $\mu_k$ be a probability measure on $N_k$ such that
		\be \label{e:muk}
		\frac{\mu_k(e')}{c_k(e')} \le  C^2 \frac{\mu_k(e'')}{c_k(e'')}
		\q \hbox{for all $e',e'' \in N_k$ with $d(e',e'') \le 4 A^{-k}$}.
		\ee
		Then there exists a 
		probability measure $\mu_{k+1}$ on $N_{k+1}$   such that  the following hold.
		\begin{enumerate}[\rm (1)]
			\item 
			For all $g',g'' \in N_{k+1}$ with $d(g',g'') \le 4 A^{-k-1}$ we have
			\be  \label{e:mure1}
			\frac{\mu_{k+1}(g')}{c_{k+1}(g')} \le  C^2 \frac{\mu_{k+1}(g'')}{c_{k+1}(g'')}.
			\ee
			\item 
			Let $\delta \in (0,1)$ be the constant in Lemma \ref{L:ESP}.
			For all points $e \in N_k$ and $g \in S_k(e)$,
			\be \label{e:mure2}
			C^{-1}\frac{\mu_{k}(e)}{c_{k}(e)} \le    \frac{\mu_{k+1}(g)}{c_{k+1}(g)} 
			\le (1-\delta) \frac{\mu_{k}(e)}{c_{k}(e)}.
			\ee
			\item The construction of the measure $\mu_{k+1}$ from the measure 
			$\mu_k$ can be regarded as the transfer of masses from the points of
			$N_k$ to those of $N_{k+1}$, with no mass transferred over a distance greater 
			than $(1+4/A)A^{-k}$.
		\end{enumerate}
	\end{lemma}
	
	\proof  	By the metric doubling property
	\begin{equation}\label{e:in1}
	\sup_{k \in \bZ} \sup_{x \in N_k} \abs{S_k(x)} = S < \infty,
	\end{equation}
	where $S_k(x)$ is as  defined in Lemma \ref{L:gdc}(f).
	We choose
	\[
	C=C_1 S,
	\]
	where $C_1$ is chosen such that \eqref{e:ce1}, \eqref{e:ce2} and \eqref{e:ce3} hold.
Let $k \ge k_0+2$, $k \in \bZ$, and let  $\mu_k$ be any  probability measure on $N_k$ such that \eqref{e:muk} holds.

	The transfer of mass is accomplished in two steps.
	In the first step we distribute the mass $\mu_k(e)$ to all its successors $S_k(e)$ such that the mass of $g \in S_k(e)$ is proportional to $c_{k+1}(g)$; that is
	\[
	f_0(g)= \frac{c_{k+1}(g)}{\sum_{g' \in S_k(e)} c_{k+1}(g')} \mu_k(e),
	\]
	for all $e \in N_k$ and $g \in S_k(e)$.

	By \eqref{e:in1}, Lemma \ref{l-ced} and Lemma \ref{L:ESP},
	we have
	\begin{equation} \label{e:in2}
	C^{-1}\frac{\mu_{k}(e)}{c_{k}(e)} \le    \frac{f_0(g)}{c_{k+1}(g)} \le (1-\delta) \frac{\mu_{k}(e)}{c_{k}(e)},
	\end{equation}
	for all points $e \in N_k$ and $g \in S_k(e)$.
	If the measure $f_0$ on $N_{k+1}$ satisfies condition (1) of the Lemma, we set $\mu_{k+1}=f_0$. This is the desired measure. Condition (2) is satisfied by \eqref{e:in2}, and (3) is obviously satisfied by Lemma \ref{L:gdc}(f).
	The second step is not necessary in this case. 
	
	But if $f_0$ does not satisfy condition (1) of the Lemma, then we proceed as follows at the second step. Let $p_1,\ldots,p_T$ be the indexed pairs of points
	$\set{g',g''}$ with $g',g'' \in N_{k+1}$ and $0 < d(g',g'') \le 4 A^{-k-1}$. Take the pair $p_1= \set{g_1',g_1''}$. If
	$\frac{f_0(g_1')}{ c_{k+1}(g_1')}\le C^2 \frac{f_0(g_1'')}{ c_{k+1}(g_1'')}$
	and
	$\frac{f_0(g_1'')}{ c_{k+1}(g_1'')} \le C^2 \frac{f_0(g_1')}{ c_{k+1}(g_1')}$,
	then we set $f_1=f_0$. Assume one of the inequalities is violated, say
	$
	\frac{f_0(g_1')}{ c_{k+1}(g_1')} > C^2 \frac{f_0(g_1'')}{ c_{k+1}(g_1'')}.
	$ Then we construct a measure $f_1$ from $f_0$ such that
	\begin{align*}
	f_1(g_1') &=f_0(g_1') - \alpha_1,\\
	f_1(g_1'') &= f_0(g_1'') + \alpha_1,\\
	f_1(g) &= f_0(g), \hspace{4mm} g \neq g_1',g_1'',
	\end{align*}
	where $\alpha_1 >0$ is chosen such that
	\[
	\alpha_1 \left( \frac{C^2}{c_{k+1}(g_1'')}+ \frac{1}{c_{k+1}(g_1')} \right) = \frac{f_0(g_1')}{ c_{k+1}(g_1')}- C^2 \frac{f_0(g_1'')}{ c_{k+1}(g_1'')}.
	\]
	It is clear that $\frac{f_1(g_1')}{ c_{k+1}(g_1')} = C^2 \frac{f_1(g_1'')}{ c_{k+1}(g_1'')}$.
	
	The next step is the construction of a measure $f_2$ from $f_1$ in exactly the same way that $f_1$ was constructed from $f_0$. Here we consider the pair $p_2$. A measure $f_3$ is next constructed from $f_2$ and so on. We claim that $\mu_{k+1}=f_T$ is the desired measure in the lemma.  If $\sX$ is non-compact, $\mu_{k+1}(g):=\lim_{j \to \infty}  f_j(g), g \in N_{k+1}$ (the existence of this limit is an easy consequence of the metric doubling property).
	
	We first verify that for all $e \in N_k$, for all $g \in S_k(e)$ and for all $s=0,1,\ldots,T$, we have
	\begin{equation}\label{e:in3}
	C^{-1}\frac{\mu_{k}(e)}{c_{k}(e)} \le    \frac{f_s(g)}{c_{k+1}(g)} \le (1-\delta) \frac{\mu_{k}(e)}{c_{k}(e)}.
	\end{equation}
	By \eqref{e:in2}, it is clear that \eqref{e:in3} holds for $s=0$. We now show \eqref{e:in3} by induction. Suppose \eqref{e:in3} holds for $s=j$, we will verify it for $s=j+1$.
	Let $p_{j+1}= \{g',g''\}$, $e'=P_{k+1}(g'),e''=P_{k+1}(g'')$. If $f_j=f_{j+1}$, there is nothing to prove. But if $f_{j+1} \neq f_j$, then assume, say, that
	\begin{equation} \label{e:in4}
	\frac{f_j(g')}{ c_{k+1}(g')} > C^2 \frac{f_j(g'')}{ c_{k+1}(g'')}.
	\end{equation}
	By \eqref{e:in4} and the construction, we have
	\begin{equation}\label{e:in4p}
	f_{j+1}(g') < f_j(g'), \hspace{4mm} f_{j+1}(g'') > f_j(g'').
	\end{equation}
	Therefore by the induction hypothesis \eqref{e:in3} for $s=j$ and \eqref{e:in4p}, we have
	\begin{equation*}
	\frac{f_{j+1}(g')}{c_{k+1}(g')} \le (1-\delta) \frac{\mu_k(e')}{c_k(e')},  \hspace{4mm}  \frac{f_{j+1}(g'')}{c_{k+1}(g'')} \ge C^{-1} \frac{\mu_k(e')}{c_k(e')}.
	\end{equation*}
	Therefore it suffices to verify that
	\begin{equation}\label{e:inAB}
	\frac{f_{j+1}(g')}{c_{k+1}(g')} \ge C^{-1} \frac{\mu_k(e')}{c_k(e')},  \hspace{4mm}  \frac{f_{j+1}(g'')}{c_{k+1}(g'')} \le  (1-\delta) \frac{\mu_k(e'')}{c_k(e'')}.
	\end{equation}
	Suppose the first inequality in \eqref{e:inAB} fails to be true, then by construction, \eqref{e:in4p} and  the induction hypothesis \eqref{e:in3} for $s=j$, we have
	\begin{equation}
	C^{-1} \frac{\mu_k(e')}{c_k(e')} > \frac{f_{j+1}(g')}{c_{k+1}(g')} = C^2 \frac{f_{j+1}(g'')}{c_{k+1}(g'')} > C^2  \frac{f_{j}(g'')}{c_{k+1}(g'')} \ge C \frac{\mu_k(e'')}{c_k(e'')},
	\end{equation}
	which implies $\frac{\mu_k(e')}{c_k(e')}  > C^2 \frac{\mu_k(e'')}{c_k(e'')}$. However $\frac{\mu_k(e')}{c_k(e')}  \le C^2 \frac{\mu_k(e'')}{c_k(e'')}$, by
	the assumption on $\mu_k$, since
	\[
	d(e',e'') \le d(e',g')+d(g',g'')+d(e'',g'') \le 2 A^{-k}+ 4 A^{-k-1} \le 4 A^{-k} .
	\]
	This proves the first inequality in \eqref{e:inAB}. The proof of the second inequality in \eqref{e:inAB} is similar. Indeed, assume to the contrary that  $\frac{f_{j+1}(g'')}{c_{k+1}(g'')} >  (1-\delta) \frac{\mu_k(e'')}{c_k(e'')}$; then we have
	\begin{equation}
	(1-\delta) \frac{\mu_k(e')}{c_k(e')} \ge\frac{f_{j}(g')}{c_{k+1}(g')} > \frac{f_{j+1}(g')}{c_{k+1}(g')} = C^2 \frac{f_{j+1}(g'')}{c_{k+1}(g'')} > C^2 (1-\delta)\frac{\mu_k(e'')}{c_k(e'')},
	\end{equation}
	which again implies $\frac{\mu_k(e')}{c_k(e')}  > C^2 \frac{\mu_k(e'')}{c_k(e'')}$.
	Therefore \eqref{e:in3} follows by induction. In particular, $\mu_{k+1}=f_T$ satisfies condition (2) of the lemma.
	
	We now verify condition (1) for $\mu_{k+1}=f_T$. For this, it suffices to prove the following assertion: if
	\begin{equation} \label{e:in6}
	C^{-2} \frac{f_j(g'')}{c_{k+1}(g'')}\le \frac{f_j(g')}{c_{k+1}(g')} \le C^2 \frac{f_j(g'')}{c_{k+1}(g'')}
	\end{equation}
	holds for a pair of points $g',g'' \in N_{k+1}$ such that $0 < d(g',g'') \le 4 A^{-k-1}$, then the same inequalities hold when $f_j$ is replaced by $f_{j+1}$.
	
	We now prove this.
	If $p_{j+1}= \set{g',g''}$, then $f_{j+1}=f_j$ and there is nothing to prove. If $\set{g',g''} \cap p_{j+1} = \emptyset$, then again there is nothing to prove.
	Let $p_{j+1} = \set{g_1,g_2}$.
	Without loss of generality, we assume $p_{j+1} \cap \set{g',g''}= \set{g_1}$ where $g_1=g''$ and $f_{j}(g'')/c_{k+1}(g'')> C^2 f_{j}(g_2)/c_{k+1}(g_2)$. Then
	\begin{equation}\label{e:in7}
	\frac{f_{j+1}(g'')}{c_{k+1}(g'')} = C^2  \frac{f_{j+1}(g_2)}{c_{k+1}(g_2)},\hspace{4mm} f_{j+1}(g'') < f_{j}(g''), \hspace{4mm} f_{j+1}(g')=f_j(g').
	\end{equation}
	Therefore, only the second inequality in \eqref{e:in6} can fail for $f_{j+1}$. Suppose that this happens, that is
	\begin{equation}\label{e:in8}
	\frac{f_{j+1}(g')}{c_{k+1}(g')} > C^2 \frac{f_{j+1}(g'')}{c_{k+1}(g'')}.
	\end{equation}
	Let $e'= P_{k+1}(g')$ and $e_2= P_{k+1}(g_2)$.
	Then by \eqref{e:in8}, \eqref{e:in7} and \eqref{e:in3}
	\begin{equation}
	(1-\delta)\frac{\mu_k(e')}{c_k(e')} \ge \frac{f_{j+1}(g')}{c_{k+1}(g')} > C^2 \frac{f_{j+1}(g'')}{c_{k+1}(g'')} = C^4 \frac{f_{j+1}(g_2)}{c_{k+1}(g_2)} \ge C^3 \frac{\mu_k(e_2)}{c_k(e_2)},
	\end{equation}
	which implies that $\frac{\mu_k(e')}{c_k(e')} > C^2 \frac{\mu_k(e_2)}{c_k(e_2)}$. However since $d(e',e_2)\le d(e',g') + d(g',g'')+d(g_1,g_2)+ d(g_2,e_2) \le 2(A^{-k}+4A^{-k-1}) \le 4 A^{-k}$, we have a contradiction and hence \eqref{e:in8} is false. This shows \eqref{e:in6} for the case  $f_{j}(g'')/c_{k+1}(g'')> C^2 f_{j}(g_2)/c_{k+1}(g_2)$. The case  $f_{j}(g'')/c_{k+1}(g'')<  C^{-2} f_{j}(g_2)/c_{k+1}(g_2)$ is analyzed similarly and therefore the assertion given by \eqref{e:in6} is proved.
	It remains to observe that this assertion proves condition  (1) of the lemma for the measure $\mu_{k+1}=f_T$. Along the path from $f_0$ to $f_T$, we ``correct'' the measure at all pairs of points where condition (1) is violated, and the assertion given by \eqref{e:in6} shows that once a pair is corrected, it remains corrected when further changes are made.
	
  It remains to verify condition (3). Note that by Lemma \ref{L:gdc}(f), there was a mass transfer over a distance of at most $A^{-k}$ while passing from $\mu_k$ to $f_{0}$. Therefore it suffices to verify that while passing from $f_0$ to $f_T=\mu_{k+1}$ there  is a transfer over a distance of at most $4A^{-k-1}$.
	
	We will now verify this.   It suffices to verify that there are no pairs 
	$p_l=\set{g_1,g_2}$, $p_{n}=\set{g_2,g_3}$, $l,n \in \bZ \cap [1,T], l \neq n$, such that some mass is transferred from $g_1$ to $g_2$ (in the transition from $f_{l-1}$ to $f_{l}$) and some mass is transferred from $g_2$ to $g_3$ (in the transition from $f_{n-1}$ to $f_{n}$). Assume the opposite.
	First note that the assertion given by \eqref{e:in6} can be modified as follows. If the second inequality in \eqref{e:in6} is true for $f_j$ it remains true for $f_{j+1}$. The same argument as before goes through. Using this modified version of the assertion, and the assumption that there are mass transfers from
	$g_1$ to $g_2$ and  
	from $g_2$ to $g_3$, we have
	\begin{equation} \label{e:in9}
	\frac{f_0(g_1)}{c_{k+1}(g_1)} > C^2 \frac{f_0(g_2)}{c_{k+1}(g_2)}, \hspace{4mm} \frac{f_0(g_2)}{c_{k+1}(g_2)} > C^2 \frac{f_0(g_3)}{c_{k+1}(g_3)}.
	\end{equation}
	If $e_1= P_{k+1}(g_1)$, $e_3=P_{k+1}(g_3)$, then
	\[
	d(e_1,e_3) \le d(e_1,g_1)+d(g_1,g_2)+d(g_2,g_3)+ d(g_3,e_3) \le 2(A^{-k}+4A^{-k-1}) \le 4 A^{-k}.
	\]
	Consequently by assumption, $\mu_k(e_1)/c_k(e_1) \le C^2 \mu_k(e_3)/c_k(e_3)$. 
	However the inequalities \eqref{e:in9} and \eqref{e:in3} imply the opposite inequality $\mu_k(e_1)/c_k(e_1) > C^2 \mu_k(e_3)/c_k(e_3)$. 
	We have arrived at the desired contradiction and the proof of the lemma is complete.
	\qed
 \begin{remark}  \label{r:subset}
 	{\rm 
 		 Lemma \ref{L:ind} and its proof above remain valid, if $N_k$ is replaced by any $M_k \subset N_k$ and $N_{k+1}$ by ${M_{k+1}}= \bigcup_{y \in M_k} S_k(y)$.
 }
 \end{remark}
	
	We now adapt the method in \cite{VK} to construct the doubling measure.

	\begin{proposition}[Measure in a cube] \label{P:ballmeas}
Let $(\sX,d,m,\sE,\sF)$ be a MMD space that satisfies Assumption \ref{A:main2}.   
		There exist  constants $C_0>1$,  $A\ge 8$ and  $0< \beta_1 \le \beta_2$ such that for any 
 		integer  $l \ge k_0-1$, 
 there exists a  $(C_0,A, \beta_1, \beta_2)$-capacity good measure $\nu=\nu_{l}$ on $Q_{l,0}$.
	\end{proposition}
	
\proof
	Choose $A \ge 8$ large enough such that the conclusions of Lemmas   \ref{l-ced} and \ref{L:ESP} hold.
 Set $M_k=N_k \cap Q_{l,0}$ for $k \ge l+3$, so that $M_{k+1}=\bigcup_{y \in M_k} S_k(y)$ (cf. Remark \ref{r:subset}).
 	Let $\mu_{l+3}$ be the probability measure on $M_{l+3}$ such that $\mu_{l+3}$ is proportional to $c_{l+3}$; that is
	\[
	\mu_{l+3}(x)= \frac{c_{l+3}(x)}{\sum_{y \in M_{l+3}} c_{l+3}(y)}, \q \mbox{for all $x \in M_{l+3}$.}
	\]
	We use Lemma \ref{L:ind} and Remark \ref{r:subset}
	to inductively construct probability measures $\mu_k$ on $M_k$ for all $k \ge l+3$. 
	We define the measure $\nu=\nu_{l}$ as a weak (sub-sequential) limit of the measures $\mu_k$ as $k \to \infty$ (the existence of such a limit follows from the compactness of $Q_{l,0}$).  We claim that 
	\be \label{e:measb1}\mbox{
		$\nu$ is $(C_0,A,\beta_1,\beta_2)$-capacity good for some $C_0,\beta_1,\beta_2>0$.}
	\ee
	For each $x \in Q_{l,0}$ and $k \ge l+3$, let $e_{x,k} \in M_k$ be the unique point in $M_k$ such that $e_{x,k} \in Q_{k}(x)$, so that by Lemma \ref{L:gdc}(c),
	\[
 d(x,e_{x,k}) \le C_A A^{-k}.
	\]
	If $s < A^{-4} \diam(Q_{l,0}) \le A^{-4} 2 C_A A^{-l} \le A^{-l-3}$, then $s<A^{-l-3}$.
	In order to show \eqref{e:measb1}, we prove the following two-sided estimate on measure of balls: there exists $C_2 \ge 1$ such that
	\be \label{e:measb2}
 C_2^{-1} \mu_n(e_{x,n}) \le \nu(B(x,s)\cap Q_{l,0}) \le C_2 \mu_n(e_{x,n}), 
\q \mbox{for all $x\in Q_{l,0},s< A^{-l-3}$,}
	\ee
 	where $n$ is the unique integer such that $A^{-n-1} \le s < A^{-n}$.
	
	Note that, by Lemma \ref{L:ind}(3) the mass from $e \in M_k$ travels a distance of at most 
	\be \label{e:measb3}
	(1+4 A^{-1}) \sum_{l=k}^\infty A^{-l}= C_3A^{-k}, \q \mbox{where $C_3:= (1+4A^{-1})(1-A^{-1})^{-1}\le \frac{12}{7}$.}
	\ee
	Therefore, none of the mass outside $M_n \cap B(x,(1+C_3)A^{-n})$ falls in $B(x,s)$, 
 	and therefore
	\be  \label{e:measb4}
	\nu(B(x,s)) \le \mu_n \left(M_n \cap B(x,(1+C_3)A^{-n}) \right)
	 \q \mbox{for all $x \in Q_{l,0}, s \in (0,A^{-l-3})$. }
	\ee
 	By the triangle inequality, if $e\in M_n \cap B(x,(1+C_3)A^{-n})$, then 
	$$
	d(e,e_{x,n}) \le d(e,x) + d(x,e_{x,n}) \le  (1+C_3)A^{-n} +C_A A^{-n} \le 4 A^{-n}.
	$$
	  Therefore by \eqref{e:measb4}, \eqref{e:muk}, \eqref{e:ce1}, and the metric doubling property, we obtain the upper bound  $\nu(B(x,s)) \lesssim \mu_n(e_{x,n})$ in \eqref{e:measb2}.
	
	For the lower bound in \eqref{e:measb2}, using \eqref{e:measb3}, we have that for all $x \in Q_{l,0}$ and for all $s < A^{-l-3}$ with $A^{-n-1} \le s <A^{-n}, n \in \bZ$, the mass from $e_{x,n+2}$ travels a distance of at most
	$C_3 A^{-n-2} \le \frac{12}{7} A^{-n-2}$ from $e_{x,n+2}$. Since $d(x,e_{x,n+2})\le C_A A^{-n-2}  \le \frac{8}{7} A^{-n-2}$, we have that the mass from $e_{x,n+2}$ stays within \[B(x, 3A^{-n-2}) \subset B(x,\frac{3}{A} s)\subset B(x,s/2).\]
	Therefore 
	\be \label{e:measb5}
	\nu(B(x,s)) \ge \mu_{n+2}(e_{x,n+2}).
	\ee
	By \eqref{e:mure2} and \eqref{e:ce2}, we obtain that  $\mu_{n+2}(e_{x,n+2})$ and $\mu_n(P_{n+1}(P_{n+2}(e_{x,n+2})))$ are comparable, where $P_{n+1}, P_{n+2}$ denote the predecessor as given in Definition \ref{d:pred}. By the triangle inequality, we obtain that $d(e_{x,n},P_{n+1}(P_{n+2}(e_{x,n+2}))) \le 4 A^{-n}$, and therefore by \eqref{e:ce1} and \eqref{e:muk}, we obtain that 
	$\mu_n(P_{n+1}(P_{n+2}(e_{x,n+2})))$ and $\mu_n(e_{x,n})$ are comparable. Combining the above with \eqref{e:measb5}, we obtain the lower bound $\nu(B(x,s)) \gtrsim \mu_n(e_{x,n})$ in \eqref{e:measb2}. This completes the proof of \eqref{e:measb2}.
	
	Next, we obtain \eqref{e:measb1} from \eqref{e:measb2}. 
	Let $0<s_1 < s_2  < A^{-4} \diam(Q_{l,0})$. Let $n_1,n_2 \in \bZ$ be such that 
	$A^{-n_i-1} \le s_i < A^{-n_i} $ for $i=1,2$.
	For $x \in Q_{l,0}$, let $x_{n_i} \in M_{n_i}$ be the unique point in $M_{n_i}$ such that $x_{n_i}\in Q_{n_i}(x)$.
	By \eqref{e:measb2} and Lemma \ref{l-ced}(a),(b),(c), we have
	\[
	\frac{\nu(B(x,s_2)) \operatorname{Cap}_{B(x,As_1)} (B(x,s_1))} {\nu(B(x,s_1)) \operatorname{Cap}_{B(x,As_2)}(B(x,s_2))} \asymp  \frac{\mu_{n_2}(x_{n_2}) c_{n_1}(x)}{\mu_{n_1}(x_{n_1})c_{n_2}(x)}  \asymp \frac{\mu_{n_2}(x_{n_2}) c_{n_1}(x_{n_1})}{\mu_{n_1}(x_{n_1})c_{n_2}(x_{n_2})}. 
	\]
  	Next, by using  Lemma \ref{L:ind}(2), we obtain
	\[
	(1-\delta)^{ n_2-n_1}\lesssim \frac{\nu(B(x,s_2)) \operatorname{Cap}_{B(x,As_1)} (B(x,s_1))} {\nu(B(x,s_1)) \operatorname{Cap}_{B(x,As_2)}(B(x,s_2))} \asymp  \frac{\mu_{n_2}(x_{n_2}) c_{n_1}(x_{n_1})}{\mu_{n_1}(x_{n_1})c_{n_2}(x_{n_2})} \lesssim C_4^{n_1-n_2},
	\]
 where $C_4>1$ is the constant $C$ in Lemma \ref{L:ind} and $\delta \in (0,1)$ is as in Lemma \ref{L:ESP}.
	The desired estimate  \eqref{e:measb1} follows by setting $\beta_1= - \log (1-\delta)/\log A$ and $\beta_2=\log C_4/\log A$.
	\qed
	
	\medskip
	
	 	We are now in the position to give the 
	 	
	\ms 
	
	{\sm {\em Proof of Theorem \ref{T:meas}.}} 
	The compact case follows by choosing $l=k_0-1$ in Proposition \ref{P:ballmeas}.
	
	It suffices to consider the non-compact case.
	For $l \le - 1, l \in \bZ$ let $\nu_{l}$ be the measure given by 
	Proposition \ref{P:ballmeas} on $Q_{l,0}$, and choose $a_n>0$ so that
	\[
	a_l \nu_{l}(B(x_0,1))=1, \q \mbox{for all $l \in \bZ, l <0$.}
	\]
 A compactness argument similar to that in  \cite{LuS} yields the existence
	of a measure $\nu$ which is a sub-sequential weak* limit of the sequence of measures $a_l \nu_{l}$ as $l \to -\infty$, bounded on
	compacts, such that it is  $(C_0, A, \beta_1, \beta_2)$-capacity good.
	\qed
	
	\subsection{A criterion for smoothness of measure}

	In this section, we will provide a useful sufficient condition for a doubling measure to be smooth.
	The definition of a smooth measure is given in Definition \ref{d:smooth}.	 
	 
	The following lemma  	follows immediately from \cite[Theorem 3.3.8]{CF} or \cite[Theorem 4.4.3]{FOT} and the countable subadditivity for capacities.  
	
	\begin{lem} \label{l:loc2glob}
		Let $(\sX,d,\mu, \sE,\sF)$ be a MMD space.
		Let $\set{B_i: i \in I}$ be a countable family of  open balls such that  $\cup_{i \in I} B_i=\sX$. 
		Let $U \subset \sX$.  
		Then $U$ has zero capacity for $(\sE,\sF)$ if and only if $U_i:=U \cap B_i$  has zero capacity for the part Dirichlet form 
		$(\sE^{B_i}, \sF^{B_i})$ for all $i \in I$.
	\end{lem}
	
	\begin{prop} \label{p:smooth}
		Let $(\sX,d,m,\sE,\sF)$ be a  MMD space that satisfies Assumption \ref{A:main2}. 
		Let $\mu$ be a $(C_0,A,\beta_1,\beta_2)$-capacity good measure on $\sX$ for some $C_0,A>1$ and $0< \beta_1 \le \beta_2$.
		Then $\mu$ is a smooth Radon measure on $\sX$.
	\end{prop}
	
	\proof
 By Theorems \ref{T:3.6} and \ref{T:Rgreen}, $(\sE, \sF)$ has regular Green functions. 
	Let $A$ denote the constant in Lemma \ref{l:nashwill}. 
	Using Lemma \ref{l:ratiocap}, we may assume that \eqref{e:m02} holds with this $A$. Let $B=B(x_0,r)$ denote any ball such that $r<\diam(\sX,d)/A^4$.
	For $x \in \sX, s< \diam(\sX,d)/A^4$, we set 
	$\Psi(x,s)= \frac{\mu(B(x,s))}{\operatorname{Cap}_{B(x,As)}(B(x,s))}$.

	We will show that $x \mapsto \int_B g_B(x,y) \, \mu(dy)$ is bounded uniformly in $B$. Note that $\int_B g_B (x, y) \mu (dy)$ is well defined for every $x\in B$ in view of the definition (see Definition \ref{D:goodgreen}) of regular Green function on $B$. By Lemma \ref{l:good-VD}, the measure $\mu$
satisfies the reverse volume doubling property (RVD) \cite[Exercise 13.1]{Hei}: that is  there exist $c_0 \in (0,1),  \alpha \in (0,\infty)$ such that
	\begin{equation} \label{e:rvd}
	\mu(B(x,R)) \ge c_0 \left(\frac{R}{r}\right)^\alpha \mu(B(x,r)), \q \mbox{for all $x \in \sX, 0<r \le R < \diam(\sX,d)$.} 
	\end{equation}
	In particular by letting $r \to 0$ in the above equation, we obtain $\mu(\{x\})=0$ for all $x \in \sX$.

	Fix $x \in B$ and set $B_i=B(x,A^{1-i}r), A_i=B_i \setminus B_{i+1}$ for $i \in \bN_{\ge 0}$.
	\begin{align} \label{e:sm0}
	\int_B g_B(x,y)\,\mu(dy) &\le  \int_B g_{B(x,Ar)}(x,y)\,\mu(dy) \mbox{ (by domain monotonicity)} \nonumber\\
	&= \sum_{i=0}^{\infty} \int_{B \cap A_i} g_{B(x,Ar)}(x,y)\,\mu(dy) \q \mbox{(since $\mu(\set{x})=0$ by (RVD))}\nonumber\\
	&\lesssim \sum_{i=0}^{\infty} \int_{B \cap A_i} \operatorname{Cap}_{B_0}(B_{i+1})^{-1}\,\mu(dy) \q \mbox{ (by (HG), \eqref{e:max}, Lemma \ref{L:3.5})}\nonumber\\
	&\lesssim  \sum_{i=0}^{\infty} \int_{B \cap A_i} \sum_{j=0}^i\operatorname{Cap}_{B_j}(B_{j+1})^{-1}\,\mu(dy) \q\mbox{ (by Lemma \ref{l:nashwill})}\nonumber\\
	&\lesssim \sum_{j=0}^\infty \operatorname{Cap}_{B_j}(B_{j+1})^{-1} \sum_{i=j}^\infty \int_{B \cap A_i}\,d\mu \lesssim \sum_{j=0}^\infty \operatorname{Cap}_{B_j}(B_{j+1})^{-1} \mu(B_j)\nonumber\\
	&\lesssim \sum_{j=0}^\infty \Psi(x,A^{-j}r) \lesssim \Psi(x_0,r) \q  \mbox{ (by  \eqref{e:m02}, Lemmas \ref{l:good-VD} and \ref{L:3.14}(b)).}
	\end{align}
 We claim that $\mu|_B $ is a smooth measure on $B$. Suppose not. Then there is a 
compact subset $K\subset B$ which is $\sE^B$-polar (which is equivalent to being $\sE$-polar)
so that $\mu (K)>0$.  Let $h(x):= \int_K g_B(x, y) \mu (dy)$. 
By \eqref{e:sm0}, $h(x)$ is bounded on $B$.
By Lemma \ref{L:3.5} and the maximum principle \eqref{e:max},
 $g_B(x, y)>0$ for $x\in B$ and $y\in K$ with $0<d(x, y)<r_1$ for some
 $r_1>0$. Thus  $\{x\in B:  h(x)>0\}$ has positive $\sE$-capacity. 
On the other hand, by the analog of \eqref{e:4.5} mentioned after \eqref{e:max1}, \eqref{e:sm0} and Fubini's theorem, $h(x)= \int_K g_B (x, y) \mu (dy)$ is  regular harmonic in $U$ with respect to $X^B$ for any open subset $U$ of $B$ with $K \cap \overline{U}=\emptyset$ and thus a bounded harmonic function in $B \setminus K$. Since $K$ is $\sE$-polar, $h$ is bounded harmonic in $B$.
Let $\{D_n; n\geq 1\}$ be an increasing sequence of relatively compact  open
subsets of $B$ that increases to $B$ and satisfies $K\subset D_1$.
 By Remark \ref{R:2.7} and Proposition \ref{P:3.2}, 
$h(x) = \bE^x[ h(X_{\tau_{D_n}})]$ for  $\sE$-q.e.~$x\in D_n$. 
Hence for  $\sE$-q.e.~$x\in B$, 
\begin{equation}\label{e:h}
h(x) = \lim_{n\to \infty} \bE^x [ h(X_{\tau_{D_n}})]
	= \lim_{n\to \infty} \int_K \bE^x g_B (X_{\tau_{D_n}}, y) \mu (dy).
\end{equation}
For every $y_0\in B$,   $x\mapsto g_B(x, y_0)$ is $X^B|_{B\setminus \sN_B} $-excessive, 
where $\sN_B$ is a Borel properly exceptional set for $X^B$ appearing in
Definition \ref{D:goodgreen}(iv) for $g_B$. Thus
$t\mapsto g_B (X^B_t, y_0)$ is a non-negative $\bP^x$-supermartingale 
for $x\in B\setminus \sN_B$, and  consequently, $\lim_{n\to \infty} g_B (X_{\tau_{D_n}}, y_0)$
exists $\bP^x$-a.s.  By the maximum principle \eqref{e:max}, EHI and H\"older estimate
\eqref{e:3.6}, $\bP^x$-a.s.~$\lim_{n\to \infty} g_B (X_{\tau_{D_n}}, y)$ exists in $\bR$ for any $y \in B(y_0,r_2)$ and is uniformly H\"older continuous in $y\in B(y_0, r_2)$ for some $r_2>0$. 
For any bounded compactly supported Borel function $\varphi \geq 0$ on $B$, by the strong Markov property and the bounded convergence theorem.
\begin{align*}
\lim_{n\to \infty} \bE^x \int_B g_B (X_{\tau_{D_n}}, y) \varphi (y)\, m(dy)
&= \lim_{n\to \infty}  \bE^{x} \left[ \int_0^{\tau_B}\varphi(X_s) \,ds 
 \circ \theta_{ \tau_{D_n}} \right]  \\
  &=\lim_{n\to \infty} \bE^x \int_{\tau_{D_n}}^{\tau_B} \varphi (X_s) ds =0.
\end{align*}
Thus by Fatou's lemma and the H\"older regularity mentioned above,
we conclude that 
$$\lim_{n\to \infty} g_B (X_{\tau_{D_n}}, y_0) = 0  \mbox{   for every $y_0\in B$ $\bP^x$-a.s.}$$
Observe that 
$\{ g_B (X_{\tau_{D_n}}, y); n\geq 1, y\in K\}$ are uniformly bounded random variables 
by  the maximum principle \eqref{e:max}. Thus we conclude from \eqref{e:h} by the bounded convergence theorem that 
$$
h(x) = \int_K \bE^x \left[ \lim_{n\to \infty} g_B (X_{\tau_{D_n}}, y) \right] \mu (dy)=0
\quad \hbox{for  $\sE$-q.e.~} x\in B\setminus \sN_B.
$$
This contradicts  the fact that $\{x\in B: h(x)>0\}$ has  positive $\sE$-capacity. 
We have thus proved that $\mu$ is a smooth Radon measure.  
	\qed

 A smooth measure $\mu$ on $\sX$ uniquely determines a positive continuous additive functional $A^\mu=\{A^\mu_t; t\geq 0\}$
of $X$. It can be used to define a time-changed process  $Y_t:=X_{\tau_t}$, where 
$$
\tau_t:=\inf\{r>0: A^\mu_r >t\}.
$$
Let $S(\mu)$ denote the quasi support of $\mu$ (see Definition \ref{d:qsupp}) and $F$ be the topological support of $\mu$. Clearly $S(\mu)  \subset F$ $\sE$-q.e. and $\mu(F\setminus S(\mu))=0$.
Suppose $\mu$ is a smooth Radon measure.  Then the  time-changed process $Y$,
 after possibly modification on a Borel properly exceptional set for $X$,  is a $\mu$-symmetric Hunt process on $F$ and its associated Dirichlet form $(\sE^\mu, \sF^\mu)$
on $L^2(F; \mu)$  is regular.  Moreover, 
\begin{align}\label{e:defTC}
	\sF^\mu &= \set{ \phi \in L^2(F,\mu): \phi=u \q\mu \mbox{-a.e.~ for some } u\in \sF_e},  \nonumber \\
	\sE^\mu(\phi,\phi) &= \sE(H_{S(\mu)}u, H_{S(\mu)}u), \q \mbox{for $\phi \in \sF^\mu$, and an arbitrary  $u\in \sF_e$ with
	   $\phi=u \q\mu$-a.e.,}  
	\end{align}  
	where $\sF_e$ is the extended Dirichlet space of $(\sX,d,m,\sE,\sF)$ and 
	$H_{S(\mu)}u(x)= \bE^x u(X_{\sigma_{S(\mu)}})$ for $x \in \sX$. 
	See \cite[Theorem 5.2.13]{CF} or \cite[Theorem 5.1.5 and Theorem 6.2.1]{FOT}.	
	The Dirichlet form $(\sE^\mu, \sF^\mu)$ is called the {\it  trace Dirichlet form} of $(\sE, \sF)$ on $L^2(F; \mu)$. 
	  If $\mu$ has full quasi support, then $\sF^\mu_e=\sF_e$ by \cite[Corollary 5.2.12]{CF}  and 
	  \eqref{e:defTC}  	  can be simplified to
	\be  \label{e:fullqs}
	\sF^\mu =\sF_e \cap L^2(\sX;\mu), \q \sE^\mu(u, u) = \sE(u,u) \q \mbox{for all  } u \in \sF_e.
	\ee
 
\begin{remark} \rm
 The above mentioned properties for time-changed processes and Dirichlet forms in fact hold for any smooth measure $\mu$ rather than just smooth Radon measures except that the time-changed process is a right process instead of being a Hunt process on $F$ 
and the trace Dirichlet form $(\sE^\mu, \sF^\mu)$ is quasi-regular on $L^2(S(\mu); \mu)$ instead of being regular on $L^2(F; \mu)$.
See \cite[Theorem 5.2.7]{CF}.   
\end{remark}
  
\medskip

 Recall the definition of  quasi support of a smooth measure from Definition \ref{d:qsupp}.  In this work, we 
 are interested in smooth measures with full  quasi support as defined below.	

\begin{definition} [Admissible smooth measures]{ \rm 
		Let $(\sX,d,m,\sE,\sF)$ be a  MMD space. We say that a smooth Radon measure $\mu$ on $\sX$ is 
		\emph{admissible} if $\mu$ has full quasi support. In particular, the  time-changed Dirichlet form is given by \eqref{e:fullqs}.
}
\end{definition}

	\begin{proposition} \label{p:fullquasisupport}
	Let
	 $(\sX,d,m,\sE,\sF)$ be a  MMD space that satisfies Assumption \ref{A:main2}. Let $\mu$ be a  $(C_0,A,\beta_1,\beta_2)$-capacity good (hence smooth) measure on $\sX$ for some $C_0,A>1$ and $0<\beta_1 \le \beta_2$. 
	Then $\mu$ is admissible.
	\end{proposition}
	
	\proof 	Let  $\sN$  be a Borel properly exceptional set for the Hunt process $X$ 
	associated with the regular Dirichlet form $(\sE, \sF)$
	on $L^2(\sX; m)$  so that the conclusion of Theorem \ref{T:Rgreen} holds. 
	Denote by $ S(\mu)$ a quasi support of $\mu$.
	As discussed in \cite[Proposition 2.6]{BK}, it suffices to show that
	\be \label{e:sm1}
	\bP^x(\sigma_{S(\mu)}=0)=1 \qq \mbox{for quasi every $x \in \sX$.}
	\ee 
	For the reader's convenience, we recall why \eqref{e:sm1} implies that $\mu$ has full quasi support. 
	By \cite[Theorem 4.6.1(i)]{FOT} we may assume that $S(\mu)^c$
	is nearly Borel and finely open, by adjusting $S(\mu)$ on a set of
	capacity zero. Then since $S(\mu)^c$ is nearly Borel and finely open,
 	for any  
	$x \in S(\mu)^c \setminus \sN$  
	we have $\bP^x ( \sigma_{ S(\mu) } >0 ) = 1$,
	which by \eqref{e:sm1}  implies   
	  that $S(\mu)^c$ has capacity zero.

 Note that  $(\sX,d,m,\sE,\sF)$ is irreducible by Theorem \ref{T:3.6}. 
 Applying  Lemma \ref{L:1} to the part process $X^{B(x_0, R_0)}$ of $X$ killed upon leaving a ball $B(x_0, R_0)$ whose complement 
has positive capacity, we have $\bP^x(\tau_x=0)=1$ for $\sE$-q.e.~$x \in B(x_0,R_0)$. By applying Lemma \ref{L:1} to countably many such balls $B(x_0,R_0)$, we conclude that
\begin{equation} \label{e:fqs0}
	\bP^x(\tau_x=0)=1 \quad \mbox{$\sE$-q.e.~$x \in \sX$.}
\end{equation}
Fix any $x \in \sX\setminus \sN$ that satisfies \eqref{e:fqs0}. Let $t>0$ and $\epsilon>0$ be arbitrary. By \eqref{e:fqs0}, we have
\be \label{e:fqs1}
\mathbb{P}^x(T<t ) > 1-\eps, \q \mbox{for $\sE$-q.e.~$x \in \sX$, where $T= \tau_{B(x,r)}$,}
\ee
for some $r=r(x,t,\epsilon)>0$.
By decreasing $r=r(x,t,\eps)$ if necessary, we may assume that $0<r<\diam(\sX,d)/A^4$, where $A$ is the constant in capacity good condition.
Since we will use Lemma \ref{L:3.5}, by increasing $A$ if necessary we assume that $A \ge 2K+1$, where $K$ is the constant in Assumption \ref{A:main}. This increase in $A$ is possible due to Lemma \ref{l:ratiocap}.
Fixing $r=r(x,t,\eps)$ as above, we define
\begin{align*}
K_1&=  B(x, A^{-1}r)   \cap  S(\mu).  
\end{align*}
We show that  there exists a  constant $c_0\in (0, 1)$  that depends only on the constants associated with
Assumption \ref{A:main2}, and capacity good condition such that
\begin{equation} \label{e:sm2}
\bP^x(\sigma_{K_1} < T) \ge c_0.
\end{equation}
Let $e$ denote the equilibrium measure for $K_1$  such that $e(\overline{K_1})=\operatorname{Cap}_{B}(K_1)$, where $B=B(x,r)$.
To prove \eqref{e:sm2}, we observe that
\be \label{e:sm3a} \bP^z(\sigma_{K_1} < \tau_B)= \int_{\ol{K_1}} g_B(z,y)\, e(dy) \quad \mbox{for $\sE$-q.e.~$z \in B$.}\ee
To obtain \eqref{e:sm3a}, we use \cite[Theorem 4.3.3 and the $0$-order version of Exercise 4.2.2]{FOT} to conclude that both sides of \eqref{e:sm3a} are quasi-continuous versions of the $0$-order equilibrium potential for $K_1$ with respect to the part Dirichlet form on $B$.
We would like to use \eqref{e:sm3a} for $z=x$, but $x$ could belong to the $\sE$-polar set associated with \eqref{e:sm3a}. To this end, we note that
both sides of \eqref{e:sm3a} are $X^B|_{B \setminus \sN}$-excessive from  \cite[Lemma A.2.4(ii)]{CF} and Theorem \ref{T:Rgreen} respectively.
 By the absolute continuity property of $X|_{\sX \setminus \sN}$ from Theorem \ref{T:Rgreen} and  
\cite[Theorem A.2.17(iii)]{CF}, we conclude that 
\be \label{e:sm3b} \bP^z(\sigma_{K_1} < \tau_B)= \int_{\ol{K_1}} g_B(z,y)\, e(dy) \quad \mbox{for every }  z \in B \setminus \sN.
\ee
It is crucial that the properly exceptional set $\sN$ does not depend on $B$. By \eqref{e:sm3b} and \eqref{e:max}, 
\be \label{e:sm3}
\bP^x(\sigma_{K_1} < T)= \int_{\overline{K_1}} g_B(x,y)\,e(dy) \ge g_B(x,A^{-1}r) \operatorname{Cap}_B(K_1).
\ee
By  \eqref{e:sm0}, \eqref{e:m02} and Lemma \ref{L:3.5}, there exists $C_1>0$ such that, 
\be \label{e:sm4}
\int_{B(x,r/A)} g_B(y,z)\,\mu(dz) \le C_1 g_B(x,r/A) \mu(B(x,r/A)) \quad\mbox{for all $y \in B(x,r)$.}
\ee
Since $\mu(S(\mu)^c)=0$, we obtain
\be \label{e:fqs2}
\mu(K_1)= \mu\left(B(x,A^{-1}r)  \right).
\ee 
We recall the following inequality for capacity: for any Radon measure $\nu$ on $B$ with $\int_Bg_B( \cdot,z)\,\nu(dz) \le 1$  $\sE$-q.e.~on $B$ and $\nu(B \setminus K_1)=0$, 
\[
\nu(K_1) \le \operatorname{Cap}_B(K_1).
\]
See \cite[p.441, Solution to Exercise 2.2.2]{FOT} and note also \cite[Exercise 4.2.2]{FOT}.
By considering the measure $\nu(\cdot)=\mu(K_1\cap \cdot)/(C_1 g_B(x,r/A) \mu(B(x,r/A)))$,  \eqref{e:sm4} and the above inequality, we obtain
\begin{align} \label{e:sm5}
\operatorname{Cap}_B(K_1)^{-1} &\le \nu(K_1)^{-1} = C_1 g_B(x,r/A) \mu(B(x,r/A))/\mu(K_1) \nonumber \\
&\lesssim    g_B(x,r/A)  \q \mbox{(by \eqref{e:fqs2}).}
\end{align}
Combining \eqref{e:sm3} and \eqref{e:sm5} establishes the claim \eqref{e:sm2}.
Choosing $\eps = c_0/2$, we obtain  
\begin{align*}
\mathbb{P}^x( \sigma_{S(\mu)} \le  t ) &\ge \bP^x (\sigma_{K_1} <T)- \bP^x (T \ge t) \q \mbox{(since $\set{\sigma_{K_1} <T} \subset \set{ \sigma_{S(\mu)} \le  t } \cup \set{T \ge  t}$)} \\
& > c_0 -\eps = \half c_0 \q \mbox{(by \eqref{e:fqs1} and \eqref{e:sm2}).}
\end{align*}
Since $t>0$ is arbitrary, the Blumenthal 0-1 law \cite[Lemma A.2.5]{CF} gives $\mathbb{P}^x( \sigma_{S(\mu)}=0 )=1$.
\qed

	\section{Quasisymmetry and stability}\label{S:7}

Although the assumption that all MMD spaces are strongly local is in force in this section,  we remark that 
 Lemma \ref{l:ehiqs} and Lemma \ref{l:picsqs}(a) in fact hold for general Dirichlet forms as well.

\ms 

The following is a straightforward consequence of the definition of quasisymmetry.
	
	\begin{lem} (\cite[Lemma 5.3]{BM1}) \label{l:ehiqs}
		Let $(\sX,d_1,\mu,\sE,\sF^\mu)$ be a MMD space  and let $d_2$ be a metric on $\sX$ quasisymmetric to $d_1$. If $(\sX,d_2,\mu,\sE,\sF^\mu)$ satisfies the EHI, then so does $(\sX,d_1,\mu,\sE,\sF^\mu)$.
	\end{lem}
	
	The next definition is a slight modification of \cite[Definition 5.4]{BM1}, 
	the change being made so that it applies to both compact and non-compact spaces.
	
	\begin{definition} \label{d:regscale} {\rm
			We say that a function $\Psi:\sX \times [0,\infty) \to [0,\infty)$ on a metric space $(\sX,d)$ is a \emph{regular scale function} if $\Psi(x,0)=0$ for all $x \in \sX$ and 
		there exist constants $C_1,   \beta_1,\beta_2>0$ such that, for all $x,y \in \sX$ and 
	 		finite $ 0< s \le r  \leq  \diam(\sX,d)$, 
		we have with    $R:=d(x,y)$,  $\Psi(y,s)>0$ and  
		\be \label{e:Psireg1}
		C_1^{-1} \Big(   \frac{r}{R \vee r}   \Big)^{\beta_2}  \Big( \frac{R \vee r}{s} \Big)^{\beta_1}
		\le  \frac{ \Psi(x,r)}{\Psi(y,s)} 
		\le C_1 \Big(   \frac{r}{R \vee r}   \Big)^{\beta_1}  \Big( \frac{R \vee r}{s} \Big)^{\beta_2}.
		\ee 
	 	}\end{definition}

\medskip

Given a  regular scale function 
$\Psi$ on $(\sX,d)$, we now define a metric $d_\Psi$. 
This is proved as in  \cite{BM1} -- the proof there still works when $\diam(\sX,d) < \infty$.

	\begin{proposition} (\cite[Proposition 5.7]{BM1}) \label{p:metric}
		Let $\Psi$  be a regular scale function on a metric space $(\sX,d)$. 
		There exists a metric $d_\Psi:\sX \times \sX \to [0,\infty)$ satisfying the following properties:
		\begin{enumerate}[\rm (a)]
			\item There exist $C,\beta >0$ such that for all $x , y \in \sX$,
		\be \label{e:defbeta}
			C^{-1} \Psi(x,d(x,y)) \le	d_\Psi(x,y)^\beta \le C \Psi(x,d(x,y)).
			\ee
			\item $d$ and  $d_\Psi$ are quasisymmetric.
			\item Assume in addition that $(\sX,d)$  (or equivalently $(\sX,d_\Psi)$) is uniformly perfect. Fix $A>1$.  Let $B_\Psi$ and $B$ denote metric balls in $(\sX,d_\Psi)$ and $(\sX,d)$ respectively. If $x \in \sX$ and $r,s>0$ satisfy, either $B_\Psi(x,s) \subset B(x,r) \subset B_\Psi(x,As)  
	\subsetneq \sX$		
			 or $B(x,r) \subset B_\Psi(x,s) \subset B(x,Ar)   \subsetneq \sX$, then there is a constant $C_1>1$ (which does not depend on  $x \in \sX, r>0, s>0$)  such that
			\be \label{e:qsm}
			C_1^{-1} s^\beta \le \Psi(x,r) \le C_1 s^\beta,
			\ee
			where $\beta >0$ is as given in \eqref{e:defbeta}.
		\end{enumerate}
	\end{proposition}

	We now introduce Poincar\'e, cutoff energy inequalities, and capacity bounds with respect to a regular scale function $\Psi$
	on $(\sX, d)$. 
	This is again a slight modification of  \cite[Definitions 5.8 and 5.13]{BM1}, so as to include both bounded and unbounded spaces.
	Recall that a \emph{cutoff function} $\vp$
	for $B_1 \subset B_2$ is any function $\vp \in \sF^\mu$ such that $0\le \vp \le 1$ in $\sX$, $\vp \equiv 1$ 
	in an open neighbourhood of $\overline{B_1}$, and $\operatorname{supp} \vp \subset B_2$. 
 Recall also that $\mu_{\<f\>}$ is the energy measure of $f\in \sF^\mu$; see Section \ref{S:2}.
	
	\begin{definition} \label{d:pi-cs}  \rm
	Let  $\Psi$ be a  regular scale function on $(\sX,d)$, and $(\sX,d,\mu,\sE,\sF^\mu)$ a MMD space. 
	\begin{description}
	\item{(i)}
		We say that $(\sX,d,\mu,\sE,\sF^\mu)$ satisfies the \emph{Poincar\'e inequality} 
		\hypertarget{pi}{$\operatorname{PI}(\Psi)$}, if there exist constants $C,A_1,A_2 \ge 1$ such that 
for all $x\in \sX$, $R \in (0, \diam(\sX,d)/A_2)$ and $f \in \sF^\mu$,  $\mu(B(x,R))<\infty$ and
\be  \tag*{$\operatorname{PI}(\Psi)$}
\int_{B(x,R)} (f - \ol f)^2 \,d\mu  \le C \Psi(x,R) \, 
 \mu_{\<f\>} ( B(x,A_1R)) , 
\ee
where $\ol f= \frac1{\mu(B(x,R))} \int_{B(x,R)} f\, d\mu$. 

\item{(ii)}We say that $(\sX,d,\mu,\sE,\sF^\mu)$ satisfies the 
\emph{cutoff energy inequality} \hypertarget{cs}{$\operatorname{CS}(\Psi)$}, 
if there exist $C_1,C_2>0,A_1,A_2>1$ such that the following holds.
For all $R\in (0,\diam(\sX,d)/A_2)$, $x \in \sX$ with $B_1=B(x,R)$ and $B_2=B(x,A_1R)$, there exists a cutoff function $\vp$
for $B_1 \subset B_2$ such that for any $u \in \sF^\mu \cap L^\infty$,
\be  \tag*{$\operatorname{CS}(\Psi)$}
 \int_{\sX} u^2 d \mu_{\< \vp\>}  \le C_1 \mu_{\<u \>} (B_2 \setminus B_1  ) 
+ \frac{C_2}{\Psi(x,R)} \int_{B_2 \setminus B_1} u^2 \,d\mu.
\ee 

\item{(iii)} We say that  $(\sX,d,\mu,\sE,\sF^\mu)$ satisfies the \emph{capacity estimate} \hypertarget{cap}{$\operatorname{cap}(\Psi)$} 	
if there exist positive constants $C_1, A_1,  A_2>1$ such that for all $R\in (0,\diam(\sX,d)/A_2)$ and  $x \in \sX$ 
\be \tag*{$\operatorname{cap}(\Psi)$}
C_1^{-1} \frac{\mu(B(x,R))}{\Psi(x,R)} \le   \operatorname{Cap}_{B(x,A_1R)} \left(B(x,R) \right) \le C_1 \frac{\mu(B(x,R))}{\Psi(x,R)}.
\ee
If $\Psi(r)=r^\beta$, we denote $\operatorname{PI}(\Psi), \operatorname{CS}(\Psi), \operatorname{cap}(\Psi)$ by $\operatorname{PI}(\beta), \operatorname{CS}(\beta), \operatorname{cap}(\beta)$ respectively.	
\end{description}
 \end{definition}
	
	The following lemma shows that the Poincar\'e and cutoff energy inequalities take a much simpler form 
	with respect to the  metric $d_\Psi$. 
	
	\begin{lem} (\cite[Lemma 5.9]{BM1})\label{l:picsqs}
		Let $(\sX,d,\mu,\sE,\sF^\mu)$ be a uniformly perfect MMD space and let $\Psi$ be a regular scale function.  Let $d_\Psi$ be the metric constructed in Proposition \ref{p:metric} with $\beta>0$ as given in \eqref{e:defbeta}.
 		 Then
		\begin{enumerate}[(a)]
			\item $(\sX,d,\mu,\sE,\sF^\mu)$ satisfies \hyperlink{pi}{$\operatorname{PI}(\Psi)$} if and only if $(\sX,d_\Psi,\mu,\sE,\sF^\mu)$ satisfies \hyperlink{pib}{$\operatorname{PI}(\beta)$}.
			\item $(\sX,d,\mu,\sE,\sF^\mu)$ satisfies \hyperlink{cs}{$\operatorname{CS}(\Psi)$} if and only if $(\sX,d_\Psi,\mu,\sE,\sF^\mu)$ satisfies \hyperlink{csb}{$\operatorname{CS}(\beta)$}.
		\end{enumerate}
	\end{lem}

The following comparison of annuli follows readily from the definition.

\begin{lem} (\cite[Lemma 1.2.18]{MT}) \label{l:cann}
	Let the identity map $\operatorname{Id}:(\sX,d_1) \to (\sX,d_2)$ be an $\eta$-quasisymmetry 
	for some distortion function $\eta$. 
	Then for all $A>1,x \in \sX, r>0$, there exists $s>0$ such that, with $B_i$ denoting balls in $(\sX,d_i)$
	\be \label{e:ann1}
	B_2(x,s) \subset B_1(x,r) \subset B_1(x,Ar) \subset B_2(x, \eta(A)s).
	\ee
	In \eqref{e:ann1}, $s$ can be defined as
	\[
	s= \sup \set{0\le s_2 < 2\diam(\sX,d_1): B_2(x,s_2) \subset B_1(x,r) }
	\]
	Moreover, for all $A >1$, $x \in \sX$ and $r >0$, there exists $t>0$ such that
	\be \label{e:ann2}
	B_1(x,r) \subset B_2(x,t) \subset B_2(x,At) \subset B_1(x,A_1r),
	\ee
	where $A_1= {1}/{\eta^{-1}(A^{-1})}$. In \eqref{e:ann2}, $t$ can be defined as
	\[
	t= A^{-1}\sup \set{0 \le r_2 < 2 A \diam (\sX,d_2): B_2(x,Ar_2) \subset B_1(x,A_1 r)}.
	\]
\end{lem}

The following is an analogue of Lemma \ref{l:picsqs} for the capacity estimate \hyperlink{cap}{$\operatorname{cap}(\Psi)$}.
\begin{lem} \label{l:qscap}
	Let $(\sX,d,\mu,\sE,\sF^\mu)$ be a MMD space that satisfies the EHI and let $\Psi$ be a regular scale function. Suppose  that $(\sX,d)$ is complete and that $\mu$ satisfies  the volume doubling property on $(\sX,d)$. Let $d_\Psi$ be the metric constructed in Proposition \ref{p:metric} with $\beta>0$ as given in \eqref{e:defbeta}. Then $(\sX,d,\mu,\sE,\sF^\mu)$ satisfies \hyperlink{cap}{$\operatorname{cap}(\Psi)$} if and only if $(\sX,d_\Psi,\mu,\sE,\sF^\mu)$ satisfies \hyperlink{cap}{$\operatorname{cap}(\beta)$}.
\end{lem}
\proof
Let $B$ and $B_\Psi$ denote balls in the metrics $d$ and $d_\Psi$ respectively.
 By Lemma \ref{l:ehiqs},  $(\sX,d_\Psi,\mu,\sE,\sF^\mu)$ also satisfies the EHI. 
  Let the identity map $\on{Id}:(\sX,d_\Psi) \to (\sX,d)$ be an $\eta$-quasisymmetry.  Note that $\mu$ satisfies the volume doubling property with respect to the metrics $d$ and $d_\Psi$.
 
Let $(\sX,d,\mu,\sE,\sF^\mu)$ satisfy \hyperlink{cap}{$\operatorname{cap}(\Psi)$}.  Set $A_1= \eta(2)$ and choose $\eps\in(0,\frac 1 2]$ so that $\eta(4 \eps) \le \frac{1}{2A_1}$. By Lemma \ref{l:ratiocap}, we may assume that 
 \be \label{e:qsc1}
 \operatorname{Cap}_{B(x,A_1 r)}(B(x,r)) \asymp \frac{\mu(B(x,r))}{\Psi(x,r)}, \q \mbox{for all $x \in \sX, 0<r \lesssim \diam(\sX,d)$.}
 \ee
 By Lemma \ref{l:cann} and Proposition \ref{p:metric}(c), for all $x \in \sX$, $0< s < \eps \diam(\sX,d_\Psi)$, there exists $r>0$ such that	$B(x,r) \subset B_\Psi(x,s) \subset B_\Psi(x,2s) \subset B(x,\eta(2)r)$  and $s^\beta \asymp \Psi(x,r)$. Note that  $B(x,\eta(2)r)\neq \sX,$  since $r \le  \operatorname{diam}(B_\Psi(x,s),d)< \eta(4 \eps) \diam(\sX,d)$ by \cite[Proposition 10.8]{Hei}. By the volume doubling property $\mu(B(x,r)) \asymp \mu(B_\Psi(x,s))$. By domain monotonicity and \eqref{e:qsc1}, we have
 \be \label{e:qsc2}
  \operatorname{Cap}_{B_\Psi(x,2s)}(B_\Psi(x,s)) \ge  \operatorname{Cap}_{B(x,A_1 r)}(B(x,r)) \asymp  \frac{\mu(B(x,r))}{\Psi(x,r)} \asymp  \frac{\mu(B_\Psi(x,s))}{s^\beta}, 
 \ee
 for all $x \in \sX, 0<s \lesssim \diam(\sX,d_\Psi)$.
 
 Set $A_2 = 1/\eta^{-1}(A_1^{-1})$.  By Lemma \ref{l:cann} and Proposition \ref{p:metric}(c), for all $x \in \sX, s \in (0,(2A_2)^{-1}\diam(\sX,d_\Psi))$, there exists $r>0$ such that 
 $B_\Psi(x,s) \subset B(x,r) \subset B(x,A_1 r) \subset B_\Psi(x,A_2 s) \neq \sX$ and $\Psi(x,r) \asymp s^\beta$. By the volume doubling property, $\mu(B(x,r)) \asymp \mu(B_\Psi(x,s))$.
 By Lemma \ref{l:ratiocap},  \cite[Proposition 10.8]{Hei}, domain monotonicity and \eqref{e:qsc1}, we have
 \begin{align} \label{e:qsc3}
  \operatorname{Cap}_{B_\Psi(x,2s)}(B_\Psi(x,s)) &\asymp \operatorname{Cap}_{B_\Psi(x,A_2s)}(B_\Psi(x,s)) \nonumber \\ &\le \operatorname{Cap}_{B(x,A_1 r)}(B(x,r)) 
 \asymp  \frac{\mu(B(x,r))}{\Psi(x,r)} \asymp  \frac{\mu(B_\Psi(x,s))}{s^\beta}
 \end{align}
  for all $x \in \sX, 0<s \lesssim \diam(\sX,d_\Psi)$. By \eqref{e:qsc2} and \eqref{e:qsc3}, $(\sX,d_\Psi,\mu,\sE,\sF^\mu)$ satisfies \hyperlink{cap}{$\operatorname{cap}(\beta)$}.
  
  The converse follows from a similar argument.
\qed

	We will now apply these results in the context of a change of measure on a MMD space.
	Let  $(\sX,d,m,\sE,\sF)$ be a MMD space 
	which satisfies the EHI and   one (and hence all) of the three equivalent conditions in Theorem \ref{T:ehitomd}.
	Let $(\sE,\sF_e)$ be its  corresponding extended Dirichlet space, and 
	$\mu$ be the measure constructed in Theorem \ref{T:meas}. 
	By Propositions \ref{p:smooth} and \ref{p:fullquasisupport},  $\mu$ is a positive Radon measure charging no set of capacity zero and possessing full quasi-support. 
	Let $(\sE^\mu,\sF^\mu)$ denote the  time-changed Dirichlet space with respect to $\mu$ as defined in \eqref{e:defTC}. 
  	We have   
	$\sF^\mu=\sF_e \cap L^2(\sX,\mu)$,
    $\sE^\mu(f,f)= \sE(f,f)$ for all $f \in \sF^\mu$, and $\sF^\mu_e=\sF_e$
 (cf. \cite[Theorem 5.2.2, (5.2.17) and Corollary 5.2.12]{CF}). Moreover,  the Dirichlet form $(\sE^\mu, \sF^\mu)$ on $L^2(\sX; \mu)$ shares the same quasi notions as the original 
 Dirichlet form $(\sE, \sF)$ on $L^2(\sX; m)$; see \cite[Theorem 5.2.11]{CF}. 
	
	\begin{thm} \label{T:PI-CS}
		Let $(\sX,d)$ be complete and metric doubling.
		Suppose that $(\sX,d,m,\sE,\sF)$ is a MMD space 
		which satisfies the EHI.  
		Let $\mu$ be  a $(C_0, A, \beta_1, \beta_2)$-capacity good measure with $A \ge 2^{1/4}$. 
	 	 Denote $D=\diam(\sX,d)$.
		Then the function $\Psi:\sX \times [0,\infty) \to [0,\infty)$ defined by $\Psi(x,0)=0$ and 
		\be
		\label{e:Psidef} \Psi(x,r) = \begin{cases}
 			\frac{\mu(B(x,r))}{ \operatorname{Cap}_{B(x, r/A^4)}(B(x,r/A^5))}  &\mbox{ if } 0<r <D,\\
				\frac{\mu(B(x,D))}{ \operatorname{Cap}_{B(x, D/A^4)}(B(x,D/A^5))}  &\mbox{ if  }
				r\geq D \hbox{ and }  D<\infty, 
		\end{cases} 
		\ee
		is a regular scale function on $(\sX, d)$.  Furthermore, the MMD space  $(\sX,d,\mu,\sE,\sF^\mu)$ satisfies  the Poincar\'e inequality \hyperlink{pi}{$\operatorname{PI}(\Psi)$}, the cutoff 
		energy inequality  \hyperlink{cs}{$\operatorname{CS}(\Psi)$} and the capacity estimate 
		\hyperlink{cap}{$\operatorname{cap}(\Psi)$}.
	\end{thm}
	\proof
	By volume doubling (Lemma \ref{l:good-VD}) and Lemma \ref{L:3.14}(c), there exists $C_2>0$ such that for all finite $0 < r \le D$ and for all $x,y \in \sX$ with $d(x,y) \le r$, we have
	\be \label{e:Psi_comp}
	C_2^{-1} \Psi(x,r) \le \Psi(y,r) \le C_2  \Psi(x,r).
	\ee
     Let $x,y \in \sX$ and set $R:=d(x,y)$.
	If $R\le r$ the inequalities in \eqref{e:Psireg1} are immediate from  \eqref{e:m02} and \eqref{e:Psi_comp}.
	If $s \le r<R$, then writing
	\bes
	\frac{ \Psi(x,r)}{\Psi(y,s)} =  \frac{ \Psi(x,r)}{\Psi(x,R)}  . \frac{ \Psi(y,R)}{\Psi(y,s)} .  \frac{ \Psi(x,R)}{\Psi(y,R)}, 
	\ees
	and bounding each of the three terms on the right using  \eqref{e:m02} and 
	\eqref{e:Psi_comp} give \eqref{e:Psireg1}. Thus $\Psi$ is a regular scale function.

	By Lemma \ref{l:capcomp}, the MMD space  $(\sX,d,\mu,\sE,\sF^\mu)$ satisfies \hyperlink{cap}{$\operatorname{cap}(\Psi)$}.

		Let $d_\Psi$ and $\beta>0$ be as given by Proposition \ref{p:metric}. By Lemma \ref{l:qscap}, the MMD space  $(\sX,d_\Psi,\mu,\sE,\sF^\mu)$ satisfies  
	\hyperlink{cap}{$\operatorname{cap}(\beta)$}.
	By Lemma \ref{l:ehiqs} and Proposition \ref{p:metric}(b), $(\sX,d_\Psi,\mu,\sE,\sF^\mu)$ satisfies the EHI. 
	
	By Lemma \ref{l:metric} the space $(\sX,d_\Psi)$ is uniformly perfect, and hence the measure $\mu$
	on $(\sX,d_\Psi)$ satisfies the reverse volume doubling property (RVD) as defined in \eqref{e:rvd} \cite[Exercise 13.1]{Hei}.

	Therefore by \cite[Theorem 1.2]{GHL}, since $(\sX,d_\Psi,\mu,\sE,\sF^\mu)$  satisfies  
	the EHI and  \hyperlink{cap}{$\operatorname{cap}(\beta)$}, 
it satisfies 
	\hyperlink{pib}{$\operatorname{PI}(\beta)$}  and \hyperlink{csb}{$\operatorname{CS}(\beta)$}. 
	We now conclude using Lemma \ref{l:picsqs}. \qed

 \ms 

The following gives equivalent characterization of the EHI for a MMD space $(\sX,d,m,\sE,\sF)$.

	\begin{thm} \label{T:main-new}
		Let $(\sX,d)$ be a complete, metric doubling, connected metric space with a strongly local 
		regular Dirichlet form $(\sE,\sF)$ on $L^2(\sX; m)$. 
		The following are equivalent: 
		\begin{enumerate}[\rm (a)]
		\item  $(\sX,d,m,\sE,\sF)$ satisfies the EHI.  
		
		\item There exist an admissible smooth doubling Radon measure $\mu$ on $(\sX,d)$ and a regular scale function $\Psi$ 
		such that the time-changed MMD space  $(\sX,d,\mu,\sE,\sF^\mu)$ satisfies the Poincar\'e inequality \hyperlink{pi}{$\operatorname{PI}(\Psi)$}  and the cutoff 
		energy inequality  \hyperlink{cs}{$\operatorname{CS}(\Psi)$}. 
		
		\item There exist an admissible smooth doubling Radon measure $\mu$ on $(\sX,d)$, a  metric $d_\Psi$ on $\sX$ that is quasisymmetric to $d$, and $\beta>0$, such that the time-changed MMD space   $(\sX,d_\Psi,\mu,\sE,\sF^\mu)$ satisfies Poincar\'e inequality \hyperlink{pib}{$\operatorname{PI}(\beta)$}  and the cutoff 
		energy inequality  \hyperlink{csb}{$\operatorname{CS}(\beta)$}.
          \end{enumerate}
	\end{thm}
	
	\proof 
	(a) $\Rightarrow $(b) This is immediate  from Lemma \ref{l:good-VD}, Propositions \ref{p:smooth}, \ref{p:fullquasisupport}, Theorems \ref{T:meas} and \ref{T:PI-CS}. 
	
	\smallskip
	
	(b)$\Rightarrow$(c) By Lemma \ref{l:metric}(b), $(\sX,d)$ is uniformly perfect. Let $d_\Psi$ and $\beta>0$ be as given by Proposition \ref{p:metric}. 
	Quasisymmetry of $d_\Psi$ follows from Proposition \ref{p:metric}(b). 
	Then \hyperlink{pib}{$\operatorname{PI}(\beta)$} and \hyperlink{csb}{$\operatorname{CS}(\beta)$}  for $(\sX,d_\Psi,\mu,\sE,\sF^\mu)$  follow from Lemma \ref{l:picsqs}. 
	
	\smallskip
	
	(c)$\Rightarrow$(a)
	 	 By Lemma \ref{l:metric}(b),	$(\sX,d_\Psi)$ is uniformly perfect. 
	Thus $\mu$ satisfies the reverse volume doubling property \eqref{e:rvd} \cite[Exercise 13.1]{Hei}. Since   $(\sX,d)$  is metric doubling, so is $(\sX,d_\Psi)$ \cite[Theorem 10.18]{Hei}.
	So by \cite[Proposition 5.11 and Remark 5.12]{BM1},  we obtain the condition  $\operatorname{(Gcap_\le)_\beta}$  in \cite{GHL}. 
	Then by  the implication  $\operatorname{(Gcap_\le)_\beta}$  plus  \hyperlink{pib}{$\operatorname{PI}(\beta)$}  to  the EHI in
	\cite[Theorem 1.1]{GHL}, we obtain the EHI for $(\sX,d_\Psi,\mu,\sE,\sF^\mu)$. 
	Since $d_\Psi$ and $d$ are quasisymmetric, the desired  EHI follows from Lemma \ref{l:ehiqs}.
	\qed

 \begin{remark}\label{R:7.10} 
{\rm
	\begin{enumerate}[(i)]
		\item 
	 Note that conditions (b) and (c) in the Theorem above do not include the requirement that 
 $(\sX,d,m,\sE,\sF)$  satisfies the conditions (HC) or (Ha) introduced in Section \ref{S:lreg}.
 (It would be undesirable to include (Ha) or (HC), since we do not know if they are stable.) 
Thus (b) or (c)  does not immediately give the existence of Green's functions; 
however 
the existence of regular Green functions 
does follow from the implications (b), (c) $\Rightarrow$ (a) and Theorems \ref{T:3.6} and
 \ref{T:Rgreen}.

The proof in \cite{GHL} that   $\operatorname{(Gcap_\le)_\beta}$ plus  \hyperlink{pib}{$\operatorname{PI}(\beta)$} implies
the EHI does not require the existence of Green's functions.

\item  The result that (a) implies (c) in Theorem \ref{T:main-new} can be sharpened as follows. If $(\sX,d,m,\sE,\sF)$ satisfies the EHI  then for any $\beta>2$ there exists a  metric $d_\Psi$ on $\sX$ that is quasisymmetric to $d$, and an admissible smooth doubling Radon measure $\mu$ such that the time-changed MMD space $(\sX,d_\Psi,\mu,\sE,\sF^\mu)$ satisfies Poincar\'e inequality \hyperlink{pib}{$\operatorname{PI}(\beta)$}  and the cutoff 
energy inequality  \hyperlink{csb}{$\operatorname{CS}(\beta)$}. The condition $\beta>2$ is sharp in the sense that any $\beta$ in property (c)  necessarily satisfies $\beta \ge 2$ and there are examples for which $\beta=2$ is not possible. These results are contained in \cite{KM}.
\end{enumerate}
} \end{remark}

\sm {\em 	Proof of Theorem \ref{T:main-stable}.}	
The condition that $\sE(f,f) \asymp \sE'(f,f)$ for all $f \in \sF$ implies that the associated
energy measures satisfy   
$\mu_{\<f\> } \asymp \mu'_{\<f\>}$;
 see \eqref{e:4.13}.
Hence 
the condition (b) in Theorem \ref{T:main-new} 
holds for $\sE'$ by Theorems \ref{T:ehitomd} and \ref{T:main-new}, and therefore the implication (b) $\Rightarrow$ (a) in Theorem \ref{T:main-new}
implies that the EHI holds for $\sE'$. \qed

\medskip

\medskip

 The following is an extension of Theorem \ref{T:main-stable}, where the symmetrizing measures for the Dirichlet forms may 
 be different.   
 
 	\begin{thm} \label{T:7.11}
		Let $(\sX,d)$ be a complete, doubling metric space, 
		and let $m$ be a   Radon measure on $\sX$ with full support.
		Let  
		$(\sE,\sF)$ be a strongly local  regular Dirichlet form on $L^2(\sX; m)$. 
		Suppose that   $(\sX,d,m,\sE,\sF)$ satisfies the EHI. 
		Let $\mu$ be a smooth Radon measure of $(\sX, d, m, \sE, \sF)$ with full quasi support on $\sX$,
		and   $(\sE', \sF')$ be  another strongly local regular
	Dirichlet form on $L^2(\sX; \mu)$ such that $\sF\cap C_c(\sX) = \sF'\cap C_c(\sX)$ and 
\begin{equation}\label{e:7.11}
 C^{-1}  \sE(f,f) \le  \sE'(f,f) \le C  \sE(f,f)   \quad \hbox{ for all  } f \in \sF \cap C_c(\sX) 
\end{equation} 
for some $C \ge 1$.
		Then $(\sX,d, \mu,\sE',\sF')$ satisfies the EHI. 
\end{thm} 

\pf  Let $X$ be the Hunt process associated with the regular Dirichlet form $(\sE, \sF)$ on $L^2(\sX; m)$.
 Since $\mu$ is a smooth Radon measure with full quasi-support,  its associated positive continuous additive functional
 $A_t$ is strictly increasing up to the lifetime of $X$. Thus its  time-changed process $Y_t:=X_{\tau_t}$, with
 $\tau_t:=\inf\{r>0: A_r >t\}$, has the same family of harmonic functions as that of $X$. To see this, first note by Proposition \ref{P:2.10} it suffices to show that $Y_t$ and $X_t$ have the same family of bounded harmonic functions which in turn follows from Proposition \ref{P:2.9}.
 
 By \eqref{e:fullqs}, the Dirichlet form $(\sE^\mu, \sF^\mu)$ of the time-changed process $Y$ is regular on $L^2(\sX; \mu)$
 and has the property that $\sF^\mu=\sF_e \cap L^2(\sX; \mu)$, $\sF^\mu_e=\sF_e$
   and $\sE^\mu=\sE$ on $\sF_e$.
    Moreover,  $(\sE^\mu, \sF^\mu)$ is strongly local and  satisfies the EHI.
   Since  $\sF=\sF_e \cap L^2 (\sX;m)$  and both $m$ and $\mu$ are Radon,  
      we have by \eqref{e:fullqs} that
 $$ 
 \sF \cap C_c(\sX) = (\sF_e \cap L^2(\sX; m))\cap C_c(\sX) = (\sF^\mu_e \cap L^2(\sX; \mu))\cap C_c(\sX)
 =\sF^\mu \cap C_c(\sX) .
 $$
  Since $\sF'\cap C_c(\sX)=\sF^\mu \cap C_c (\sX)$ is   dense in $\sF'$
and $\sF^\mu$ with respect to the Hilbert  norms  $\sqrt{\sE'_1}$ and $\sqrt{ \sE^\mu_1}$,  respectively,
where 
$$
\sE'_1(u, u):=\sE'(u, u)+ \int_{\sX} u(x)^2 \mu (dx) \quad \hbox{and} \quad 
  \sE^\mu_1(u, u):=\sE^\mu(u, u)+ \int_{\sX} u(x)^2 \mu (dx) ,
 $$
 we have by \eqref{e:7.11} that $\sF'=\sF^\mu$ and 
   $$
 C^{-1}  \sE^\mu (f, f) \le  \sE'(f,f) \le C  \sE^\mu (f, f)   \quad \hbox{ for all  } f \in \sF^\mu.   
  $$
  The desired conclusion of the theorem now follows from Theorem \ref{T:main-stable} applied 
 to the MMD space $(\sX, d, \mu, \sE^\mu, \sF^\mu)$. 
 \qed 
 
   \begin{remark}\label{R:7.12}    \rm
The stability results of this paper, Theorem \ref{T:main-stable} and Theorem \ref{T:7.11},
   hold for the  ${\rm EHI}_{\rm loc}$
   as well.  We now indicate the needed modifications.
     All of the results of Section \ref{S:5}
    extend easily under the assumption  
   ${\rm EHI}_{\rm loc}$
    except that the conclusions only hold for balls of 
    small enough radii. 
 The main difference is in the construction of the measures $\nu_l$ in Proposition \ref{P:ballmeas}. Instead of the initial condition on $M_{l+3}$ for the inductive construction using Lemma \ref{L:ind}, we set the initial condition on $M_1$ to the uniform probability measure on $M_1$, where $M_1$ is as given in the generalized dyadic decomposition of  $Q_{l,0}$. Then the weak* subsequential limit as in the proof of Theorem \ref{T:meas} will be a capacity good measure (only at small enough scales using the same argument). 
 However, this property is enough so that our construction gives a smooth measure with full quasi support. 
 All the results used in Section 7 (for example, \cite[Theorem 1.2]{GHL}) will also admit local versions. 
 As noted in \cite[Remark 4.6]{GT12},  \cite[Proof of Theorem 4.2]{GT12} cannot be localized, but \cite[Theorems 6.2 and 7.3]{GK} give a localization of it.
 Although  there is no clear reference in the literature for these results,  a careful reading of the proofs in the literature 
shows that these local versions do hold, with essentially the same proof.
 \end{remark}
	
\section{Examples}\label{S:8}

\begin{example}\label{E:instab1}
 {\rm 
   The following example uses the instability of the Liouville property given in Lyons \cite{Lyo}
to show that without some regularity of the metric the EHI is not stable.

\smallskip

We begin by describing briefly Lyons' example. 
Let $\Gam = (\bV_\Gam, E_\Gam)$ be the free group with two generators $a$ and $b$, and let
 $\bV = \bV_\Gamma \times \{0,1\}$. Lyons  constructed two sets of  symmetric conductances $\{a^{(i)}_{xy}, x,y \in \bV\}$, 
$i =1,2$,  on $\bV$ such that if $E_i =\{ (x,y):  a^{(i)}_{xy} >0\}$  then $E_1 = E_2$. 
 Denote $E_G=E_1=E_2$, and let $\bG = (\bV, E_G)$ be the associated graph.
 The two sets of conductances 
have the following additional properties:
\begin{description}
	\item[(0)] $G=(\bV,E_G)$ is connected.
\item[(1)] For each $x \in \bV$, $4 \le \big|\{y: a^{(i)}_{xy}>0 \} \big| \le 8$, so  the graph $\bG$ has uniformly bounded vertex degrees. 

\item[(2)] There exists $p_0 \in (0,1)$ so that $a^{(i)}_{xy} \in \{0\} \cup [p_0, 1]$ for all $x,y$, $i=1,2$. Thus 
the conductances $a^{(i)}_{xy}$ are uniformly bounded above and below on the graph $\bG$. 
\end{description}
Define the quadratic forms, for $f : \bV \to \bR$,
\be 
  \sE^{(i)}
  (f,f) = \half \sum_{x \in \bV}\sum_{y \in \bV}  a^{(i)}_{xy} (f(y)-f(x))^2 , \q i=1, 2.
\ee
 In view of (2) above we have
\be \label{e:8.2}
  p_0  \sE^{(1)} (f,f)  \le  \sE^{(2)} (f,f)  \le p_0^{-1}   \sE^{(1)}(f,f)  \quad \hbox{for any function $f$ on } \bV.
 \ee
  Let $d$ be the graph distance on $\bV$ and 
 $m$ be the measure on $\bV$ which assigns mass 1 to each vertex $x$.   
 In view of \cite[Proposition 1.21(c)]{B2},  $(\sE^{(i)}, \sF)$ is a symmetric regular Dirichlet form on
$L^2(\bV,m)$ for $i=1, 2$, where $\sF:= L^2(\bV,m)$.

Recall that we say that the strong Liouville property (SLP) holds (for a MMD space) if every non-negative harmonic function is
constant, and the Liouville property (LP) holds if every bounded harmonic function is constant. Lyons
 constructed $\{a^{(i)}_{xy}\}$ for $i=1, 2$
so that the SLP holds for  
the MMD space $(\bV,  d, m, \sE^{(1)}, \sF)$ while the LP fails for the MMD $(\bV,  d, m, \sE^{(2)}, \sF)$. 
Thus   Lyons' example 
shows that neither SLP nor LP are stable.

Let $(\sX, \tilde d)$ be the cable system for the graph $\bG$: each edge $e \in E_G$ is replaced by a copy
of $[0,1]$ -- see for example \cite{V} for details of the construction. Let $\mu$ be the measure which assigns
a copy of Lebesgue measure to each cable. The metric $\tilde d$ is the unique length metric on $\sX$ 
which equals Euclidean distance on each cable. Write $(\tilde \sE^{(i)}, \tilde \sF)$ for the (regular and strongly local) Dirichlet forms on $\sX$
associated with the cable system.  (For further details of this construction see \cite{V, BM1}.)

Now let $d'(x,y) = 1 \wedge \tilde d(x,y)$ for $x,y \in \sX$, and write $B'(x,r)$ for balls with respect to the metric $d'$.
  Then $(\sX, d')$ is locally compact and complete, and
$ (\sX,d',\mu,\tilde \sE^{(i)}, \widetilde{\sF})$, $i=1, 2$,  are MMD spaces.
Note that (MD) fails for this space.
  The instability of the SLP and LP for the graph $\bG$  extends to the cable systems,
 so that the MMD space $ (\sX,d',\mu,\tilde \sE^{(1)}, \widetilde{\sF})$  
  satisfies the SLP while  $(\sX,d',\mu,\tilde \sE^{(2)}, \widetilde{\sF})$ fails to satisfy the LP.

\smallskip

We now consider the EHI for balls $B'(x,r/2) \subset B'(x,r)$. Note first  that if $r \le 2$ then since each vertex of $\bV$ has 
between 4 and 8 neighbours, the EHI follows from the local Harnack inequality for both 
$ (\sX,d',\mu,\tilde \sE^{(1)}, \widetilde{\sF})$ and $(\sX,d',\mu,\tilde \sE^{(2)}, \widetilde{\sF})$.
If $r > 2$ then $B'(x,r) = B'(x,r/2) = \sX$, and so if $h$ is positive and $\sE^{(1)}$-harmonic on $B'(x,r)$ then
$h$ is constant, and thus the EHI holds for the MMD space $ (\sX,d',\mu,\tilde \sE^{(1)},\widetilde{\sF})$. 
 On the other hand, as the LP fails for $\sE^{(2)}$ there exists a non-constant
bounded $\sE^{(2)}$-harmonic function $h$ on $\sX$, which we can normalize so that 
 $\inf_{\sX}  h =0$, $\sup_{\sX} h=1$. 
It is thus clear that
the EHI for the MMD space $(\sX,d',\mu,\tilde \sE^{(2)}, \widetilde{\sF})$ fails for every ball $B'(x,\delta r) \subset B'(x,r)$  with $\delta \in (0,1)$ and $r > 1/\delta$. 
  } 
  \end{example}

\begin{remark}
{\rm  
It would be interesting to have an example 
 of a strongly local MMD space that does not have (MD) property for which  EHI fails to be stable   
  but for which all balls  are relatively compact. 
It does not seem easy to modify the example above to give this.
} \end{remark}

\begin{example}\label{E:8.3} 
 \rm	
We give an example of a strongly local irreducible MMD space where harmonic functions may be discontinuous and 
(Ha) fails. However the condition (HC) does hold. 
The space consists of three parts: the closure of a domain in $\bR^2$, the standard
Sierpinski gasket, and a line segment. Let
$\sX_1$ be the compact Sierpinski gasket, with vertices $z_1=(0,0)$, 
$z_2=(1,0)$ and $z_3=(\half, \frac{\sqrt{3}}2)$,  
  $\sX_2 = [0, 1]\times    ( -1, 0]$, a unit square with the bottom $[0, 1]\times \{-1\}$ removed, 
 and let $\sX_3$ be a smooth curve outside $\sX_1\cup\sX_2$ 
that connects the vertex $z_3$ 
of the Sierpinski gasket with the point  $z_4=(1, -1/2)$ 
at the middle of the right side of the square $\sX_2$.
We identify $\sX_3$ with a closed line segment of length $l >1$. 

Let $\sX = \sX_1 \cup \sX_2 \cup \sX_3$.
  For $x,y \in \sX$ let $d(x,y)$ be the (Euclidean) length of the shortest curve in $\sX$ connecting
$x$ and $y$. 
Clearly, $(\sX, d)$ is a   locally compact separable   metric space. 
Let $m_1$ be the measure  on $\sX_1$ which assigns mass $3^{-n}$ to each triangle of
side $2^{-n}$, and for $j=2,3$,  let $m_j$ be Lebesgue measure on $\sX_j$. 
Let $m$ be the measure on $\sX$ such that $m|_{\sX_i} = m_i$ for each $i$.
Clearly, $m$ is a finite measure on $\sX$. 
 
Let  $(\sE^{(1)}, \sF^{(1)})$ be the strongly local Dirichlet form on $L^2(\sX_1,m_1)$
associated with the standard diffusion on the Sierpinski gasket -- see \cite[Chapter 3]{Kig}. 
 Let $(\sE^{(2)}, \sF^{(2)})$ be the Dirichlet form associated with the  reflected Brownian motion
 on $\sX_2$    killed upon hitting  the horizontal line segment $[0, 1]\times \{-1\}$. 
  Let $(\sE^{(3)}, \sF^{(3)})$ be the Dirichlet form associated with Brownian motion
on $\sX_3$, with reflection at the two endpoints.
  Since  $(\sE^{(1)}, \sF^{(1)})$  and $(\sE^{(3)}, \sF^{(3)})$  are resistance forms on $\sX_1$ and $\sX_3$,
every $\sE^{(i)}$-quasi-continuous $m_i$-version  of elements in $\sF^{(i)} $ is H\"older continuous   on $\sX_i$ for $i=1, 3$.
 Following \cite{Kum}
we can construct a strongly local regular Dirichlet form $(\sE,\sF)$ on $L^2(\sX,m)$ such that 
$\{ f|_{\sX_i} : f \in \sF \} = \sF^{(i)}$ for {\red $i=1,2, 3.$ } 
  
Let $X=\{X_t, t\geq 0; \bP^x, x\in \sX\}$ be the diffusion process associated with the regular Dirichlet
form $(\sE, \sF)$   on $L^2(\sX; m)$.  
The diffusion $X$ is transient   due to the killing on the line $\{x_2 =-1\}$.
This process   behaves as follows: 
\begin{description}
\item{(i)}  when $X_t$ is inside $\sX_1$, it behaves like Brownian motion on the Sierpinski gasket
$\sX_1$ until it reaches the vertex $z_3$ or the bottom $\sX_1\cap \sX_2$;

\item{(ii)} when $X_t$ is inside $\sX_2,$ it behaves like two-dimensional Brownian motion in $\sX_2$, 
 with reflection
on the two lines $\sX_2 \cap \{ x_1=0\}$ and $\sX_2 \cap \{ x_1=1\}$, and killing on $\sX_2 \cap \{ x_2=-1\}$.  
\item{(iv)} when $X_t$ is at the vertex $z_3$, it immediately enters both $\sX_1\setminus \sX_3$ and $\sX_3\setminus \sX_1$.
When $X_t$ is at a point $y \in \sX_1 \cap \sX_2$, it 
immediately enters both $\sX_1\setminus \sX_2$ and $\sX_2\setminus \sX_1$.
When  $X_t$ is at $z_4$, it is reflected into   $\sX_3$. 
\end{description} 

The unique point  $z_4 \in \sX_2 \cap \sX_3$  
is polar for reflected Brownian motion in $\sX_2$, so 
the process  $X$ starting from $\sX_2 \setminus \{z_4\}$ can only enter $\sX_3$ through 
the Sierpinski gasket $\sX_1$ via vertex $z_3$. 
 
For any $r\in (0, 1/2)$,    let $u\in \sF$ be a bounded function so that 
$u=1$ on  $B(z_4, r/2)^c \cap \sX_2$ and $u=0$ on $B(z_4, r/2)^c \cap \sX_3$.
Define  $h(x)= \bE^x [ u(X_{\tau_{B(z_4, r)}}]$, which is  a bounded function in $\sF_e$
and is regular harmonic in $B(z_4, r)$. Observe that  
   $h(x)=   \bP^x (X_{\tau_{B(z_4, r)}}\in \sX_2   )$ for $x\in B(z_4, r)$ 
  with     $h(x)=1 $ for $x\in  B(z_4, r) \cap \sX_2    \setminus \{z_4\}  $
and $h(x)=0$ on $B(z_4, r)\cap \sX_3  $.  
 Thus $h$ does not satisfy the non-scale-invariant
Harnack inequality.  In other words,  (Ha) fails for this strongly local Dirichlet form $(\sE, \sF)$.  
Let 
$$ 
  D_k  = \left\{ x =(x_1, x_2)\in \sX_2: 0< |x-z_4| < 8^{-k} \ \hbox{ or }\   x_2<  8^{-k}  -1   \right\}   
$$
 and set $F_k = \sX \setminus D_k$. 
  It is straightforward to verify that 
$(F_k)$ is  an $\sE$-nest consisting of 
 compact sets with $\cup_{k=1}^\infty F_k = \sX$.

  Denote by $\zeta$ the lifetime of $X$. It is not hard to see that $\bE^x \zeta$ is bounded on 
$\sX$. So there is a small constant  $c_0 >0$ so that for  $g_0=c_0$, $Gg_0(x)= c_0 \, \bE^x \zeta \leq 1$ on $\sX$, $Gg_0\in \sF_e$ with
$\sE (Gg_0, Gg_0)=\int_{\sX} g_0 (x)Gg_0 (x) m(dx)\leq 1$.    
  We now verify (HC).  We begin by specifying the function $r_{x_0}$. 
If $x_0$ is in exactly one of the sets $\sX_i$, $i=1,2,3$, then 
choose $r_{x_0}$ small enough so that $B(x_0, 2r_{x_0})$ is contained in that  $\sX_i$. 
The remaining possibilities are that $x_0 \in \sX_1 \cap \sX_2$, $x_0= z_3$ or $x_0 =z_4$, and in these
cases we take $r_{x_0 }= \fract17$.
Now let $r < r_{x_0}$, and let $f  \in \mathcal{B}_+(\sX)$ 
have compact support in $B(x_0, 2r)^c$, 	and satisfy
$0\le f \leq cg_0$ for some $c>0$.
We have
$$		Gf(x)= \bE^x[ Gf (X_{  \tau_{B(x_0, 2r)}})] \quad \hbox{for } x\in B(x_0, r).	$$

Now let $x_0 \in \sX$. If for some $i$ we have $B(x_0, 2r_{x_0}) \subset \sX_i$, 
 then the continuity of $Gf$ on $B(x_0,r)$ follows exactly as for the space $\sX_i$.
 Next suppose that
$x_0 \in \sX_1 \cap \sX_2$. 
Since $Gf \in \sF$  and is $\sE$-quasi-continuous on $\sX$,   $Gf|_{\sX_1} \in \sF^{(1)}$  and is $\sE^{(1)}$-quasi-continuous. 
So  $Gf$ is H\"older continuous on $\sX_1 \cap B(x_0, r)$ and, in particular, is
continuous on the the line segment $L= \sX_1 \cap \sX_2 \cap B(x_0, 2r)$. So
$Gf$ is harmonic in $B(x_0, 2r) \setminus \sX_1$ and is continuous on $L$.
 For $x\in B(x_0, r) \cap \sX_2$, by the strong Markov property of $X$, 
$$ Gf (x)= \bE^x [ Gf (X_{\tau_{B(x_0, 2r)\cap \sX_2}}  )]= \bE^x [ Gf (Y_{\tau_{B(x_0, 2r)\cap \sX_2}}  )],
$$
where $Y$ is the normally reflected Brownian motion in the rectangle $[0, 1]\times [-1, 1]$. 
Note that $Y$ is a Feller process having strong Feller property and has two-sided short time 
Gaussian type heat kernel estimates on $[0, 1]\times [-1, 1]$, and the open set $(0, 1]\times [-1, 0)$
satisfies the exterior cone condition for every boundary point on $(0, 1]\times \{0\}$.
Thus every point on $(0, 1]\times \{0\}$ is regular for $[0,1]\times [0, 1]$ with respect to the reflected Brownian motion $Y$.  
 It follows that  $Gf$ is continuous in $B(x_0, r) \cap \sX_2$, and therefore in $B(x_0,r)$
  
If $x_0 =z_3$ then $Gf$ is continuous on
$\sX_1 \cup \sX_3$, and so is continuous on $B(x_0,r)$. 
Finally let $x_0=z_4$.  
Then $B(z_4, 2r)\cap F_k$ consists of two disjoint components, a part annulus $A_1$ contained in $\sX_2$
and the set $A_2=\{ y \in \sX_3: 0\le d(z_4,y) < 2r\}$. As above, we have that $G f$ is continuous on each
component $A_i$ and so is continuous on $B(z_4, r) \cap F_k$.
Thus $\sX$ satisfies (HC). 

Let
$$ h(x) = \bP^x( T_{z_3} < \infty \}. $$
Then since $X$ can only leave $\sX_3$ via the point $x_3$, we have 
$h|_{\sX_3} \equiv 1$. On the other hand  the symmetry of $\sX_2$ implies that 
$\bP^x( T_{\sX_1} < \infty ) \simeq \half$ for points $x \in \sX_2\setminus \{ z_4\}$ close to $z_4$,
and thus $h$  is not continuous at $z_4$.

Note further that  the point $z_4$ is of positive capacity and $(\sE, \sF)$ is irreducible.   
On the other hand, the part Dirichlet form $(\sE, \sF^{B(z_4, r)})$ on $L^2( B(z_4, r), m|_{B(z_4, r)})$
is not irreducible for any $r\in (0, 1/2]$; the space $B(z_4, r)$,   which is connected, 
 has two disjoint invariant sets
$B(z_4, r)\cap (\sX_2\setminus \{z_4\})$ and $B(z_4, r) \cap \sX_3$   for the part process $X^{B(z_4, r)}$
of $X$ killed upon leaving $B(z_4, r)$. 
Thus this example also shows that a strongly local regular Dirichlet form
does not need to be irreducible even though the underlying metric space is connected.
   \end{example}

 \begin{example}\label{E:8.cylinder}  \rm	 
 We now show that the spaces  studied in \cite{BSC} provide an example of a MMD space which fails 
 the condition (HC). Let $r_k \in (0,1)$  be chosen so that $(r_k)$ is strictly decreasing and 
 $$ \sum_1^\infty  e^{- t /r_k^2 } = +\infty \q \hbox { for all } t \ge 0.
 $$
  Let $S_k$ be a circle of radius $r_k$, and $m_k$ be Lebesgue measure on $S_k$ normalized so that
 $m_k(S_k)=1$.  
 Set $\sX = \prod_{k=1}^\infty S_k$, and $m = \otimes_1^\infty m_k$.
Since each $S_k$ is compact, the space $\sX$ is compact
and therefore locally compact. 
Let $d_k$ be the usual metric on $S_k$, and define a metric $d$ on $\sX$
by taking $d(x,y) = \sup_k \{ d_k(x_k, y_k) \}$: this metric induces the product topology on $\sX$.
For $k \in \bN$ let $W^{(k)}$ be independent Brownian motions on $S_k$ and set
$W=( W^{(1)},  W^{(2)}, \cdots ) \in \sX$. Then $W$ is a conservative symmetric Hunt process with invariant measure $m$.
  Its  associated  Dirichlet form $(\sE, \sF)$ is strongly local and regular  on $L^2(\sX,m)$.  
   For more details of this construction  see \cite{BSC}.
There exist heat kernel measures $h_t(x, \cdot)$ such that for any $f \in C(\sX)$,
$$ \bE^x f(W_t) = \int h_t(x,dy) f(y). $$
By Theorem 1.2 of \cite{BSC} the measure $h_t(x, \cdot)$ is singular with respect to $m$ for all $x \in \sX$ and
$t \ge 0$; thus $W$ does not have a regular transition density withy respect to $m$. 
 The process $W$ is recurrent, but we can define a transient process by choosing a point
$z_1 \in S_1$ and killing $W$ when $W^{(1)}$ hits $z_1$. Write $\ol W$ for this killed process, and 
 $(\ol \sE, \ol \sF)$ for the associated Dirichlet  form. 
 If $\ol h_t(x, \cdot) = \bP^x( \ol W \in \cdot)$ then $\ol h_t \le h_t$ and thus is also singular with respect to $m$.
 The condition (HC) must fail for  $(\ol \sE, \ol \sF)$, since otherwise by Theorem \ref{T:2}
$\ol W$ would have a transition density with respect to $m$.
We note that this space has infinite Assoud dimension and therefore fails to satisfy metric doubling.
   \end{example}

  \begin{example}\label{E:8.5} \rm	
To give a concrete example of an irreducible strongly local MMD space that fits the setting of Theorem  \ref{T:main-stable} but  
fails to satisfy the local regularity of \cite[Assumption 2.5(i)]{BM1} in the compact setting, 
 consider $\sX$ to be the join of the Vicsek tree (compact) with the unit interval $[0,1]$, where the symmetrizing measure  $m$
is given by the Hausdorff measure on each of the pieces. 
The space $\sX$ satisfies the relatively ball connected condition. 
We take $(\sE,\sF)$ to be the strongly local regular Dirichlet form on $L^2(\sX; m)$ obtained by combining the 
Dirichlet form associated
with Brownian motion on $(0,1]$ with the Dirichlet form associated with the diffusion on the Vicsek tree, in a similar
fashion to the previous example. The argument in \cite{De2} 
 can be adapted to show  
that this example satisfies the EHI.
This example is essentially due to Delmotte \cite{De2}.
 \end{example}

	\sm {\bf Acknowledgment. }  
	We thank Naotaka Kajino for helpful comments on an earlier draft of this work, and for pointing us to
	the spaces studied in Example \ref{E:8.5}.
  We are grateful to the anonymous referee for a very careful reading and detailed comments, which helped improving the exposition of this paper. 
	The research of Martin Barlow is partially supported by an NSERC grant.
The research of Zhen-Qing Chen is partially supported by a Simons Foundation Grant. 
The research of Mathav Murugan  is partially supported  by the Canada Research Chairs program and an NSERC grant.

	\bigskip
	 {\bf Martin T. Barlow}

	   Department of Mathematics,
	University of British Columbia,
	Vancouver, BC V6T 1Z2, Canada. 
	
	E-mail:  barlow@math.ubc.ca

\bigskip	
	
	{\bf Zhen-Qing Chen}

Department of Mathematics, University of Washington, Seattle,
WA 98195, USA.

E-mail: zqchen@uw.edu

\bigskip

{\bf Mathav Murugan}

Department of Mathematics,
	University of British Columbia,
	Vancouver, BC V6T 1Z2, Canada.  

E-mail: mathav@math.ubc.ca

\end{document}